\documentclass[11pt]{article}
\usepackage{amsmath,amsfonts,amssymb,algorithm,algorithmic,psfrag}
\usepackage{amsthm}
\usepackage{array}
\usepackage{psfrag}
\usepackage{bm}
\usepackage{color}
\usepackage{bigstrut}

\usepackage{listings}
\newcommand{\R} {\mathbb R}
\newcommand{\calP}{\mathcal{P}}
\newcommand{\calQ}{\mathcal{Q}}
\newcommand{\calS}{\mathcal{S}}

\newcommand{\raw}{\rightarrow}

\newcommand{\ds}{\displaystyle}

\usepackage{color} 
\usepackage{ifpdf}
\ifpdf
    \usepackage[pdftex]{graphicx}
    \usepackage[pdftex]{hyperref}
    \hypersetup{
        unicode=false,          
        pdftoolbar=true,        
        pdfmenubar=true,        
        pdffitwindow=false,     
        pdfstartview={FitH},    
        pdftitle={AKG Article},      
        pdfauthor={Andrew Gillette},   
        pdfsubject={Mathematics},    
        pdfcreator={Andrew Gillette},  
        pdfproducer={Andrew Gillette}, 
        pdfkeywords={math, mathematics}, 
        pdfnewwindow=true,      
        colorlinks=true,        
        linkcolor=red,          
        citecolor=blue,         
        filecolor=magenta,      
        urlcolor=cyan           
    }

\else
    \usepackage{graphicx}
    \usepackage{pstricks}
    
    \newcommand{\href}[2]{#2}
\fi

\title{Numerical studies of serendipity and tensor product elements for eigenvalue problems}

\author{
Andrew Gillette, 
Craig Gross, 
and Ken Plackowski
\footnote{
Department of Mathematics,
University of Arizona,
Tucson, AZ 85721. \newline
\text{}\qquad{\it agillette@math.arizona.edu,~~~cgross95@email.arizona.edu,~~~plackow1@math.arizona.edu}}
}

\begin{document}
\maketitle

\begin{abstract}
While the use of finite element methods for the numerical approximation of eigenvalues is a well-studied problem, the use of serendipity elements for this purpose has received little attention in the literature.  
We show by numerical experiments that serendipity elements, which are defined on a square reference geometry, can attain the same order of accuracy as their tensor product counterparts while using dramatically fewer degrees of freedom.
In some cases, the serendipity method uses only 50\% as many basis functions as the tensor product method while still producing the same numerical approximation of an eigenvalue.
To encourage the further use and study of serendipity elements, we provide a table of serendipity basis functions for low order cases and a Mathematica file that can be used to generate the basis functions for higher order cases.
\end{abstract}

\section{Introduction}

Computational approximation of eigenvalues is a topic of ongoing interest across a broad spectrum of the applied mathematics community, due in part to the wide variety of application areas where it is required.
In this work, we compare two finite element methods for the computation of eigenvalues of the Laplacian: tensor product and serendipity.
While tensor product finite element methods have been used for decades to compute eigenvalues, the lesser known serendipity elements have been employed rarely, if ever, in this context, despite the fact that they are expected to require fewer computations to achieve the same order of accuracy.
%
%
%

The potential benefits of a serendipity element eigenvalue solver are obvious from a rough estimate of the degrees of freedom required for a method with $\mathcal{O}(h^p)$ error decay.
Here, $h$ indicates the maximum diameter of an affinely-mapped square mesh element and $p\geq 1$ indicates the maximum exponent of any variable appearing in a basis for the element.
The tensor-product finite element method for $H^1$-conforming problems in $\R^n$ uses $(p+1)^n$ basis functions per element, while the serendipity method uses roughly $p^n / n!$ for large~$p$.
Thus, for domains in $\R^2$, an $\mathcal{O}(h^p)$ serendipity method has about 50\% the number of basis functions as its tensor product counterpart, while for domains in $\R^3$, an $\mathcal{O}(h^p)$  serendipity method has only 17\% the number of basis functions as a tensor product method. 
As we show by numerical evidence, these computational savings are not restricted to an asymptotic regime but can be realized even in domains in $\R^2$ and for values $p\leq 6$.

The body of prior work studying finite element methods for eigenvalue approximation dates back to the 1970s~\cite{H1979} and is quite large, due to the many options available when designing finite element schemes and the many kinds of inquiries that could be made.
An excellent survey of the research in this area was given by Boffi in~\cite{B2010}.
While many works are concerned with approximation of the spectrum of the Laplacian (e.g.\ concerns about pollution and completeness of the computed spectrum), here we focus on accurate computation of individual eigenvalues to a high order of accuracy with the goal of minimizing the number of global degrees of freedom.
A similar kind of study by Wang, Monk, and Szabo compared $h-$ and $p-$refinement schemes on tetrahedra for computing resonant modes in a cavity using tetrahedral elements~\cite{WMS1996}.
This work focuses on square elements, which offer greater ability to reduce the number of global degrees of freedom than simplicial elements.

In this paper, we carry out a series of numerical experiments to compare the accuracy of serendipity and tensor product finite element methods in the context of eigenvalue computation.
We compare square and $L$-shaped domains, Dirichlet and Neumann boundary conditions, and $h$- and $p$-refinement strategies.
To ensure a fair comparison, we implement basis functions for both tensor product and serendipity elements using the construction process described in the work of Floater and Gillette~\cite{FG2014}, which uses interpolation conditions based on partial derivative data at edge and cell midpoints.
To the best of our knowledge, this is the first time such functions have been tested numerically.

Our results show that a $p$-refinement strategy with serendipity elements is preferable to the same strategy with tensor product elements in a variety of domain and boundary condition scenarios.
In particular, we find many specific instances where the serendipity elements achieve the same order of accuracy as the corresponding tensor product element with only 50\% the number of degrees of freedom.
The results also show that an $h$-refinement strategy does not always favor serendipity elements, meaning application context is essential when deciding between the use of tensor product and serendipity elements.

The remainder of the paper is organized as follows.
In Section~\ref{sec:fem-eval}, we review the eigenvalue problem for the Laplace equation with Neumann and Dirichlet boundary conditions, as well as the derivation of a Galerkin finite element method. 
Following this is a discussion of the two families of finite elements studied in this paper: tensor product and serendipity. 
In Section~\ref{sec:implement}, we state interpolation conditions that involve both values and derivative values and compute the basis functions for both the tensor product and the serendipity finite elements. 
We also discuss the relevant components for implementation via Mathematica and MATLAB.
In Section~\ref{sec:num-exp}, we provide a description of our results and a discussion of the comparison between the tensor product and serendipity elements. 
This includes comparisons of the aforementioned scenarios.
In Section~\ref{sec:conc}, we summarize our conclusions and give some directions for future research.
Finally, in Appendix~\ref{app:mathematica}, we give tables of the serendipity basis functions that we use and provide a link to a Mathematica code that can be used for further studies.

\section{Finite element methods for eigenvalue problems}
\label{sec:fem-eval}

Our focus in this work is the scalar-valued Laplace eigenvalue problem.
With Dirichlet boundary conditions, the problem is to find $\lambda\in\R$ and $u\in H^2(\Omega)$ such that:
\begin{equation}
\label{eq:lep-d-strong}
\left\{
\begin{array}{rcll}
		-\Delta u &=& \lambda u &\textrm{in } \Omega,\\
		       u &=& 0 &\textrm{on } \partial \Omega.
		       \end{array}
\right.
\end{equation}
With Neumann boundary conditions, the problem is to find $\lambda\in\R$ and $u\in H^2(\Omega)$ such that:
\begin{equation}
\label{eq:lep-n-strong}
\left\{
\begin{array}{rcll}
		-\Delta u &=& \lambda u &\textrm{in } \Omega,\\
		       \ds\frac{du}{d\textbf{n}} &=& 0 &\textrm{on } \partial \Omega.
		       \end{array}
\right.
\end{equation}
where $\textbf{n}$ is the unit vector normal to the boundary of $\Omega$.

We consider two subsets of $\R^2$ for the domain $\Omega$: the unit square $[0,1]^2$ and the L-shaped domain, $[0,2]^2 - (1, 2]^2$.
On $[0,1]^2$, the eigenvalues for the Dirichlet problem (\ref{eq:lep-d-strong}) are
\[(m^2 + n^2) \pi^2 \quad \textrm{ for }~ m, n \in\{1, 2, \ldots\} \textrm{ and } \Omega = [0, 1]^2.\]
For the Neumann problem (\ref{eq:lep-n-strong}) on $[0,1]^2$, the eigenvalues are 
\[(m^2 + n^2) \pi^2 \quad \textrm{ for }~ m, n \in\{0,1, 2, \ldots\} \textrm{ and } \Omega = [0, 1]^2,\]
the only difference being that $m$ and $n$ are allowed to have value 0.
For $(m,n)$ pairs with $m\not=n$, the corresponding eigenvalue has multiplicity at least 2 while those with $m=n$ have multiplicity 1 and are called `simple.'

On the L-shaped domain, Dauge~\cite{Dauge2003} has given benchmark computations with at least 8 digits of accuracy for the lowest non-zero eigenvalues for the Neumann problem.
The first four of these are
\begin{align*}
\lambda^{(1)} = 1.4756218450 & \qquad\lambda^{(2)} = 3.5340313683 \\
\lambda^{(3)}= 9.8696044011 & \qquad\lambda^{(4)}= 11.389479398
\end{align*}
Note that $\lambda^{(3)}=2\pi^2$, which is also an eigenvalue for the Dirchlet problem.
In our experiments, we look at approximation of $2\pi^2$ for each kind of boundary condition as well as the approximation of $\lambda^{(1)}$ for the Neumann case.

Discretization of (\ref{eq:lep-d-strong}) for numerical approximation begins with the weak form of (\ref{eq:lep-d-strong}).
Set $V:=H^1_0(\Omega)$ and find $\lambda\in\R$ and $u\in V$, $u\not=0$ such that
\begin{equation}
\label{eq:lep-weak}
\int_\Omega\nabla u\cdot\nabla v = \lambda \int_\Omega u~v\qquad\forall v\in V.
\end{equation}
A Galerkin finite element method seeks a solution to (\ref{eq:lep-weak}) that holds over a finite-dimensional subspace $V_{h,p}\subset V$: find $\lambda_{h,p}\in\R$ and $u_{h,p}\in V_{h,p}$, $u_{h,p}\not=0$ such that
\begin{equation}
\label{eq:lep-Galerk}
\int_\Omega\nabla u_{h,p}\cdot\nabla v_{h,p} = \lambda_{h,p} \int_\Omega u_{h,p}~v_{h,p}\qquad\forall v_{h,p}\in V_{h,p}.
\end{equation}
The dimension of $V_{h,p}$ is determined by the type of element used (tensor product or serendipity, in our case) in addition to the parameters $h$ and $p$.
Here, $h$ indicates the maximum diameter of a element in the mesh and $p$ indicates the maximum exponent of any variable appearing in the monomial basis for the element.
Hence, as $h\raw 0$ or $p\raw\infty$, we have $\dim V_{h,p}\raw\infty$.

We consider two possible choices for $V_{h,p}$ that are subsets of $H^1(\Omega)$ and associated to a partition of $\Omega$ into a mesh of squares.
We will follow notational conventions from the Periodic Table of the Finite Elements~\cite{femtable,AL2014} to describe the two choices in terms of the local spaces on each square element.
The first choice for a local space is $\calQ_p^-\Lambda^0(\square_2)$, more commonly known as the tensor product element of order $p$ on a square~\cite{ABB2012}.  
This element has 1 degree of freedom per vertex, $(p-1)$ degrees of freedom per edge, and $(p-1)^2$ degrees of freedom associated to the interior for a total of $(p+1)^2$ degrees of freedom per square element.
The second choice for a local space is $\calS_p\Lambda^0(\square_2)$, known as the serendipity element of order $p$ on a square~\cite{AA2011}.
The serendipity element has the same degrees of freedom associated to vertices and edges of the square, but only $\frac 12 (p-3)(p-2)$ degrees of freedom\footnote{For $p=1$, there are no interior degrees of freedom; the formula applies for any $p\geq 2$.} associated to the interior of the square.
It has a total of $\frac 12\left(p^2+3p+6\right)$ degrees of freedom per element.

In addition to the type of domain $\Omega$ (square or L-shaped), the family of element ($\calQ^-$ or $\calS$), and the order of $p$ selected, the dimension of $V_{h,p}$ depends on the maximum diameter of a mesh element.
We only consider meshes where all elements are squares of the same side length $h$, so that the maximum diameter of a mesh element is $\sqrt 2 h$.
By this convention, if $h=1/N$ for an integer $N\geq 1$, the square domain will have $N^2$ elements and the L-shaped domain will have $3N^2$ elements.
By counting the total number of vertices, edges, and elements in the mesh, we have the formula
\[\dim V_{h,p}= \text{ (\# vertices) + $(p-1)\cdot$(\# edges) + $\left(\frac{\text{\# DoFs}}{\text{interior}}\right)\cdot$(\# elements)},\]
where the number of degrees of freedom (DoFs) per interior depends on the choice of $\calQ^-$ or $\calS$, as described above.
Note that when Dirichlet boundary conditions are used, the value of degrees of freedom associated to the boundary of the domain are set to zero, which decreases the dimension of $V_{h,p}$.

The goal of the numerical experiments in this paper is to study the following question: given a domain, a set of boundary conditions, a rough guess for an eigenvalue $\lambda$, an $h$-refinement or $p$-refinement strategy, and a desire to attain a precise estimate of $\lambda$ while avoiding fruitless growth in $\dim V_{h,p}$, is it better to use $\calQ^-$ or $\calS$ elements?
Since the $\calS_p\Lambda^0$ and $\calQ^-_p\Lambda^0$ element each contain polynomials of total degree at most $p$ and $\dim\calS_p\Lambda^0 < \dim \calQ^-_p\Lambda^0$ for $p\geq 2$, we might expect that the serendipity elements would be preferable in every case.
On the other hand, perhaps the ``extra'' approximation power afforded by the larger basis in the tensor product element provides better eigenvalue estimation overall.
To make a fair comparison, we implement serendipity and tensor product elements by the same methodology, and then report their results when used in a series of computational experiments.

\section{Implementation of Serendipity Elements}
\label{sec:implement}

Here, for the first time, we compute and employ the basis functions for $\calS_p\Lambda^0(\square_2)$ with Hermite-like interpolation conditions at edge midpoints, as described in the recent work of Floater and Gillette~\cite{FG2014}.
We review the degrees of freedom for these elements here and explain how the process outlined in~\cite{FG2014} was used to derive the basis functions employed in our numerical experiments.

\paragraph*{Serendipity degrees of freedom}
The term ``serendipity element'' has appeared in various mathematical and engineering texts since the 1970s~\cite{BS2002,Ci02,H1987,M1990,SB1991,SF73}, referring to the fact that these elements seemed to achieve $O(h^p)$ accuracy with fewer degrees of freedom than their tensor product counterparts. 
Arnold and Awanou provided degrees of freedom in the classical finite element sense for the $H^1$-conforming version of these spaces in~\cite{AA2011}: for a $d$-dimensional face $\square_d$ of of an $n$-cube $\square_n$, the order $p$ serendipity degrees of freedom for a scalar function $u$ are
\begin{equation}
\label{eq:srdp-dofs}
u\longmapsto\int_{\square_d} uq,\qquad \forall~q\in\calP_{p-2d}(\square_d),
\end{equation}
where $\calP_{p-2d}(\square_d)$ denotes the space of polynomials in $n$-variables of degree $\leq p-2d$ on face $\square_d$.
For $n=2$ and $p\geq 2$, the space of polynomials associated to the degrees of freedom (\ref{eq:srdp-dofs}) is denoted $\calS_p\Lambda^0(\square_2)$ and given by
\begin{equation}
\label{eq:srdp-monoms}
\calS_p\Lambda^0(\square_2) = \calP_{p}(\square_2)\oplus\text{span}\{x^py, xy^p\}.
\end{equation}
It is shown in~\cite{AA2011} that the degrees of freedom (\ref{eq:srdp-dofs}) are unisolvent for (\ref{eq:srdp-monoms}), but a consideration of how to construct suitable basis functions for the implementation of these elements in applications was not provided.

\paragraph*{Basis functions for serendipity elements}
We use a procedure outlined by Floater and Gillette in~\cite{FG2014} to construct basis functions for the $\calS_p\Lambda^0(\square_2)$ element.
To the best of our knowledge, these functions have not been constructed explicitly or used in numerical experiments previously.
The procedure is also used to construct bases for the $\calQ_p^-\Lambda^0(\square_2)$ element.

Given $p\geq 1$, we first we define a set of $p+1$ functions over $[-1,1]$, denoted
\[ \Phi_p[x] := \left\{\phi_1(x), \dots, \phi_{p+1}(x)\right\}.\]
Let $D$ denote the endpoints and midpoint of $[-1,1]$, i.e.\ $D=\{-1,0,1\}$, and denote the Kronecker delta function by
\[\delta_i(j)=\begin{cases} 0 & \text{if $i\not=j$,} \\ 1 & \text{if $i=j$.}\end{cases}\]
Define $\Phi_1[x]:=\{(1-x)/2,(1+x)/2\}$.
For $p\geq 2$, fix the following interpolation properties\footnote{If a set of indices on the right is empty, the property should be treated as vacuous.}:
\begin{alignat}{3}
 \phi_1(x_0) & = \delta_{-1}(x_0),\quad && \forall\,x_0\in D \label{eq:interp-cond-first}\\
 \phi_2(x_0) & = \delta_{0}(x_0), \quad && \forall\,x_0\in D \\
 \phi_{p+1}(x_0) & = \delta_{1}(x_0), \quad && \forall\,x_0\in D \\
  \phi_i^{(k)}(0) & = 0, \quad && \forall\, i\in\{1,2,p+1\},~\forall\,k\in\{1,\ldots, p-2\} \\
 \phi_i(x_0) & = 0, \quad && \forall\,x_0\in D,~\forall\,i\in\{3,\ldots,p\} \\
 \phi_i^{(i-2)}(0) & = 1, \quad && \forall\,i\in\{3,\ldots,p\} \\
\phi_i^{(k)}(0) & = 0, \quad && \forall\, i\in\{3,\ldots,p\},~\forall\,k\in\{1,\ldots, i-3\} \label{eq:interp-cond-last}
\end{alignat}
For $i=1$ to $p+1$, we find the lowest degree polynomial $\phi_i$ that satisfies the above constraints.
Since there are at most $p+1$ constraints for each $i$, this process uniquely defines a set of $p+1$ polynomials, each of degree at most $p$.
Moreover, $\phi_1$, $\phi_2$, and $\phi_{p+1}$ are the only functions in the set that have non-zero value at -1, 0, and 1, respectively, while the functions $\phi_3$ through $\phi_p$ have linearly independent constraints on their derivatives at 0.
Thus, for each $p\geq 1$, $\Phi_p[x]$ is a basis for $\calP_p([-1,1])$.
The sets $\Phi_1[x],\ldots,\Phi_5[x]$ are listed explicitly in Table~\ref{tab:Phi-bases}.

\begin{table}[h]
\begin{tabular}{c|ccccc}
$p$ & $\phi_1(x)$ & $\phi_2(x)$ & $\phi_3(x)$ & $\phi_4(x)$ & $\phi_5(x),~\phi_6(x)$ \\ \hline
1 & $-\frac 12(x-1)$ & $\frac 12 (1+x)$ \\ \hline
2 & $\frac{1}{2} (x-1) x$ & $1-x^2$ & $\frac{1}{2} x (x+1)$ \\ \hline
3 & $-\frac{1}{2} (x-1) x^2$ & $1-x^2$ & $x-x^3$ & $\frac{1}{2} x^2 (x+1)$ \\ \hline
4 & $\frac{1}{2} (x-1) x^3$ & $1-x^4$ & $x-x^3$ & $-\frac{1}{2} (x-1) x^2 (x+1)$ & $\frac{1}{2} x^3 (x+1)$ \\ \hline
5 & $-\frac{1}{2}(x-1) x^4$ & $1-x^4$ & $x-x^5$ & $-\frac{1}{2} (x-1) x^2 (x+1)$ & $-\frac{1}{6} (x-1) x^3 (x+1),$\\
 &&&&& $\frac{1}{2} x^4 (x+1)$ \\
\end{tabular}
\caption{Basis functions for $\Phi_p[x]$ with $1\leq p\leq 5$.}
\label{tab:Phi-bases}
\end{table}

By taking tensor products of the $\Phi_p[x]$ sets, we can build out bases for tensor product and serendipity spaces over $[-1,1]^n$ for any $n\geq 1$, although we we consider only $n=2$ here.
We fix the notation 
\[\Phi_{pq} :=\left\{\phi_i(x)\phi_j(y)~:~\phi_i(x)\in\Phi_p[x],~~\phi_j(y)\in\Phi_q[y]\right\},\]
where $p$ and $q$ need not be distinct.
A basis for the tensor product space $\calQ_p^-\Lambda^0(\square_2)$ can be computed immediately as
\[\text{basis for}~\calQ_p^-\Lambda^0(\square_2) = \Phi_{pp}.\]
A basis for the serendipity space $\calS_p\Lambda^0(\square_2)$ is more involved to describe but only slightly more difficult to compute.
First, an addition operation on sets of the type $\Phi_{pq}$ is defined as follows.
To build the set $\Phi_{pq}+\Phi_{rs}$, let $M=\max\{p,q,r,s\}$ and build a square array of indices $\{1,\ldots,M+1\}\times\{1,\ldots,M+1\}$.
Associate function $\phi_i(x)\phi_j(y)\in\Phi_{pq}$ to index $\{k,\ell\}$ according to the rule
\[\phi_i(x)\phi_j(y)\longmapsto
\begin{cases} 
\{M+1,j\} &~\text{if $i=p+1$, $j<q+1$}\\
\{i,M+1\} &~\text{if $i<p+1$, $j=q+1$}\\
\{M+1,M+1\} &~\text{if $i=p+1$, $j=q+1$}\\
\{i,j\} & ~\text{otherwise}
\end{cases}\]
Associate function $\phi_i(x)\phi_j(y)\in\Phi_{rs}$ to indices according to the same rule, replacing $p$ by $r$ and $q$ by $s$.
Initialize $\mathbb A_{pq,M}$ as an $(M+1)\times (M+1)$ array of zeroes, then place the functions from $\Phi_{pq}$ into $\mathbb A_{pq,M}$ according to their index assignment.
Define $\mathbb A_{rs,M}$ analogously, using functions from $\Phi_{rs}$.
The set $\Phi_{pq}+\Phi_{rs}$ is then defined to be the set of non-zero entries of $\mathbb A_{pq,M}+\mathbb A_{rs,M}$.
In practice, this re-indexing and summation procedure is carried out by insertion of rows or columns of zeros at appropriate places into the arrays storing $\Phi_{pq}$ and $\Phi_{rs}$ and then adding the arrays together.


A basis for $\calS_p\Lambda^0(\square_2)$ can then be written as a linear combination of this addition operation on some $\Phi_{rs}$ sets.
For $p=1$ through $p=6$, these linear combinations are as follows:

\begin{align}
\calS_1\Lambda^0(\square_2)~\text{basis}~&= \Phi_{1 1} \label{eq:srdp-LC-1}\\
\calS_2\Lambda^0(\square_2)~\text{basis}~&= \Phi_{2 1} + \Phi_{1 2} - \Phi_{1 1}\\
\calS_3\Lambda^0(\square_2)~\text{basis}~&= \Phi_{3 1} + \Phi_{1 3} - \Phi_{1 1}\\
\calS_4\Lambda^0(\square_2)~\text{basis}~&= \Phi_{4 1} + \Phi_{1 4} + \Phi_{2 2} - (\Phi_{2 1} + \Phi_{1 2})\\
\calS_5\Lambda^0(\square_2)~\text{basis}~&= \Phi_{5 1} + \Phi_{1 5} + \Phi_{3 2} + \Phi_{2 3} - (\Phi_{3 1} +\Phi_{1 3} + \Phi_{22})\\
\calS_6\Lambda^0(\square_2)~\text{basis}~&= \Phi_{6 1} + \Phi_{1 6} + \Phi_{4 2} + \Phi_{2 4} + \Phi_{3 3} \notag\\
 &\qquad\qquad - (\Phi_{4 1} +\Phi_{1 4} + \Phi_{2 3} + \Phi_{3 2})\label{eq:srdp-LC-6}
\end{align}
The derivation of these linear combinations is given in~\cite[Section 5]{FG2014} using different notation.
The techniques of \cite{FG2014} can produce bases in this way for $\calS_p\Lambda^0(\square_n)$ for any $p\geq 1$ and $n\geq 1$. As an example, in Appendix~\ref{app:mathematica}, we provide the 2D serendipity basis functions for $p=1$ to $4$.

\paragraph*{Implementation via Mathematica and \textsc{MATLAB}}
We use Mathematica to compute the bases for $\calQ_p^-\Lambda^0(\square_2)$ and $\calS_p\Lambda^0(\square_2)$ according to the procedure just described and the process of basis generation is summarized below.
The Mathematica function \texttt{InterpolatingPolynomial} is used to produce the sets $\Phi_p[x]$ based on the conditions (\ref{eq:interp-cond-first})-(\ref{eq:interp-cond-last}).
For example, $\phi_3(x)\in \Phi_3[x]$ should satisfy $\phi_3(-1)=\phi_3(0)=\phi_3(1)=0$ as well as $\phi_3'(0)=1$.
The unique cubic polynomial satisfying these constraints is computed by the command
\[
\texttt{InterpolatingPolynomial[\{\{-1, 0\}, \{0, 0, 1\}, \{1, 0\}\}, x]}
\]
We define a function \texttt{interpolatingList[p]} that creates the required inputs to \texttt{InterpolatingPolynomial} for each $\phi_i\in\Phi_p[x]$.
We also define a function \texttt{genTable2D[p,q,M]} that builds the array $\mathbb{A}_{pq,M}$.
Bases for $\calS_p\Lambda^0(\square_2)$ are constructed by simplifying linear combinations of appropriate \texttt{genTable2D[r,s,M]} arrays according to (\ref{eq:srdp-LC-1})-(\ref{eq:srdp-LC-6}); the value of $M$ is set to $p$ for each term in the combination so that the output is a $(p+1)\times(p+1)$ array with exactly $\dim\calS_p\Lambda^0(\square_2)$ non-zero entries.
The basis for $\calQ_p^-\Lambda^0(\square_2)$ is built by the command \texttt{genTable2D[p,p,p]}, which generates a $(p+1)\times(p+1)$ array with all entries non-zero.

Once the basis functions are created, we pass them to a finite element solver in \texttt{MATLAB} in order to compute approximate eigenvalues. 
The resulting finite element problem is given by
\[
\lambda M v = L v
\]
where $M$ is the mass matrix and $L$ is the stiffness matrix with
\begin{align*}
M &= \left[M_{i, j} \right] \textrm{ such that } M_{i, j} = \int_{\Omega} \psi_i \psi_j \; dA, \\
L &= \left[L_{i, j} \right] \textrm{ such that } L_{i, j} = \int_{\Omega} \langle \nabla \psi_i, \nabla \psi_j \rangle \; dA,
\end{align*}
where $\psi_i$, $\psi_j$ range over a basis for  $\calQ_p^-\Lambda^0(\square_2)$ or $\calS_p^-\Lambda^0(\square_2)$.
The finite element solver takes a local approach, making use of the basis functions defined over a reference element as above (specifically $[-1, 1]^2$). 
By calculating the desired entries of the mass and stiffness matrices over the reference element, scaling, and assembling on the global square or L-shaped domain, we produce global mass and stiffness matrices.

Furthermore, in the derivation of the variational form of the problem, the imposition of the Neumann conditions is encoded by the vanishing of any integrals over the boundary of the domain.
To impose the Dirichlet conditions, it is necessary to manipulate the equations in the discrete problem that solve for the coefficients corresponding to boundary nodes. 
Traditionally, this is realized by setting each of the coefficients corresponding to value interpolating nodes on the boundary equal to zero. 
As the tensor product and serendipity basis functions that we use include interpolation of some partial derivative values along the boundary, we also set to zero the coefficients of the basis functions corresponding to those conditions.

\section{Numerical Experiments and Results}
\label{sec:num-exp}

\begin{figure}[t]
$\begin{array}{ccc}
\includegraphics[width=.31\textwidth]{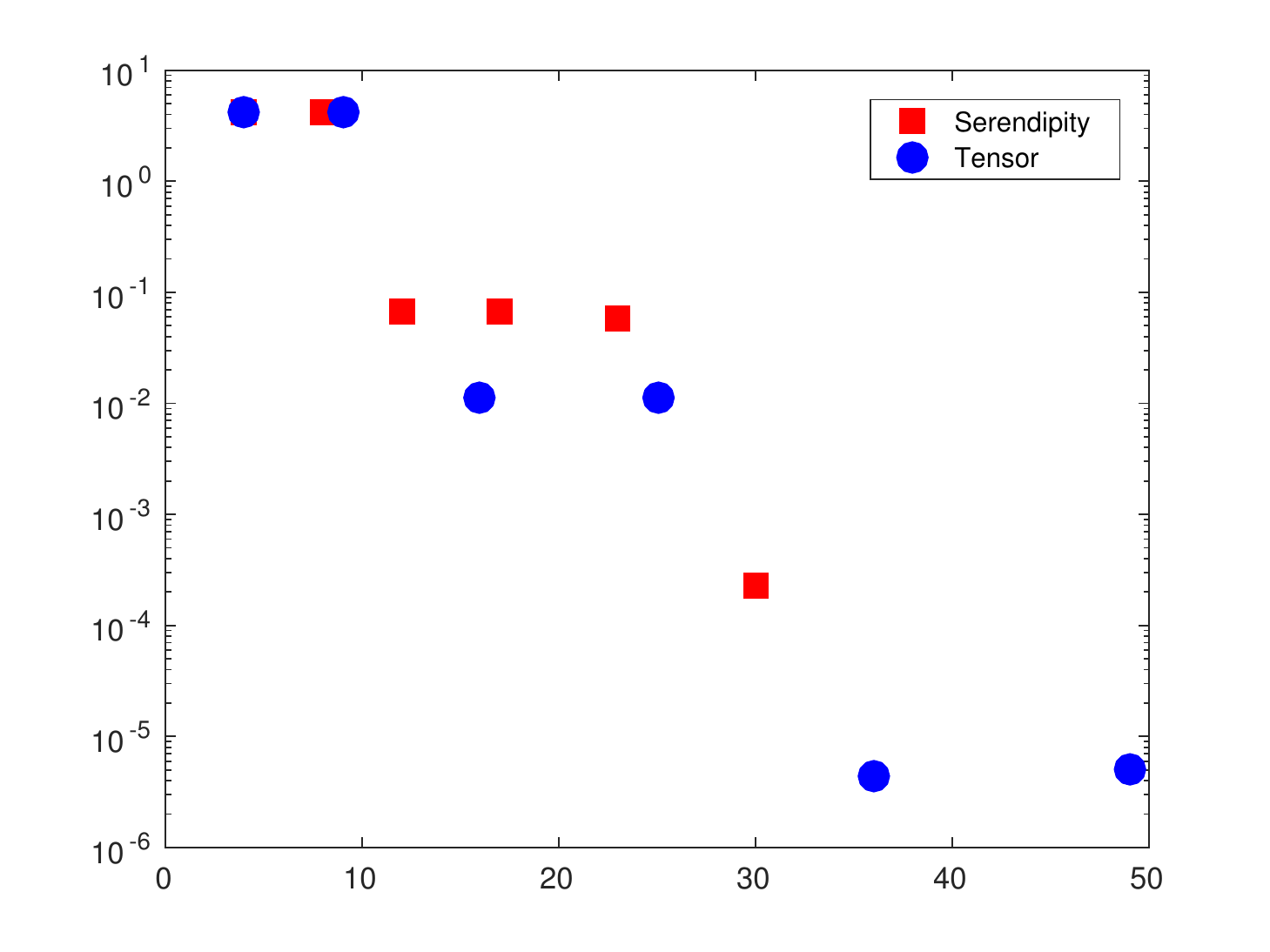} &
\includegraphics[width=.31\textwidth]{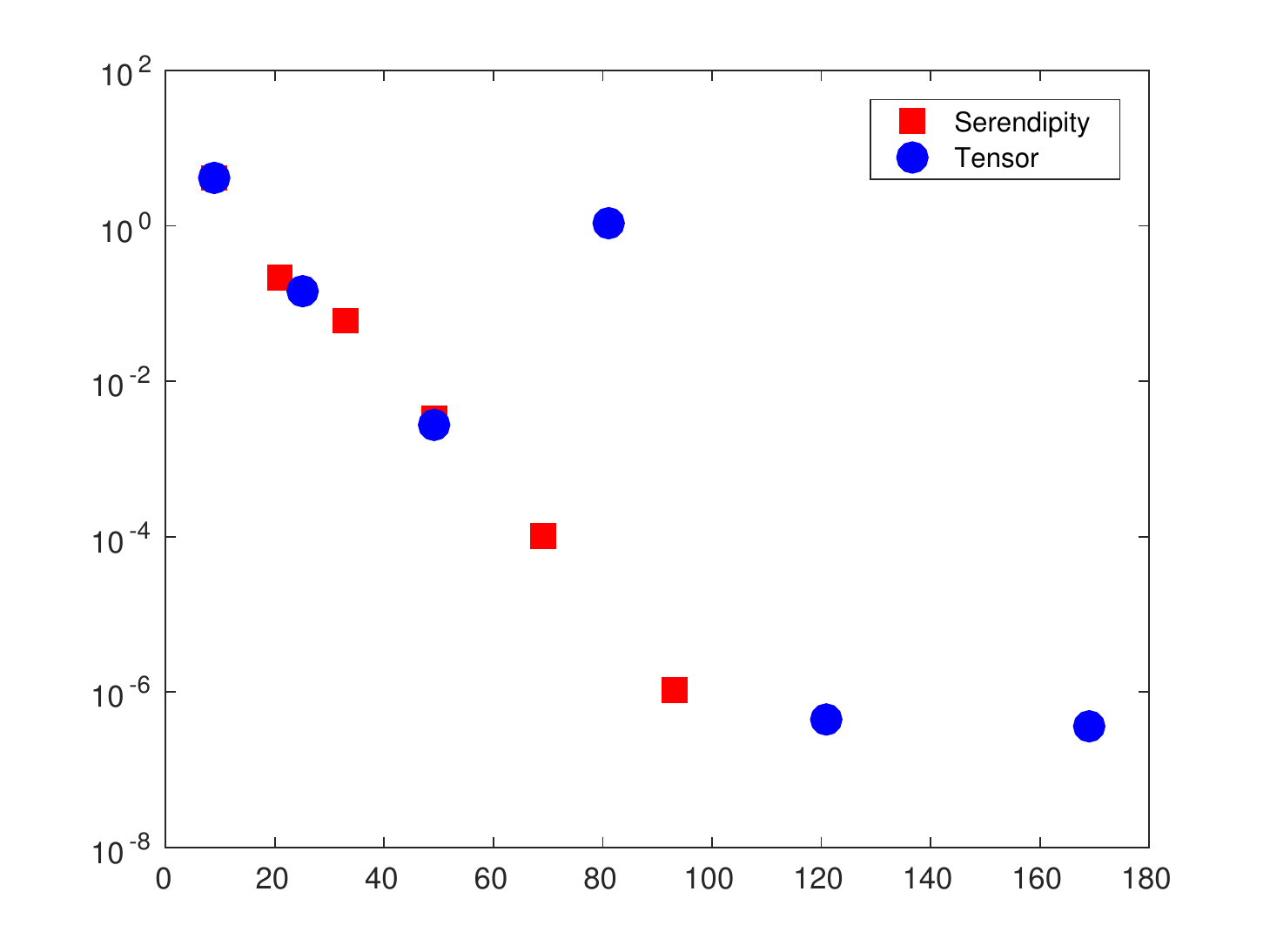} &
\includegraphics[width=.31\textwidth]{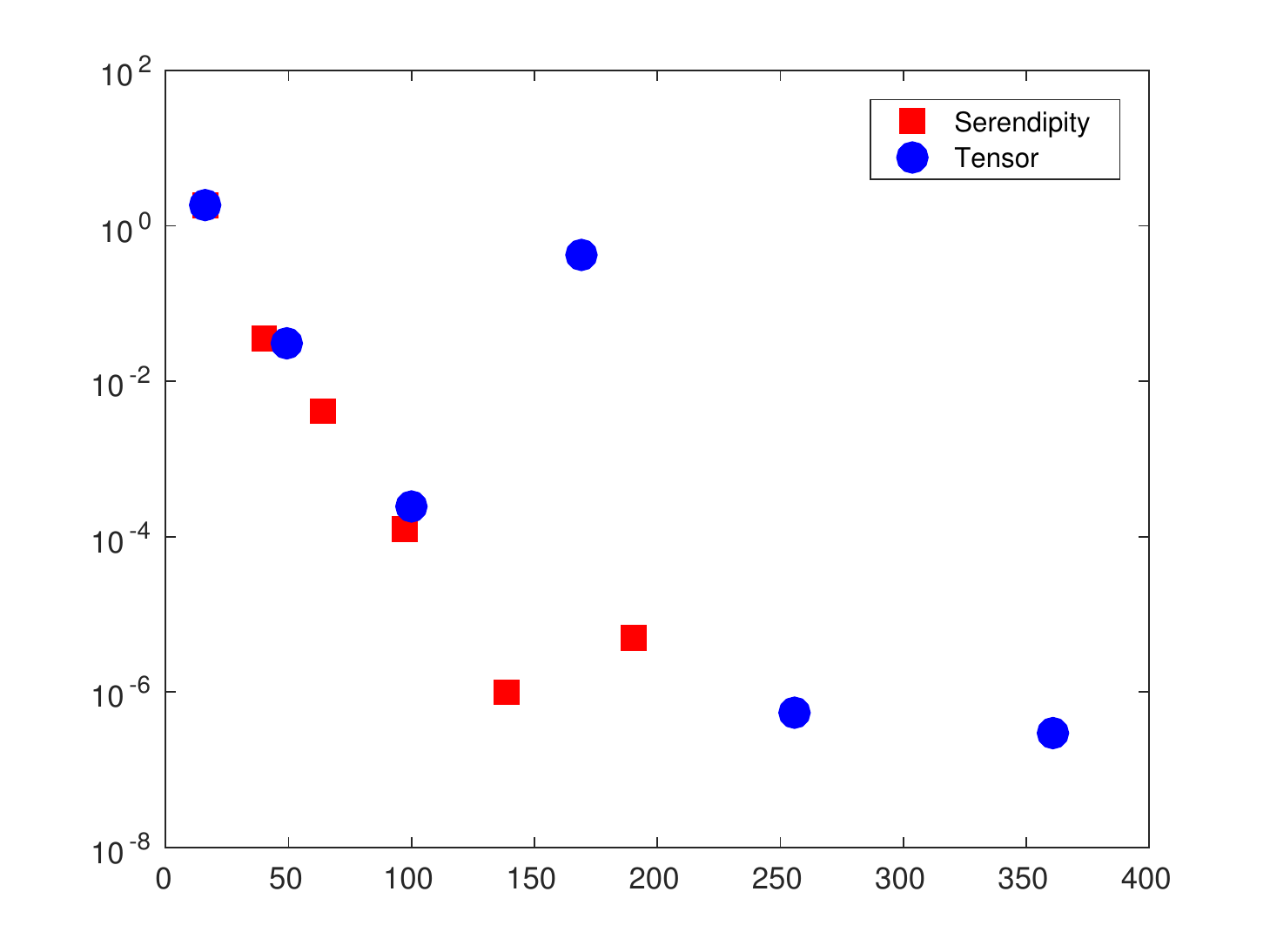} \\[-4mm]
h=1 & 
h=1/2 &
h=1/3 \\[2mm]
\includegraphics[width=.31\textwidth]{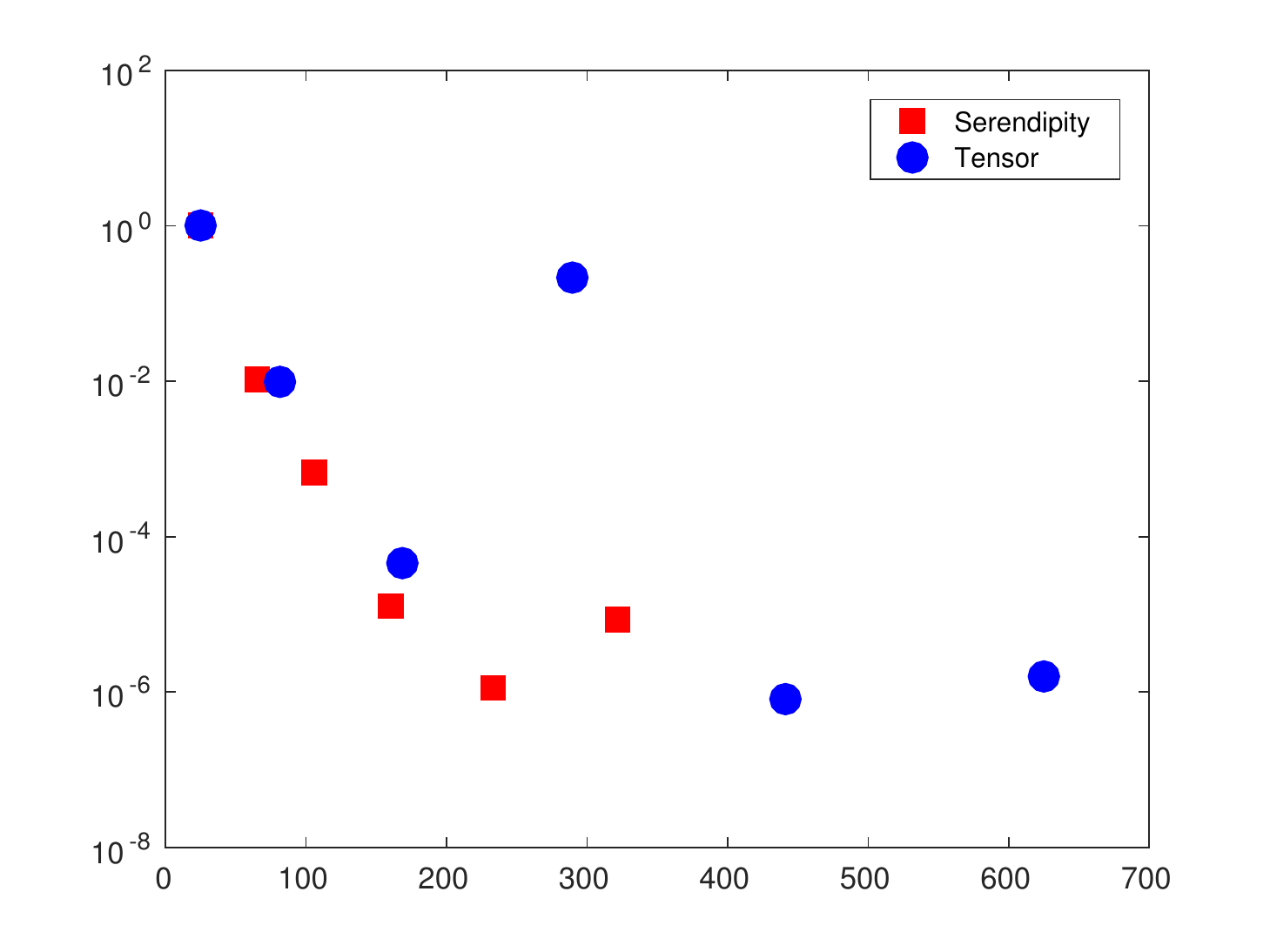} &
\includegraphics[width=.31\textwidth]{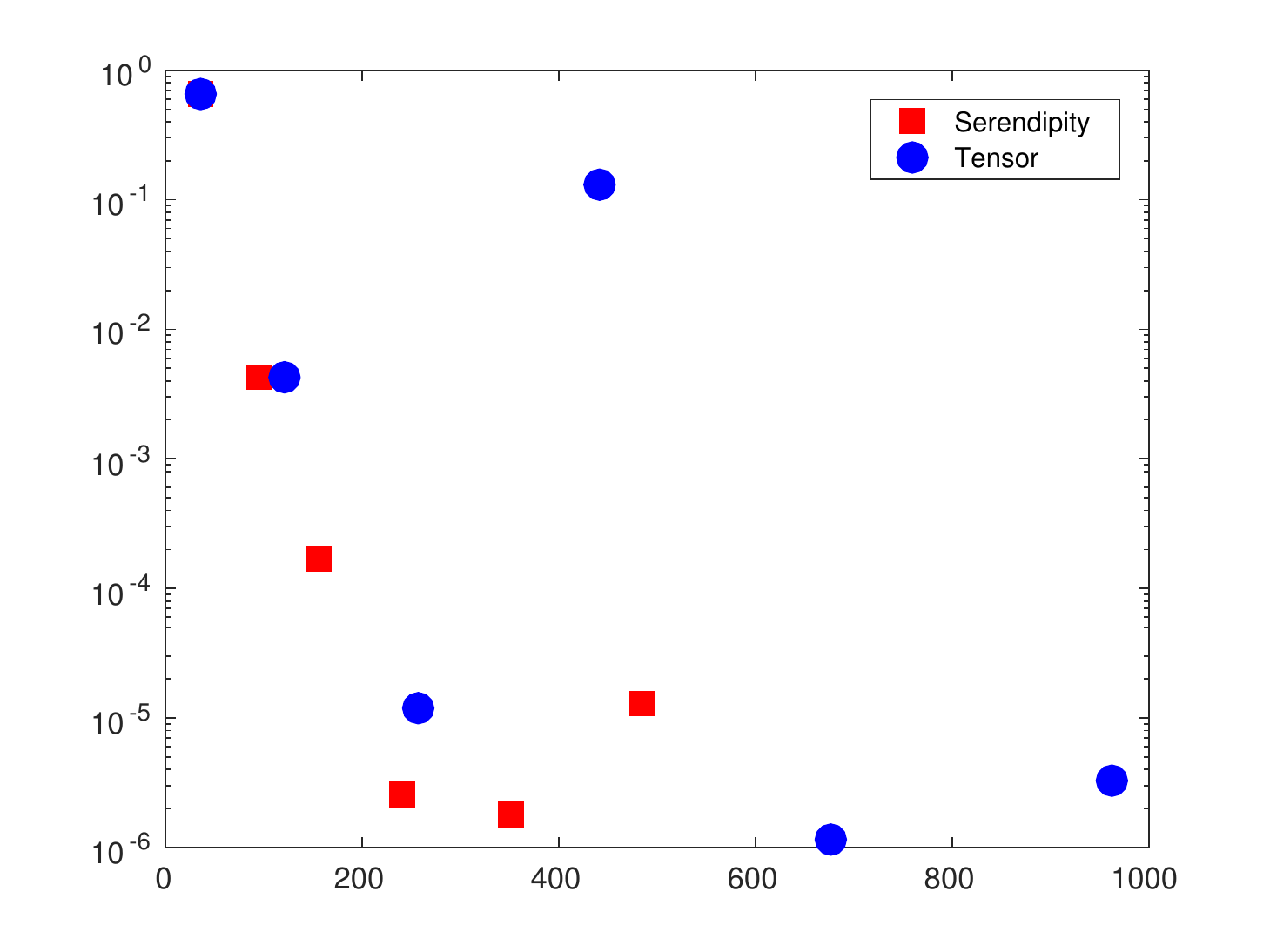} \\[-4mm]
h=1/4 &
h=1/5
\end{array}$
\caption{
Square domain, Neumann conditions, $\lambda=2\pi^2$, $p$-refinement experiments.
}
\label{fig:NeuTwoPiK}
%
%
%
$\begin{array}{ccc}
 &&\\
  &
\includegraphics[width=.31\textwidth]{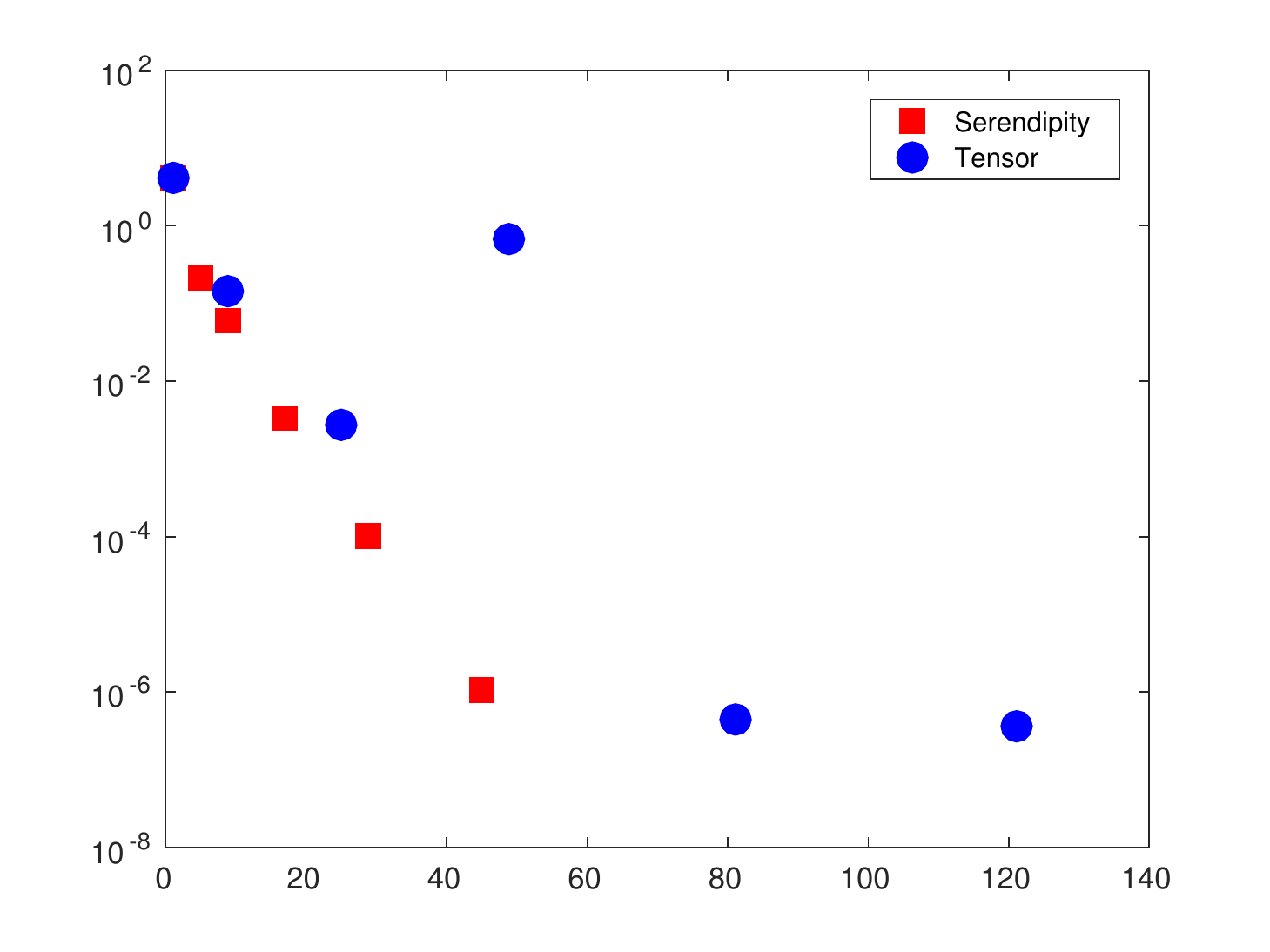} &
\includegraphics[width=.31\textwidth]{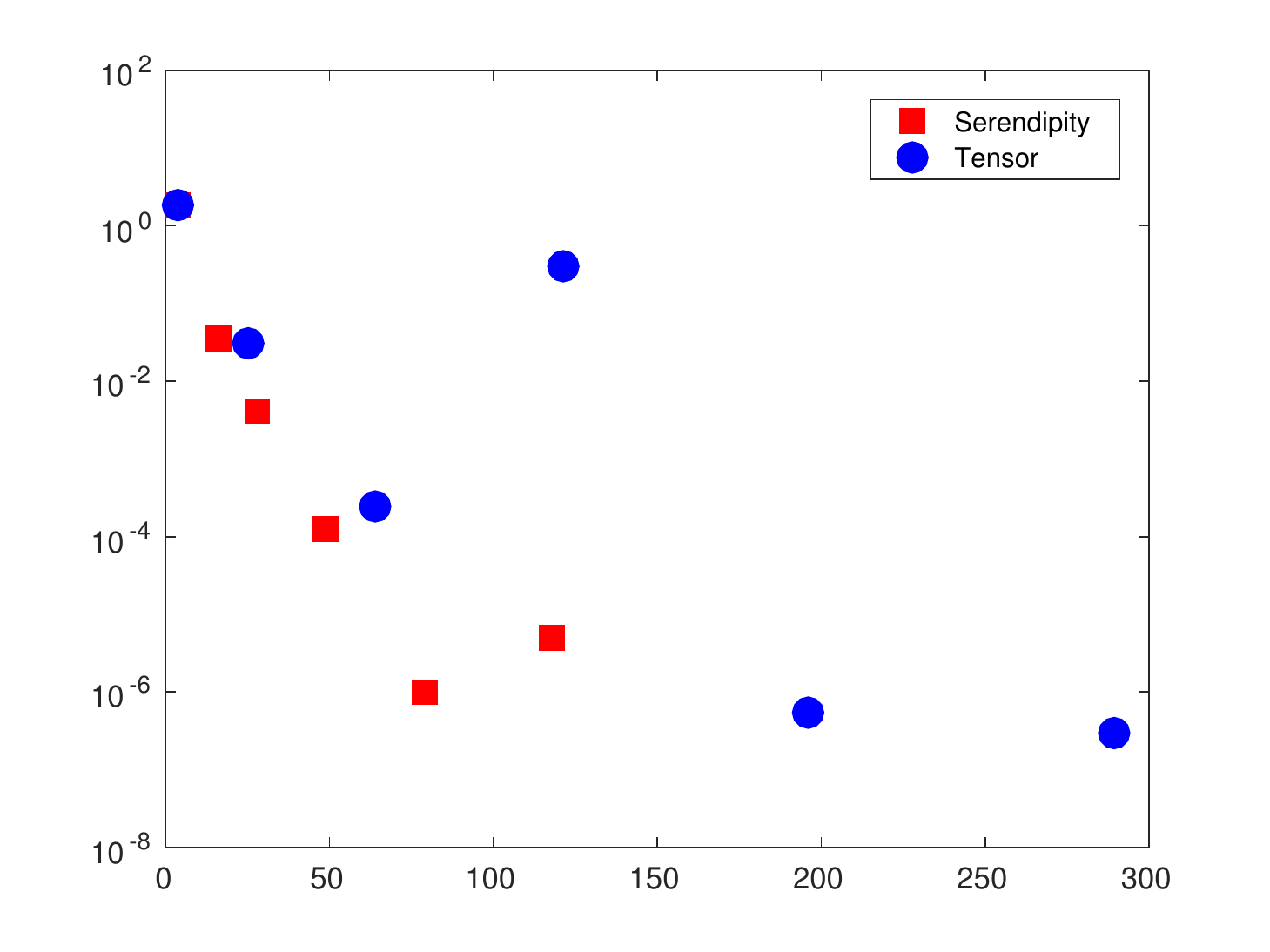} \\[-4mm]
 & 
h=1/2 &
h=1/3 \\[2mm]
\includegraphics[width=.31\textwidth]{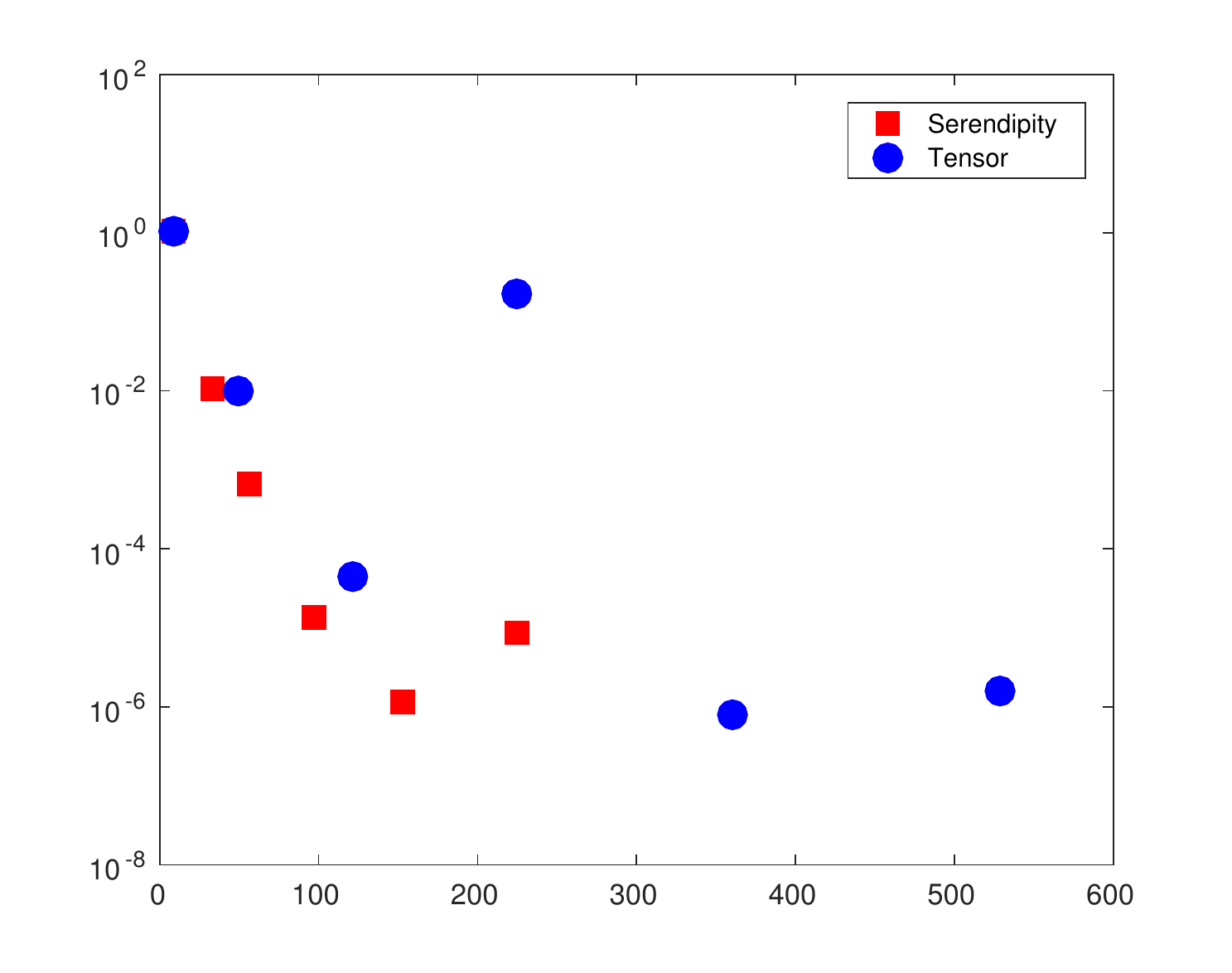} &
\includegraphics[width=.31\textwidth]{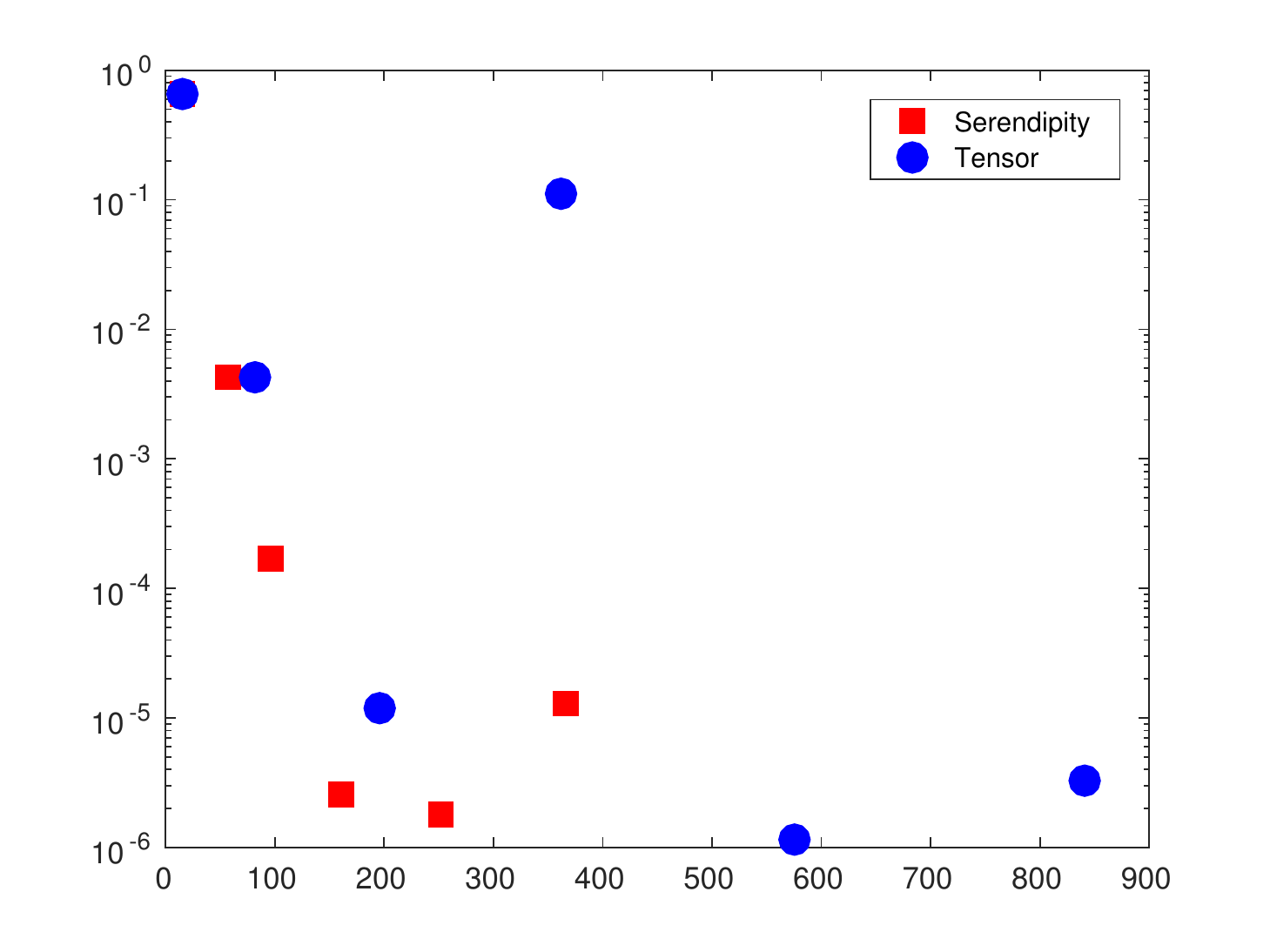} \\[-4mm]
h=1/4 &
h=1/5
\end{array}$
	\caption{
	Square domain, Dirichlet conditions, $\lambda=2\pi^2$, $p$-refinement experiments.
	}
	\label{fig:DirTwoPiK}
\end{figure}

\begin{figure}[ht]
$\begin{array}{ccc}
\includegraphics[width=.31\textwidth]{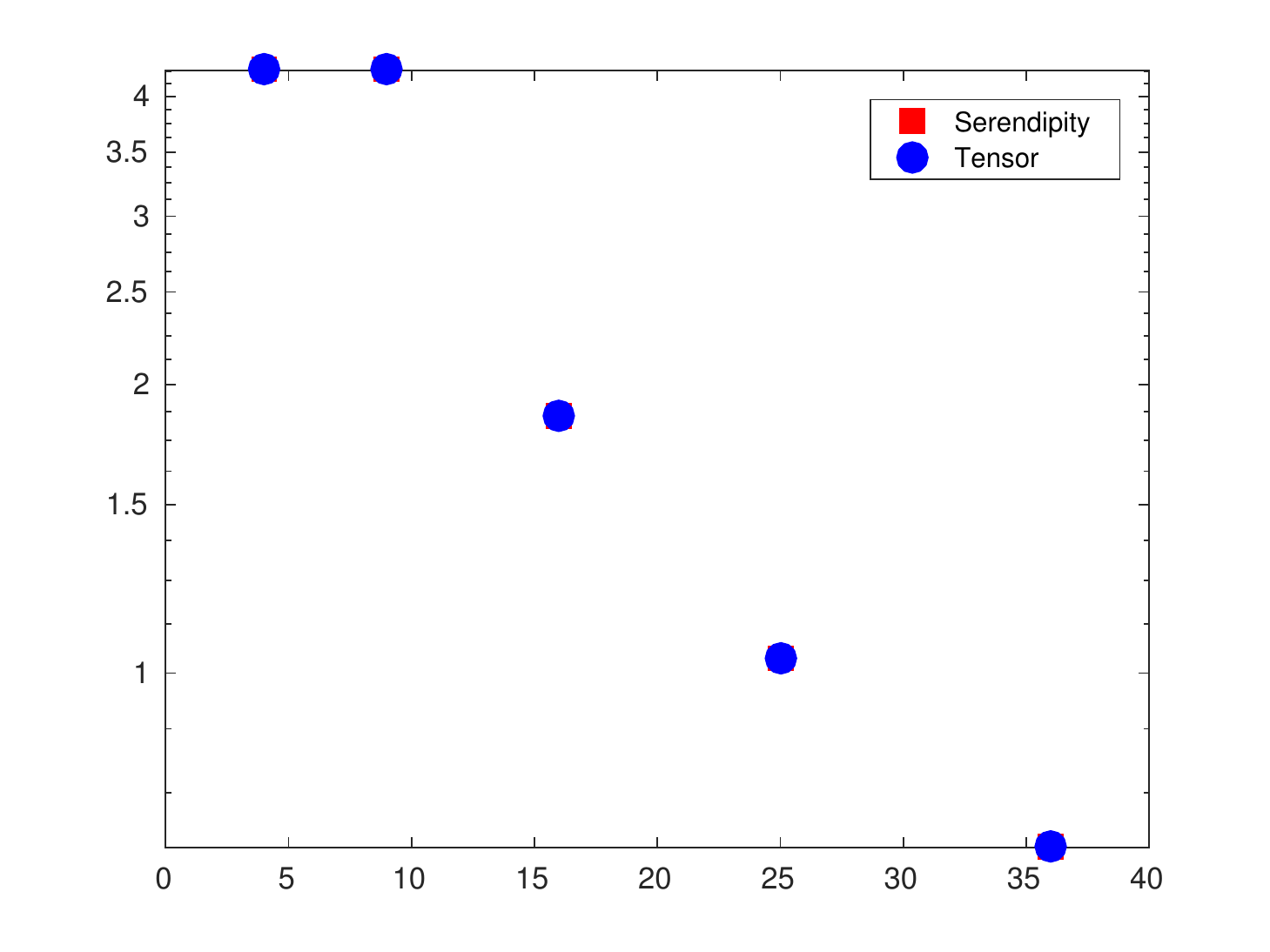} &
\includegraphics[width=.31\textwidth]{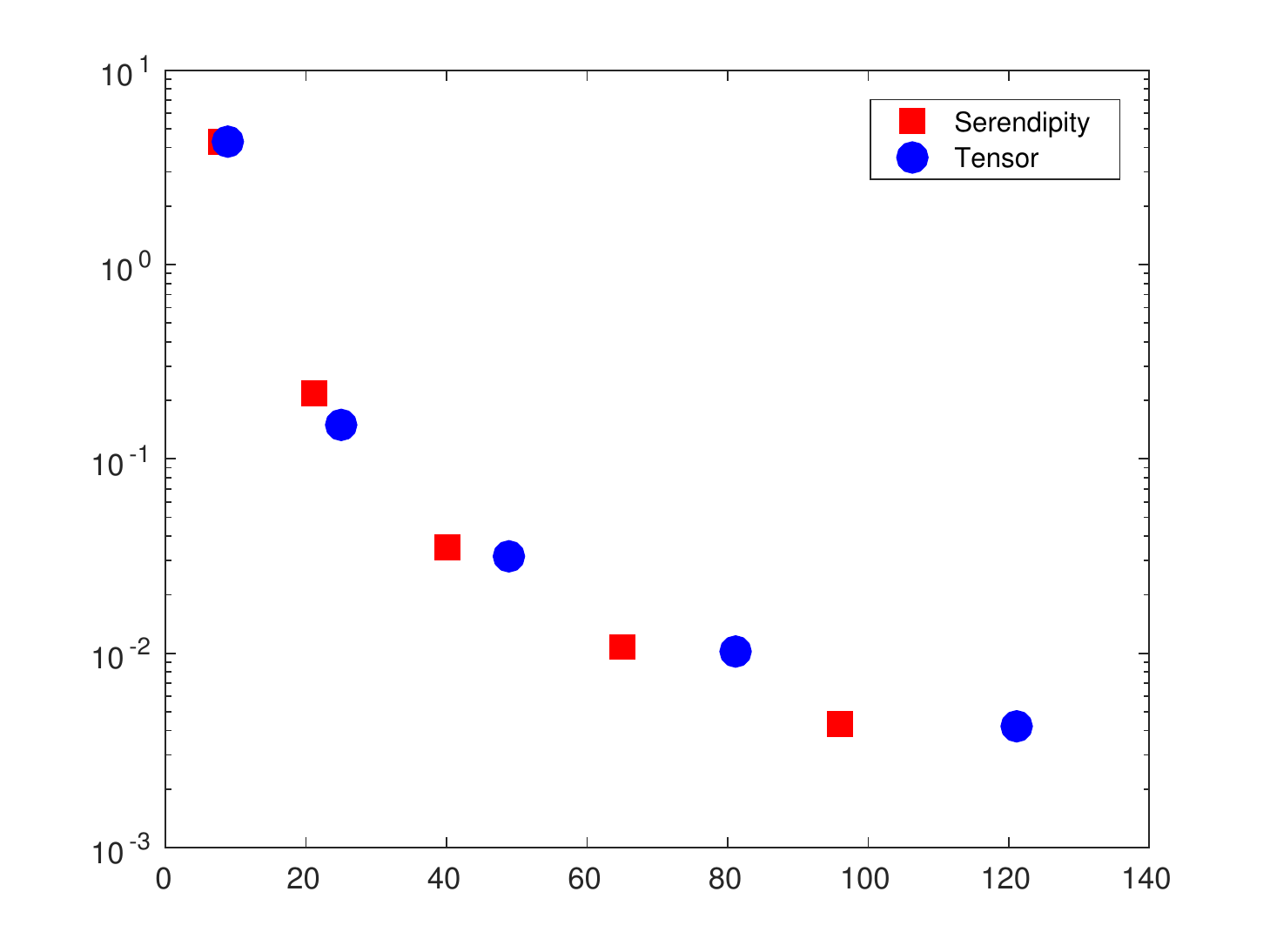} &
\includegraphics[width=.31\textwidth]{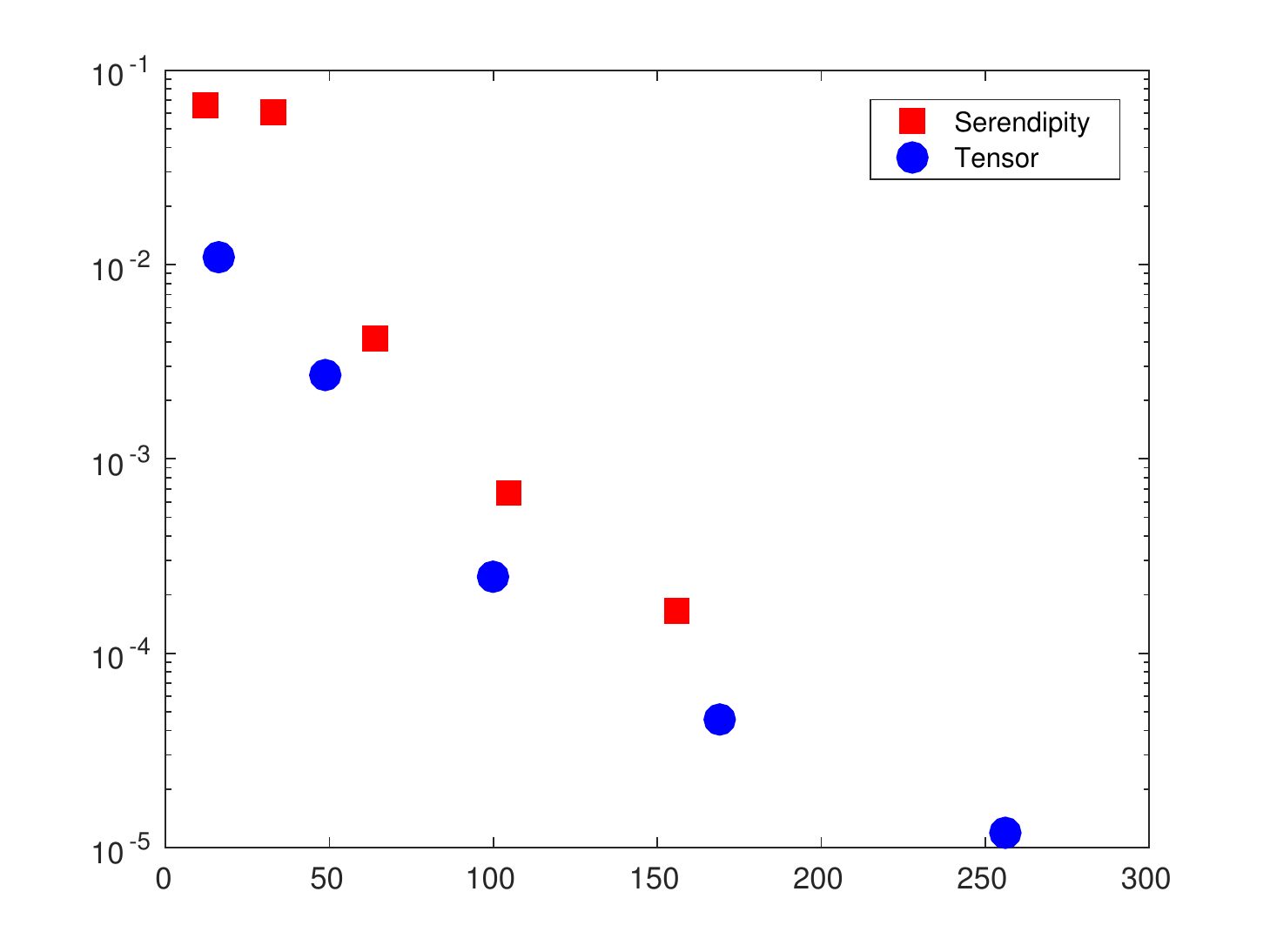} \\[-4mm]
p=1 & 
p=2 &
p=3 \\[0mm]
\includegraphics[width=.31\textwidth]{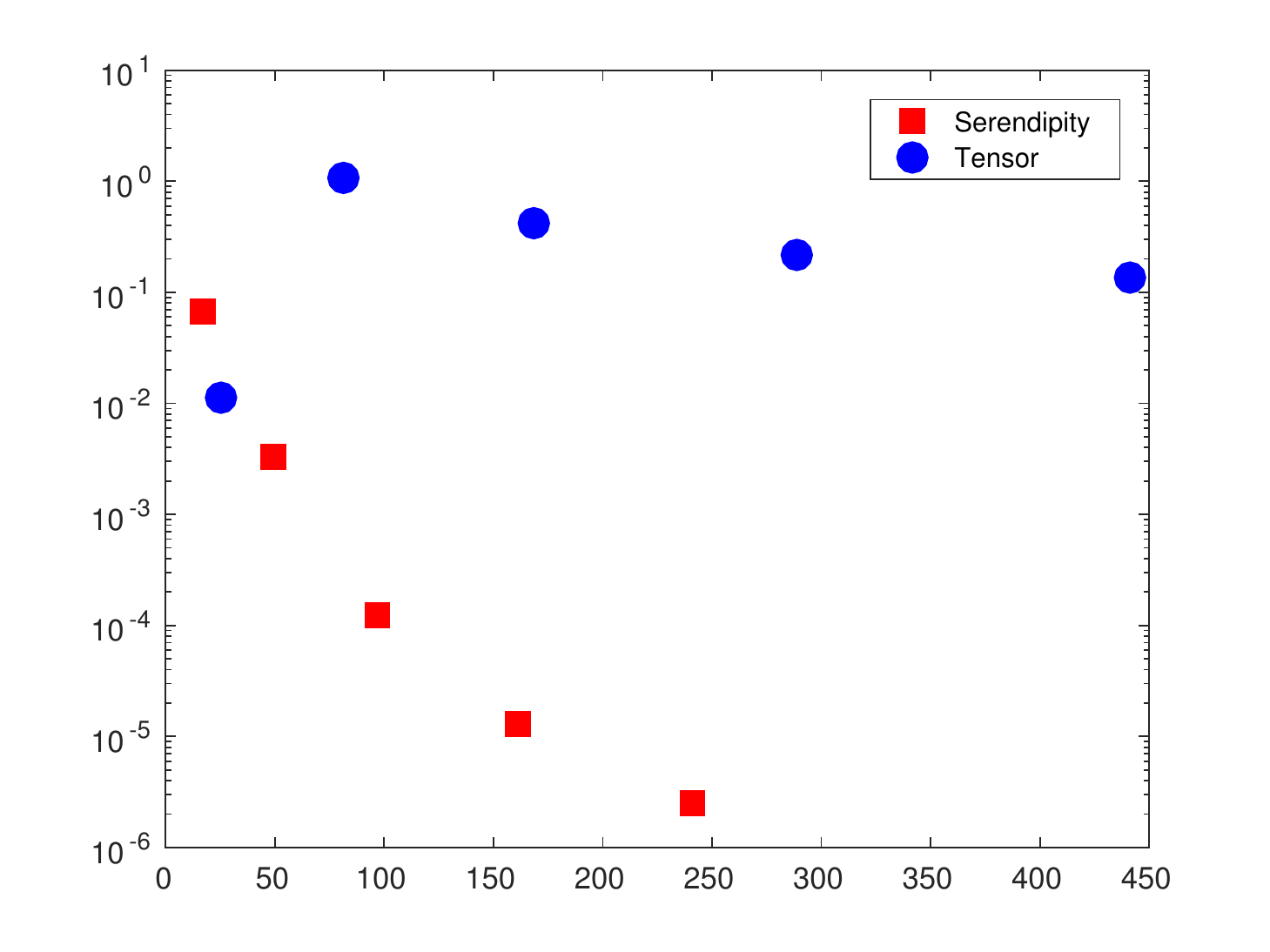} &
\includegraphics[width=.31\textwidth]{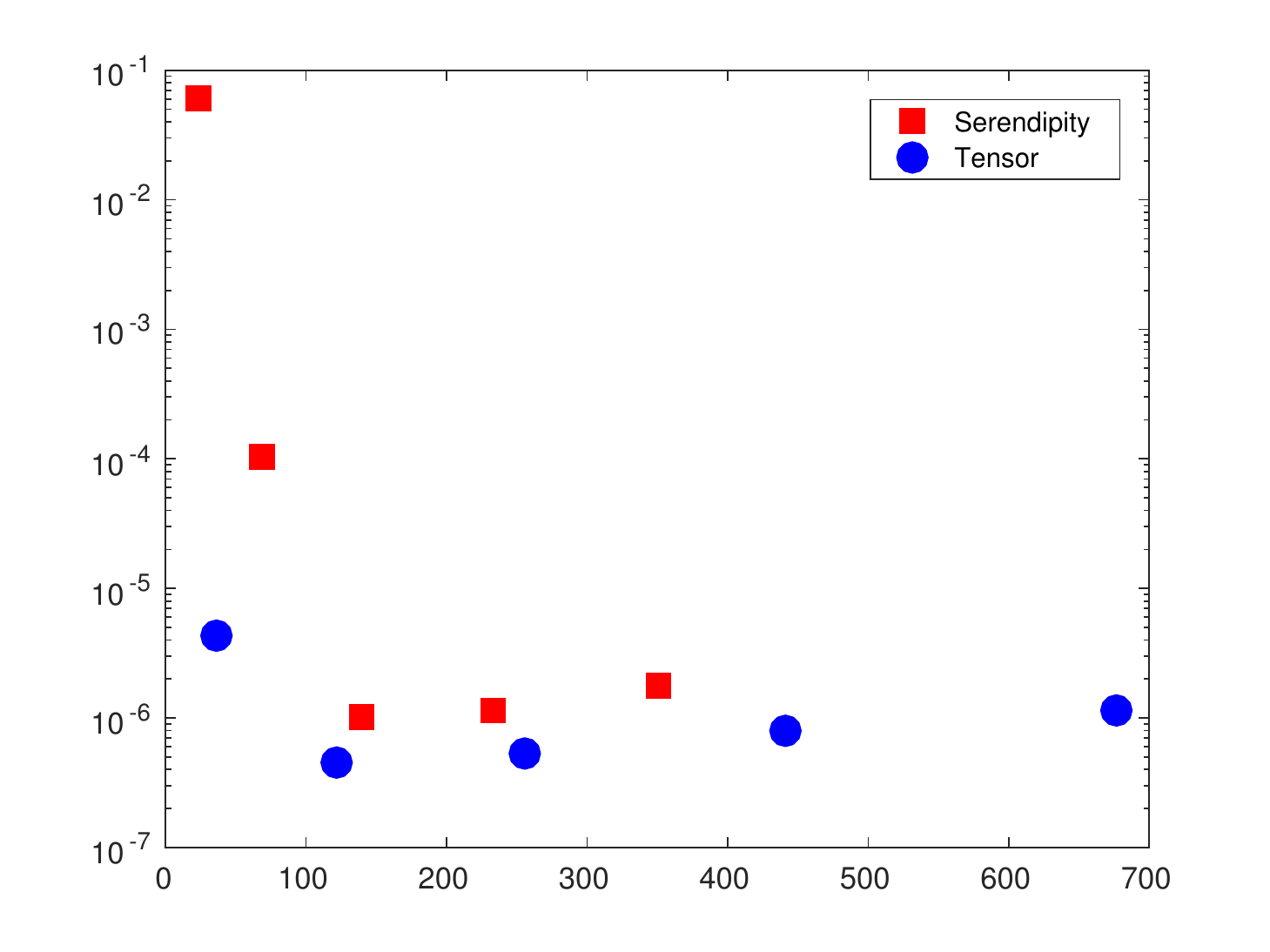} &
\includegraphics[width=.31\textwidth]{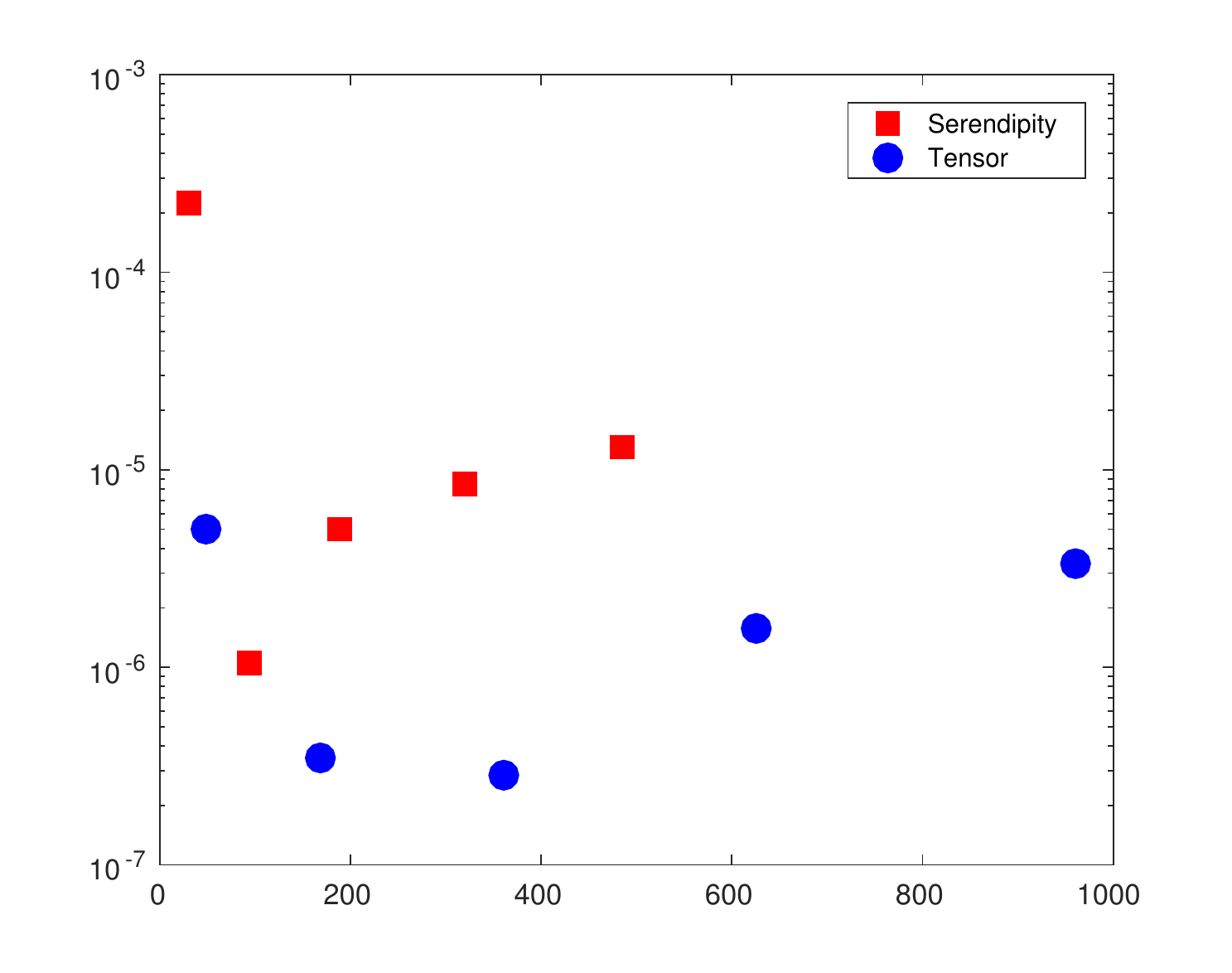} \\[-4mm]
p=4 &
p=5 &
p=6
\end{array}$
	\caption{
	Square domain, Neumann conditions, $\lambda=2\pi^2$, $h$-refinement experiments.
	}
	\label{fig:NeuTwoPiP}
\end{figure}

\begin{figure}[ht]
$\begin{array}{ccc}
\includegraphics[width=.31\textwidth]{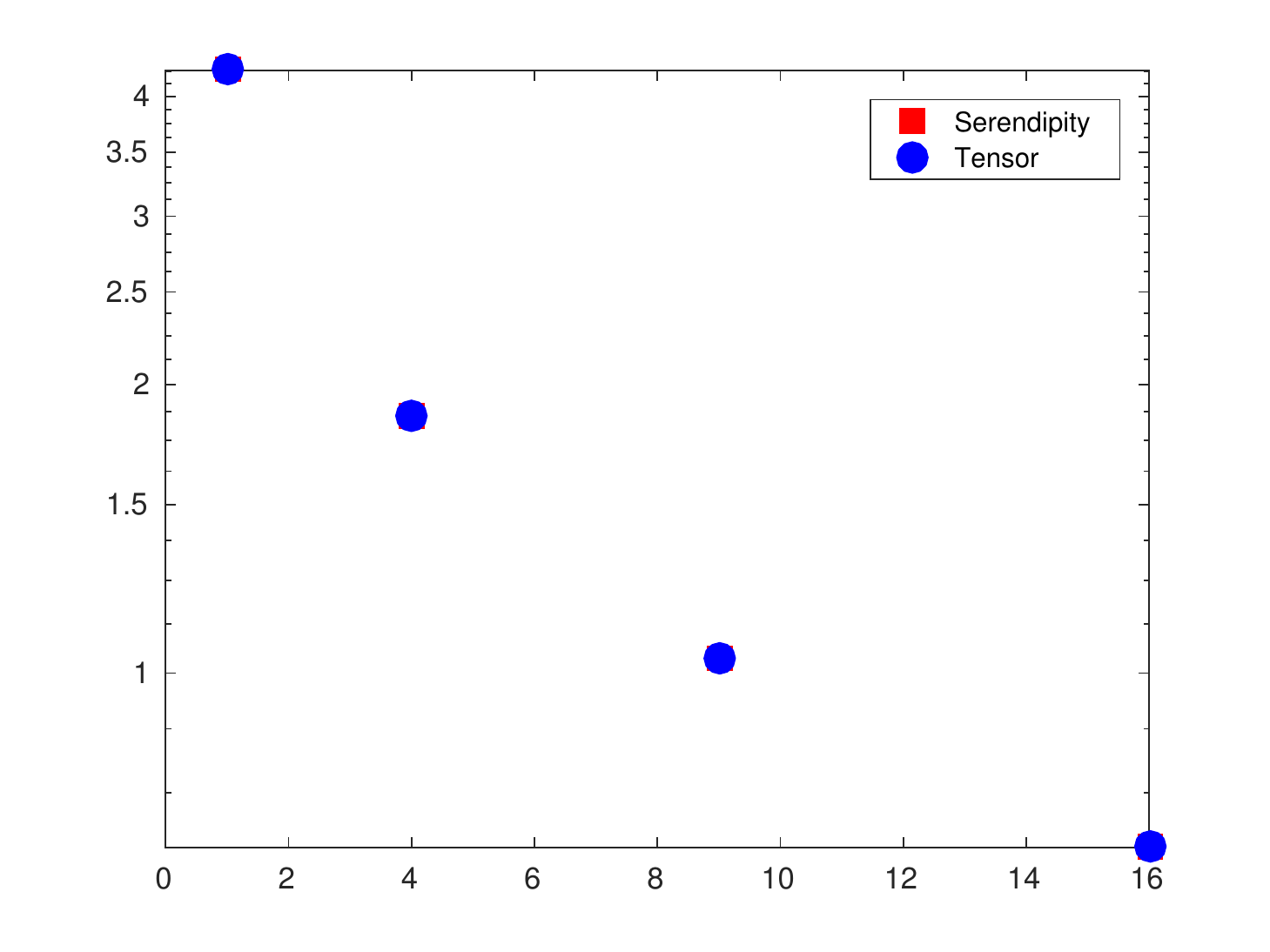} &
\includegraphics[width=.31\textwidth]{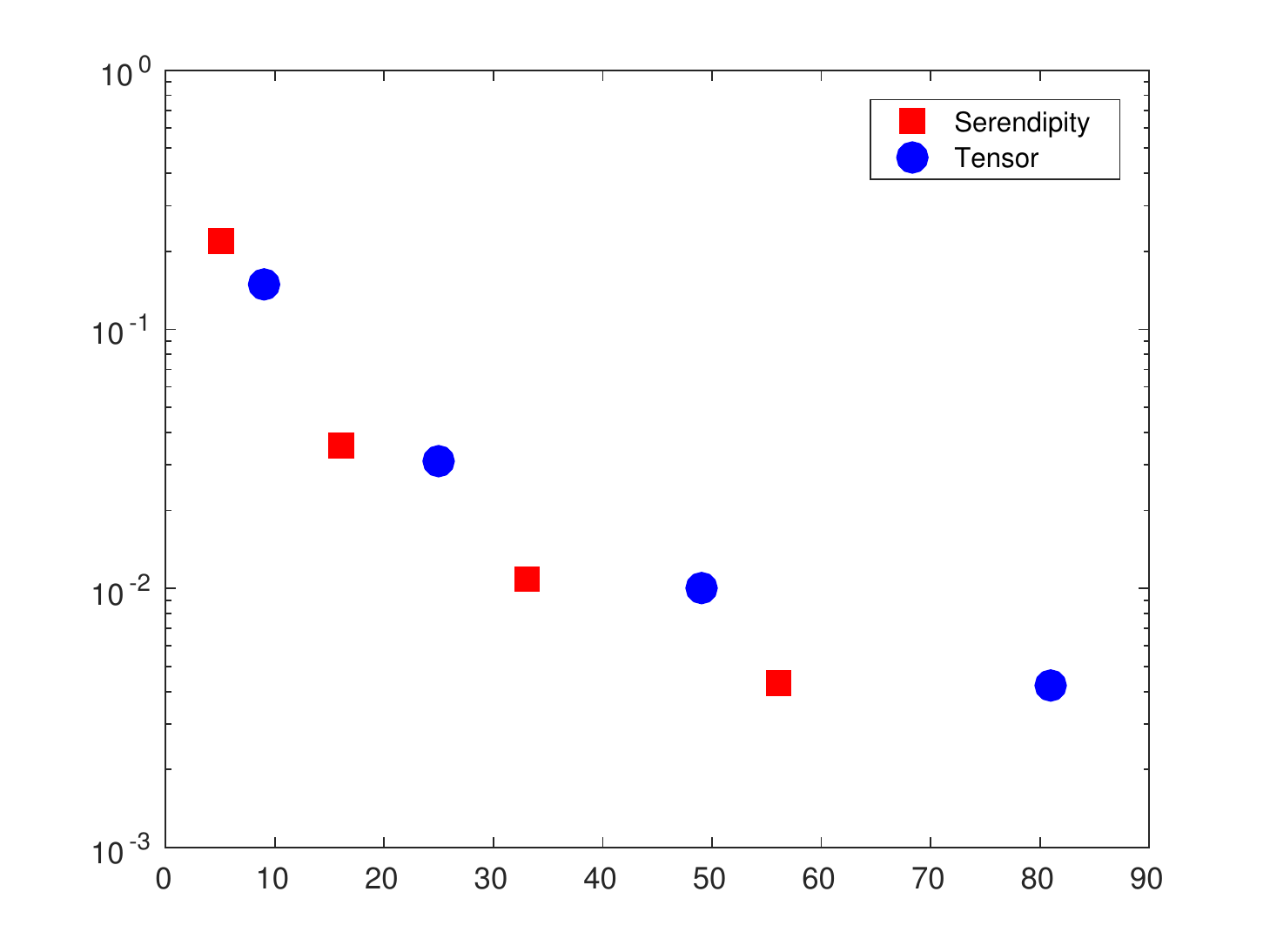} &
\includegraphics[width=.31\textwidth]{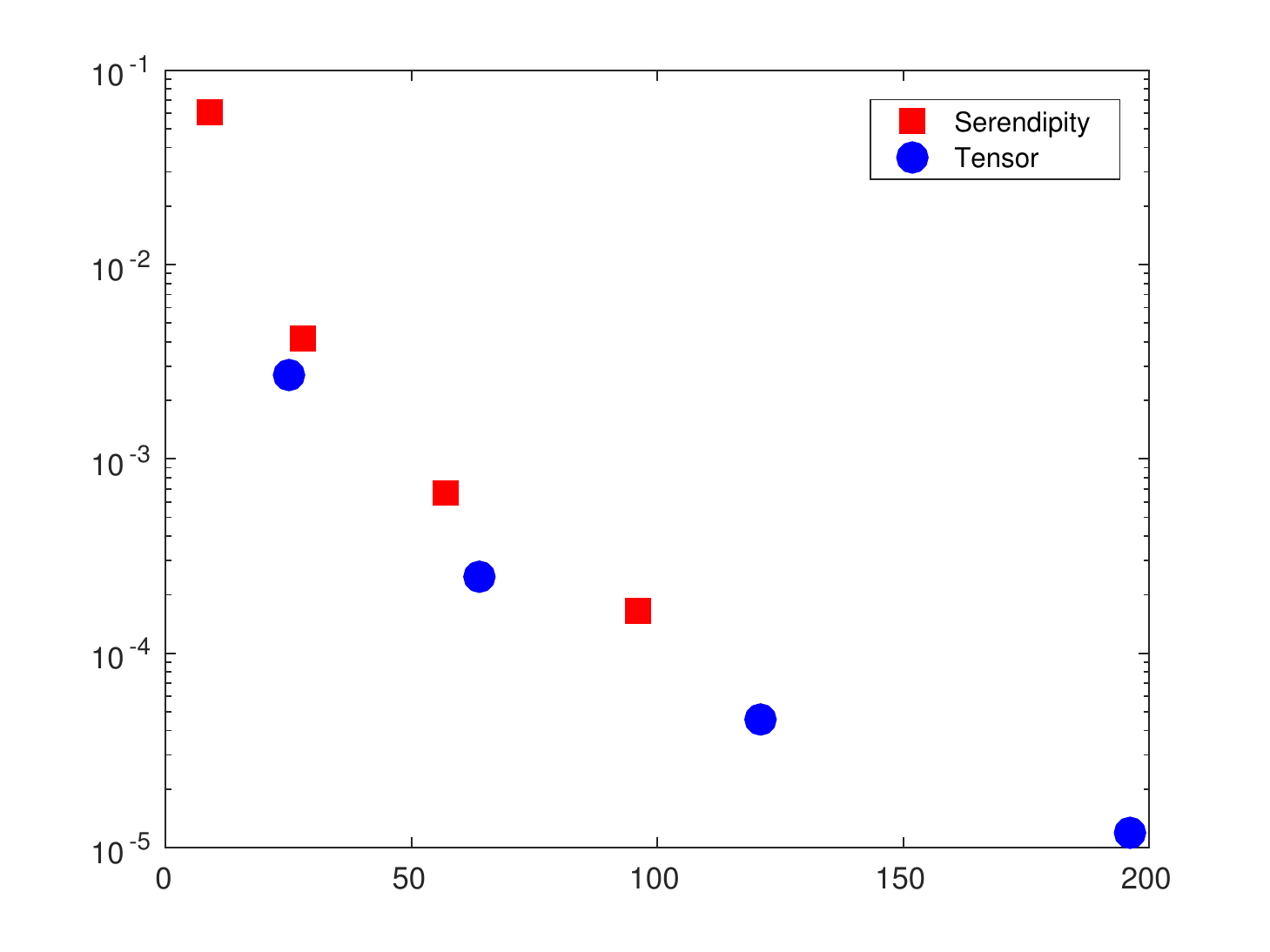} \\[-4mm]
p=1 & 
p=2 &
p=3 \\[2mm]
\includegraphics[width=.31\textwidth]{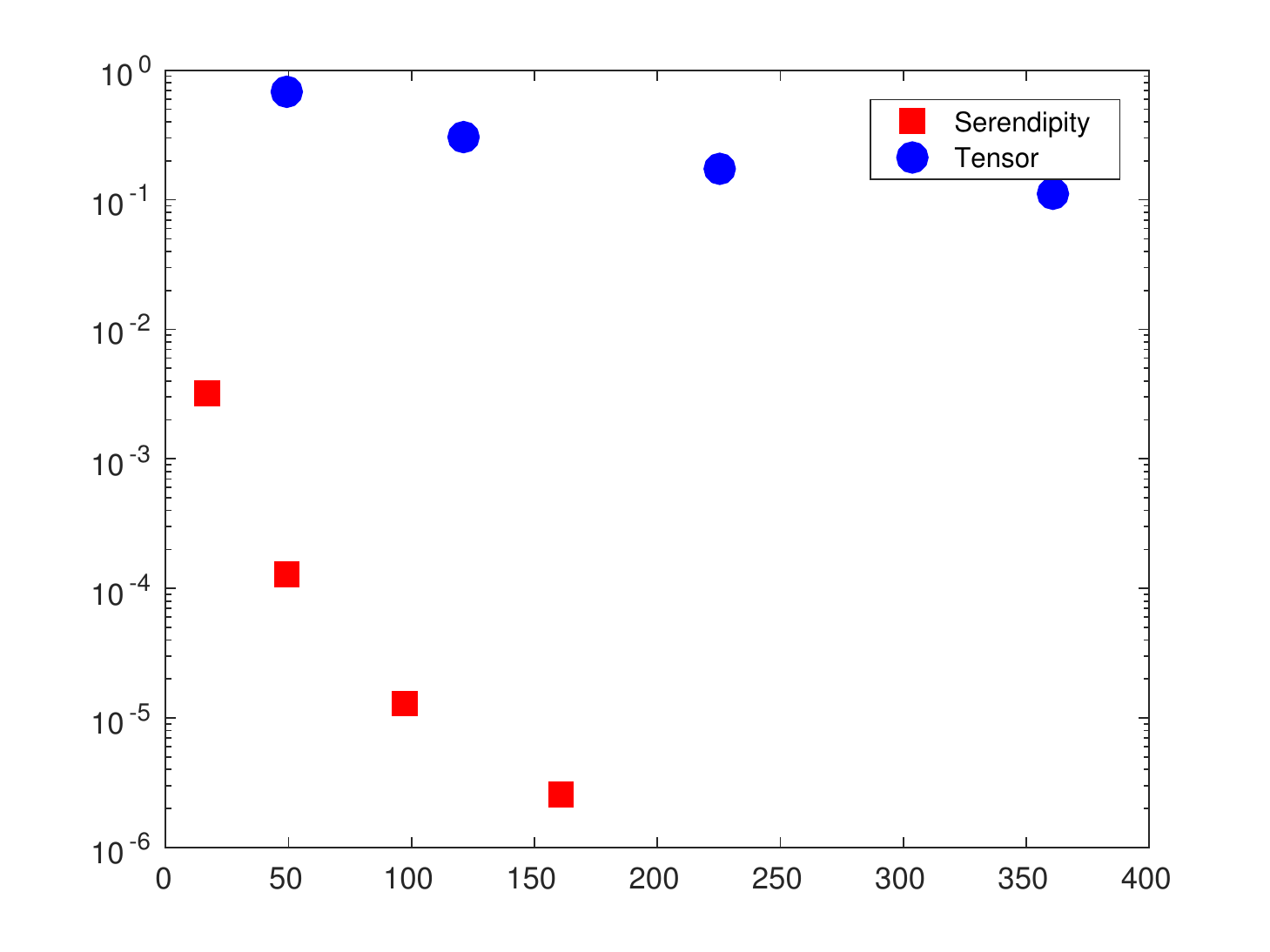} &
\includegraphics[width=.31\textwidth]{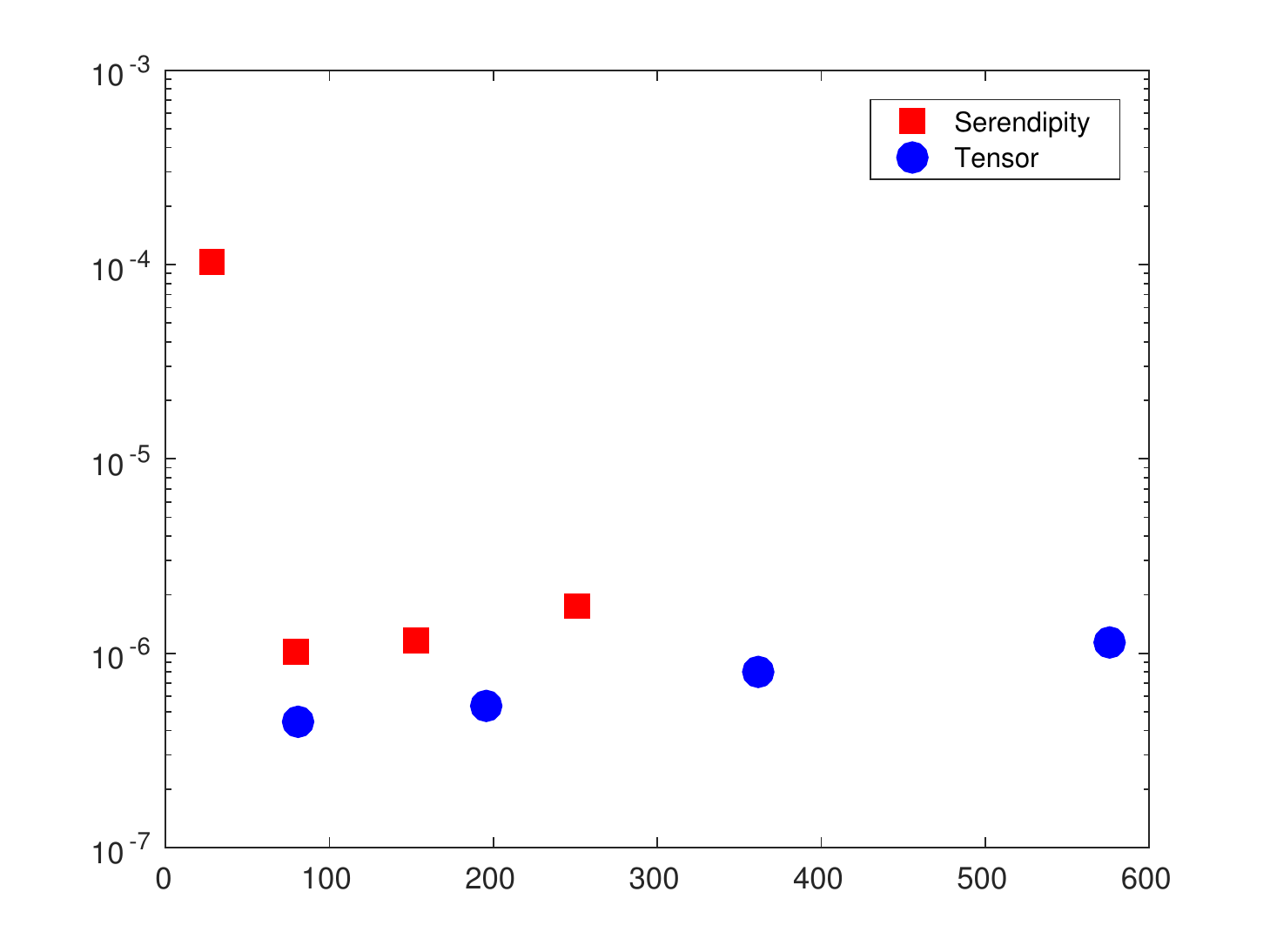} &
\includegraphics[width=.31\textwidth]{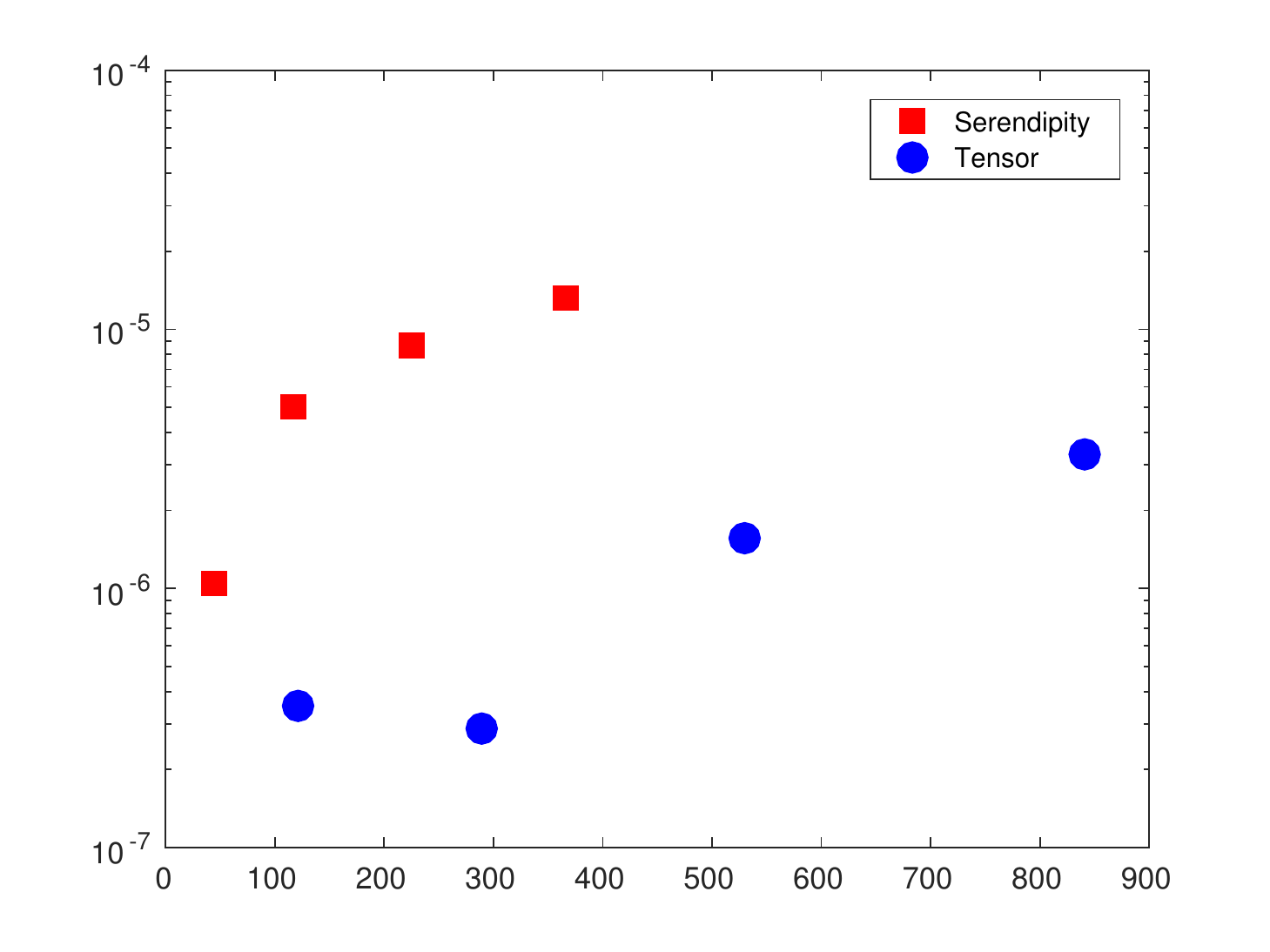} \\[-4mm]
p=4 &
p=5 &
p=6
\end{array}$
	\caption{
	Square domain, Dirichlet conditions, $\lambda=2\pi^2$, $h$-refinement experiments.
	}
	\label{fig:DirTwoPiP}
\end{figure}

\begin{figure}[ht]
$\begin{array}{ccc}
\includegraphics[width=.31\textwidth]{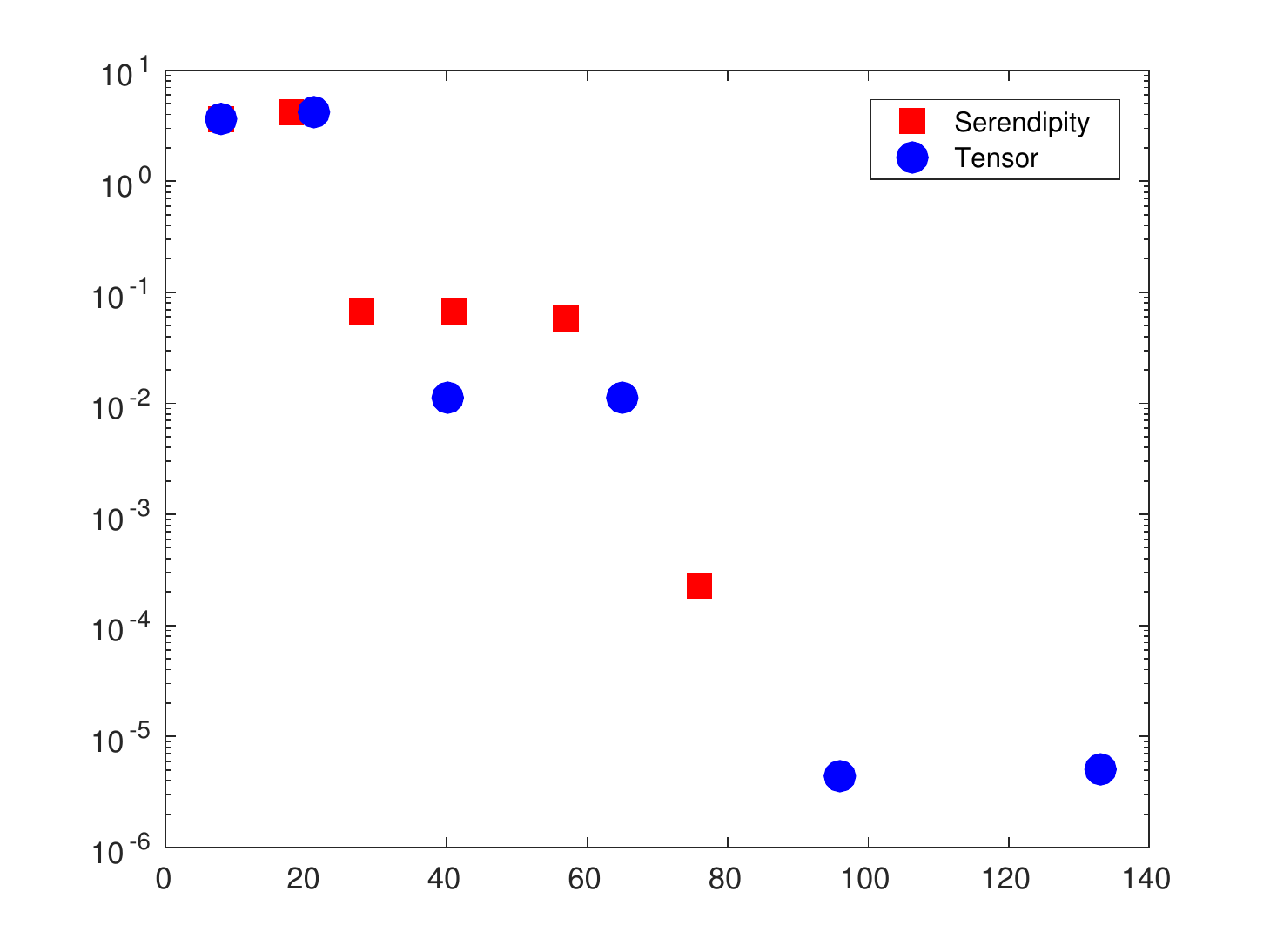} &
\includegraphics[width=.31\textwidth]{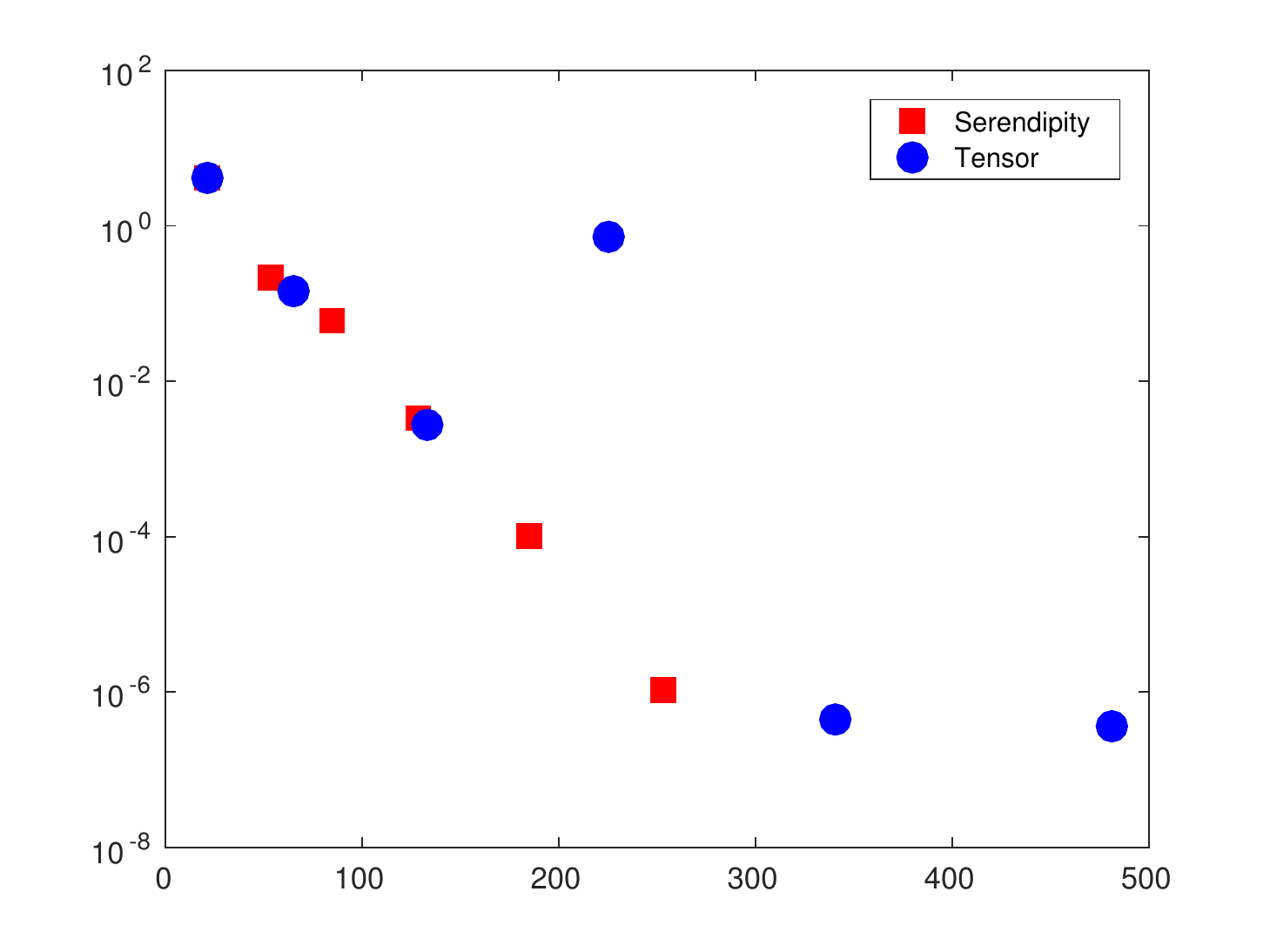} &
\includegraphics[width=.31\textwidth]{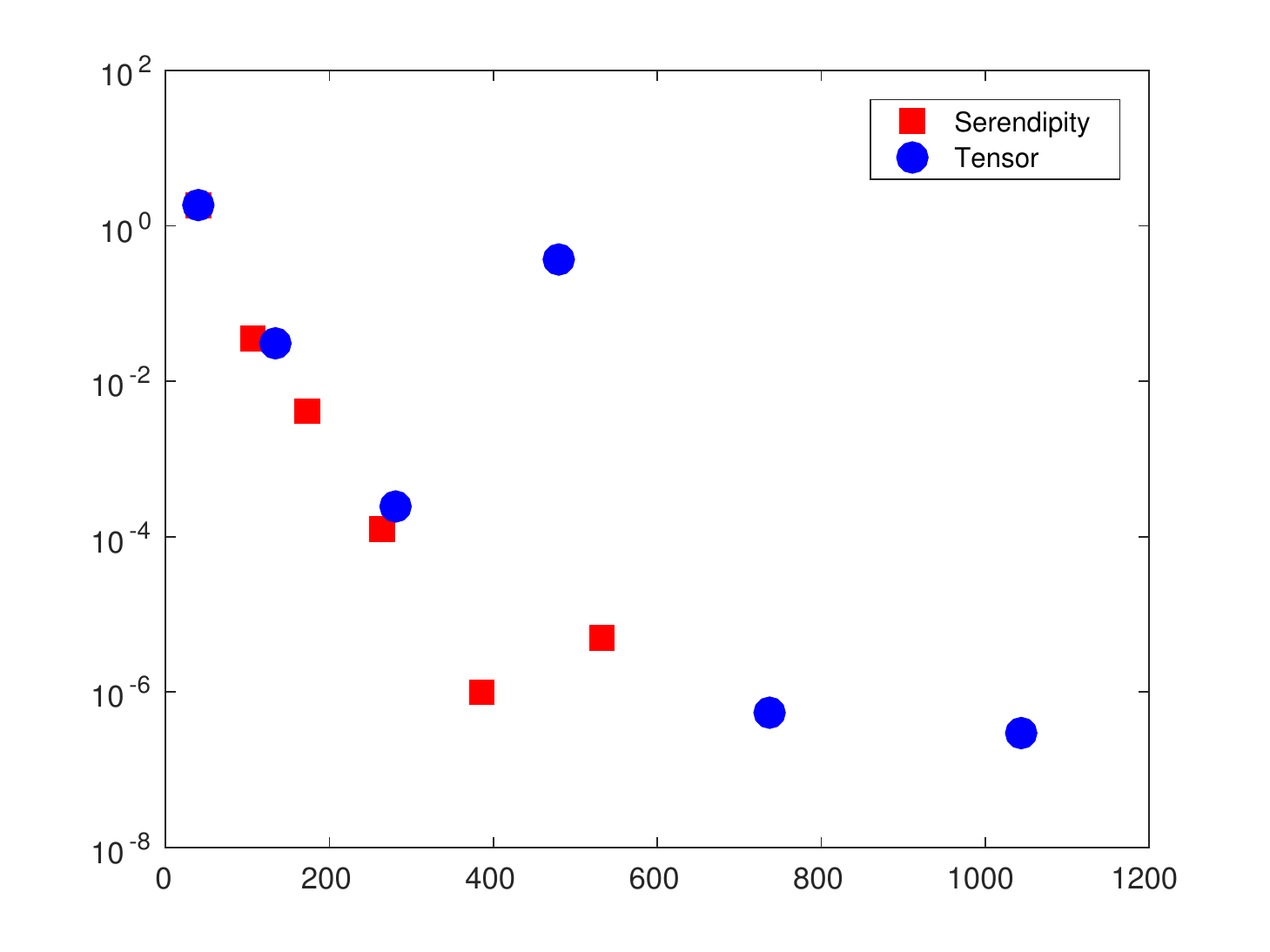} \\[-4mm]
h=1 & 
h=1/2 &
h=1/3 \\[2mm]
\includegraphics[width=.31\textwidth]{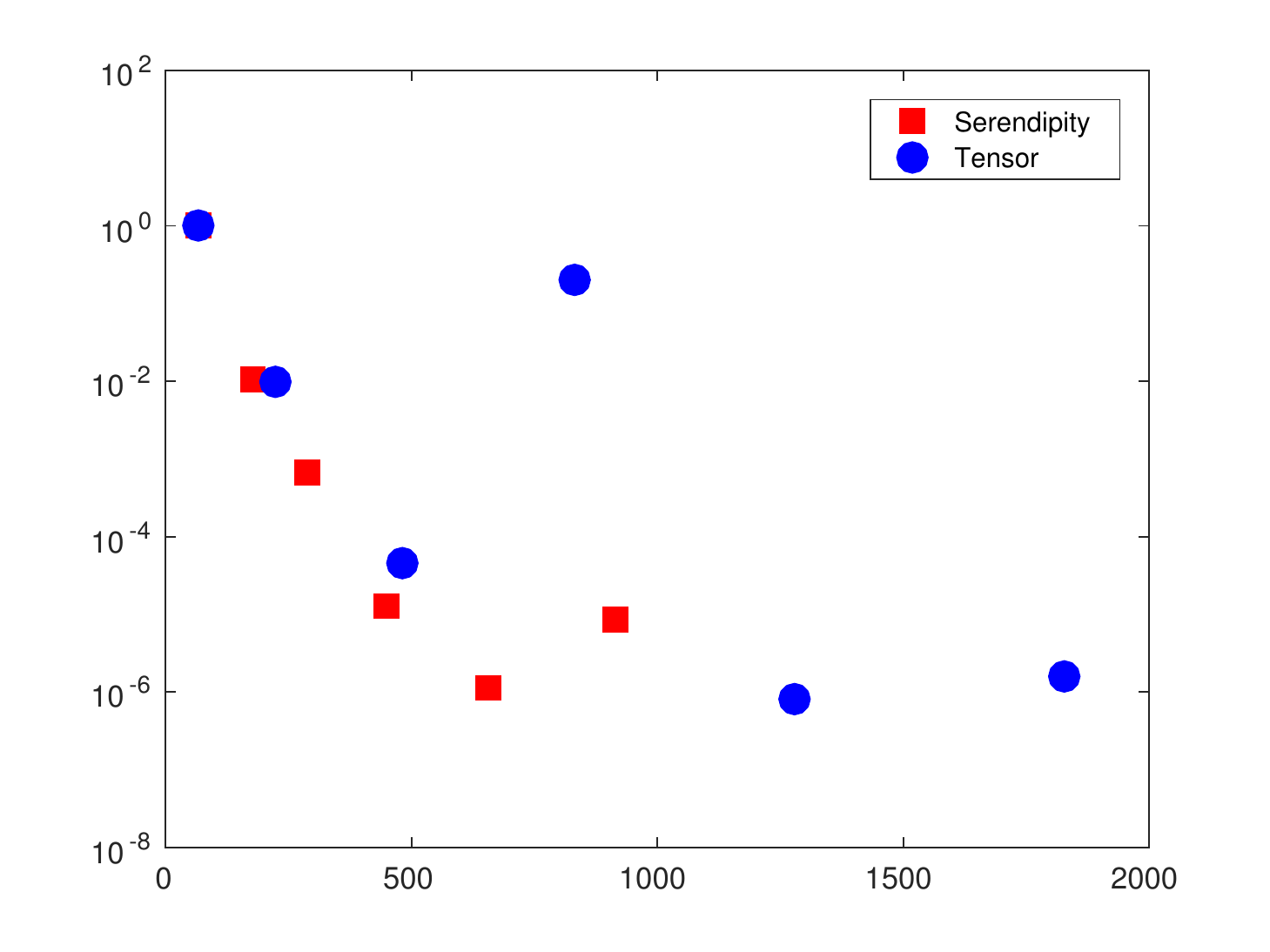} &
\includegraphics[width=.31\textwidth]{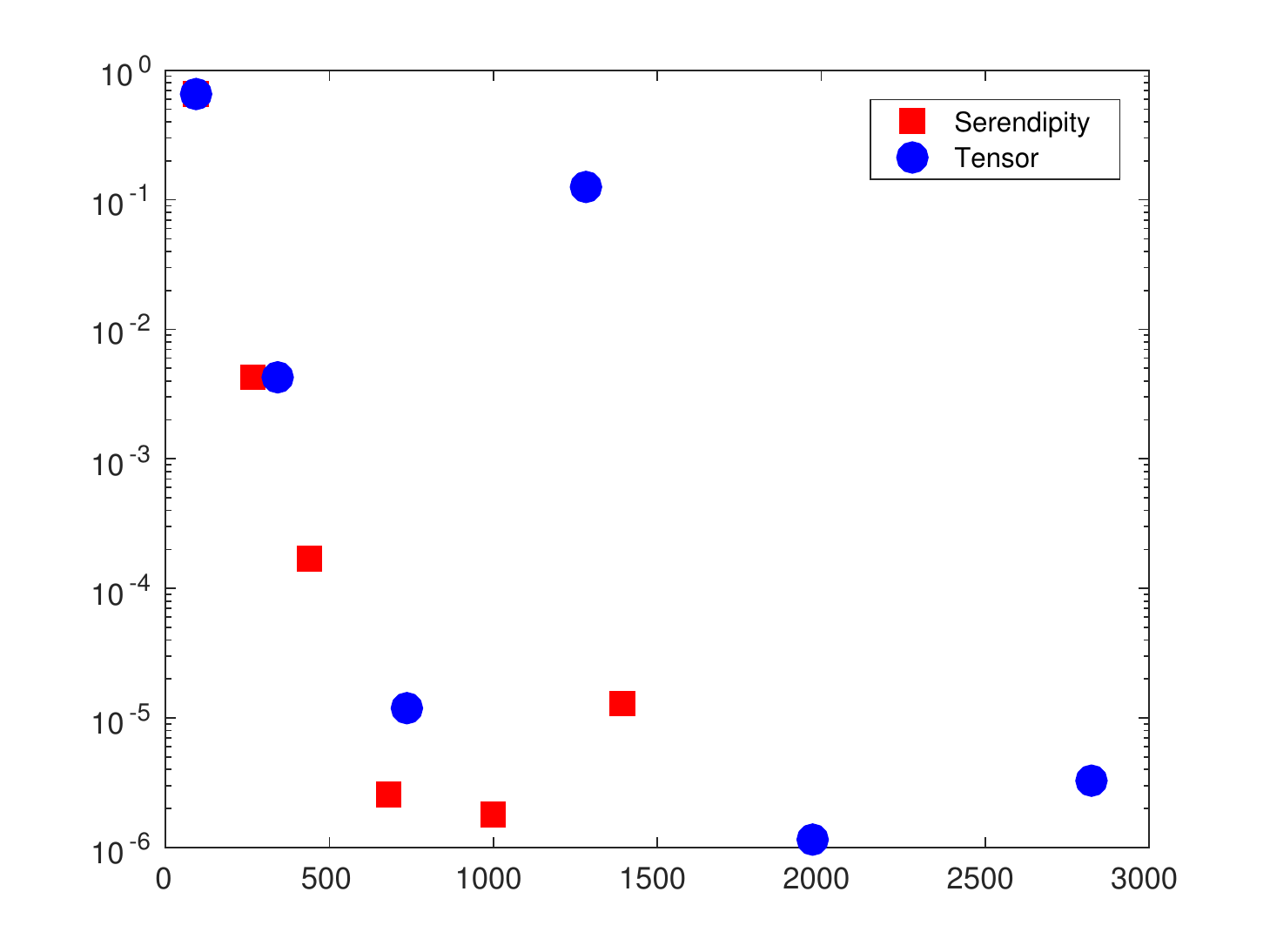} \\[-4mm]
h=1/4 &
h=1/5 
\end{array}$
    \caption{
L-shaped domain, Neumann conditions, $\lambda=2\pi^2$, $p$-refinement experiments.
}
    \label{fig:LNeuTwoPiK}
\end{figure}

\begin{figure}[ht]
$\begin{array}{ccc}
 &
\includegraphics[width=.31\textwidth]{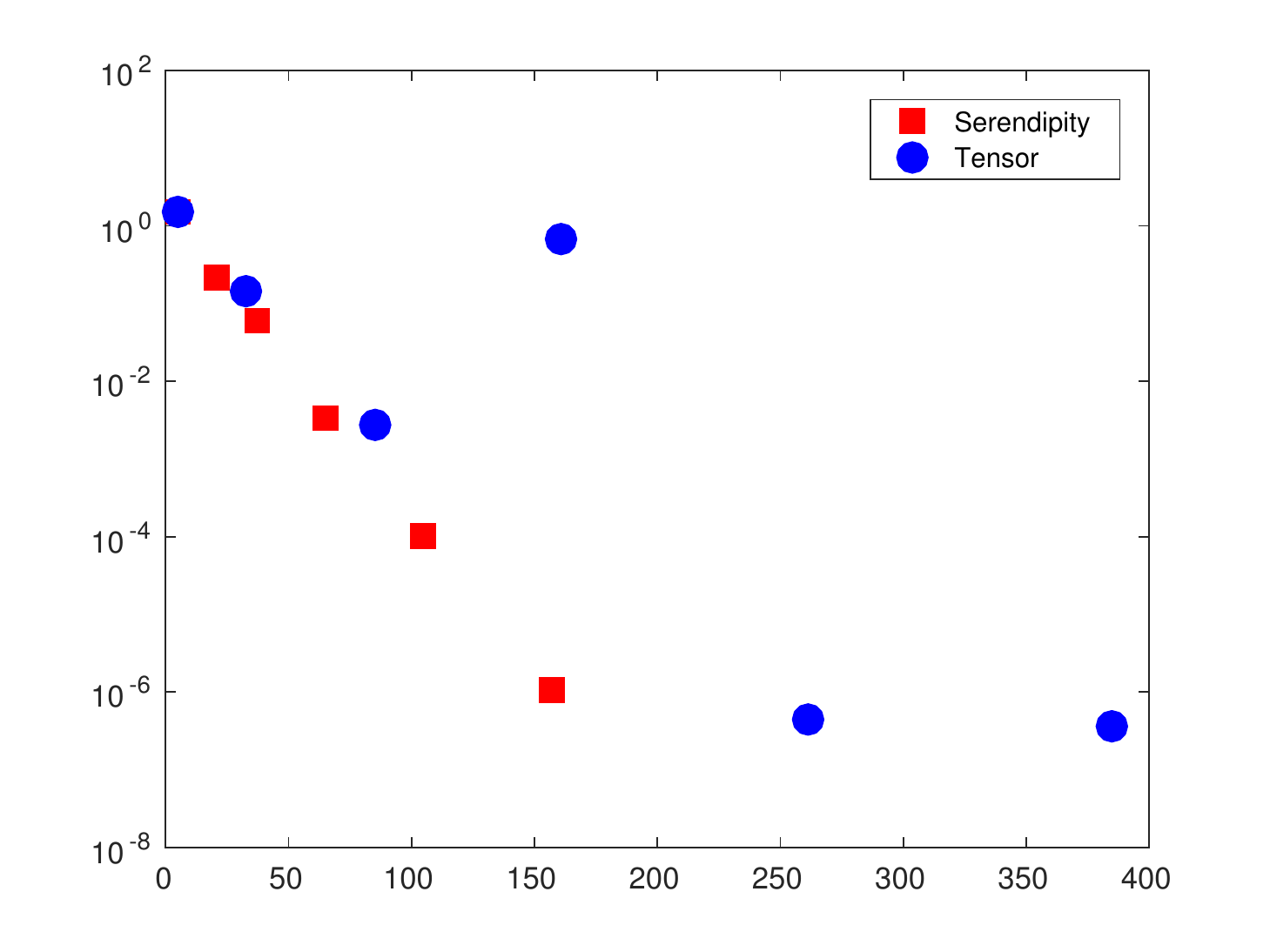} &
\includegraphics[width=.31\textwidth]{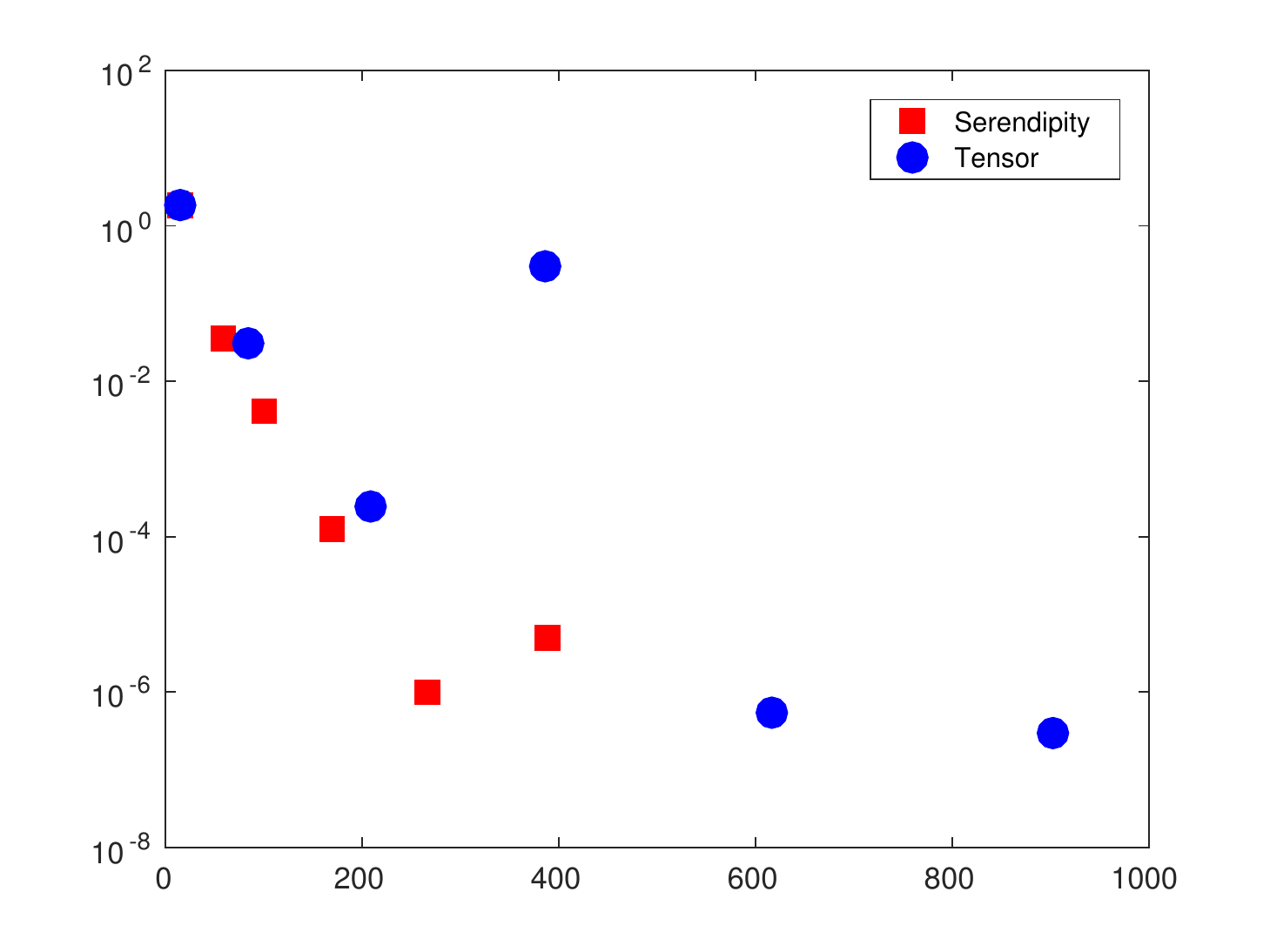}\\[-4mm]
 &
h=1/2 & 
h=1/3  \\[2mm]
\includegraphics[width=.31\textwidth]{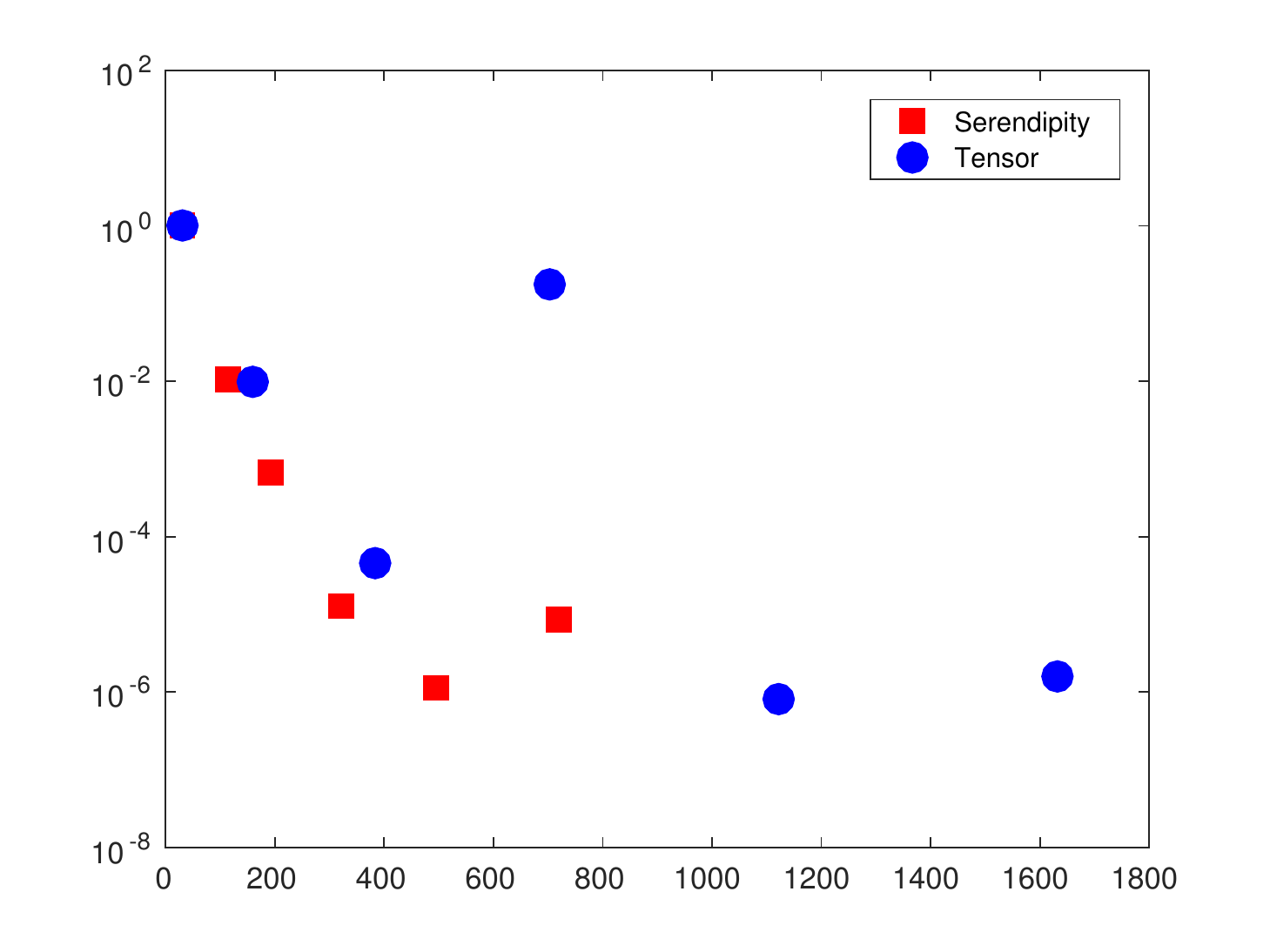} &
\includegraphics[width=.31\textwidth]{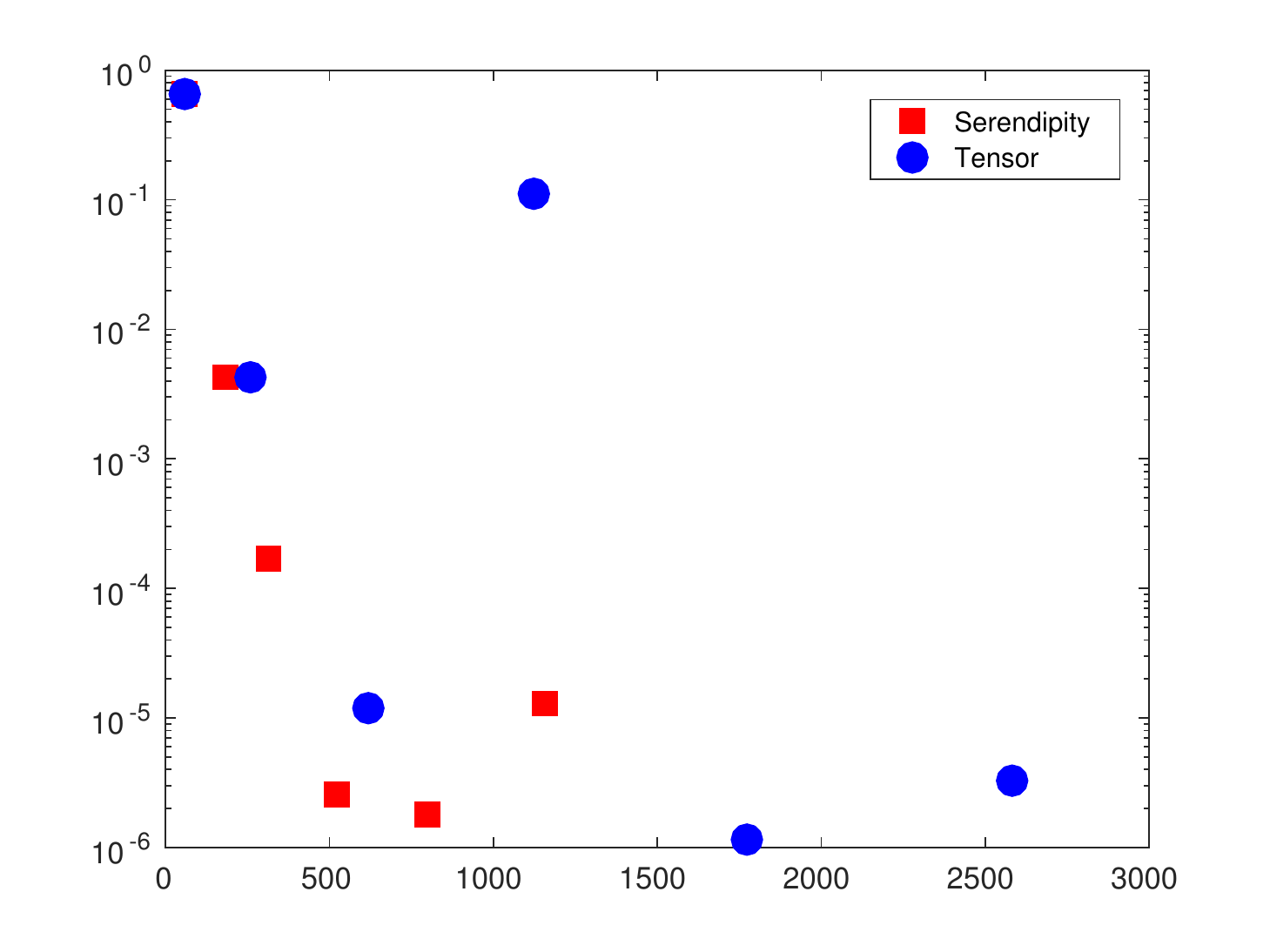} \\[-4mm]
h=1/4 &
h=1/5 
\end{array}$
    \caption{
L-shaped domain, Dirichlet conditions, $\lambda=2\pi^2$, $p$-refinement experiments.
}
    \label{fig:LDirTwoPiK}
\end{figure}

\begin{figure}[ht]
$\begin{array}{ccc}
\includegraphics[width=.31\textwidth]{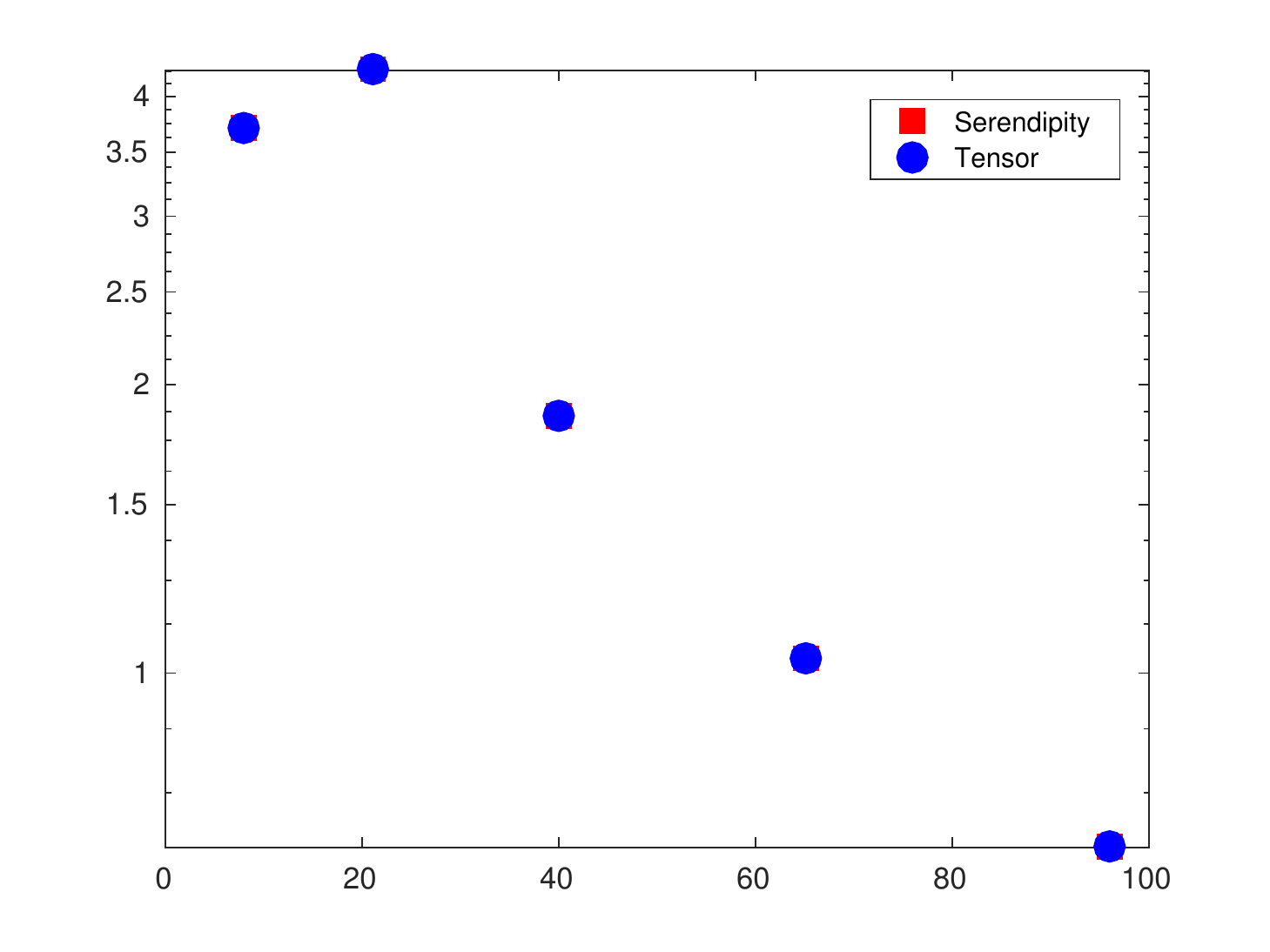} &
\includegraphics[width=.31\textwidth]{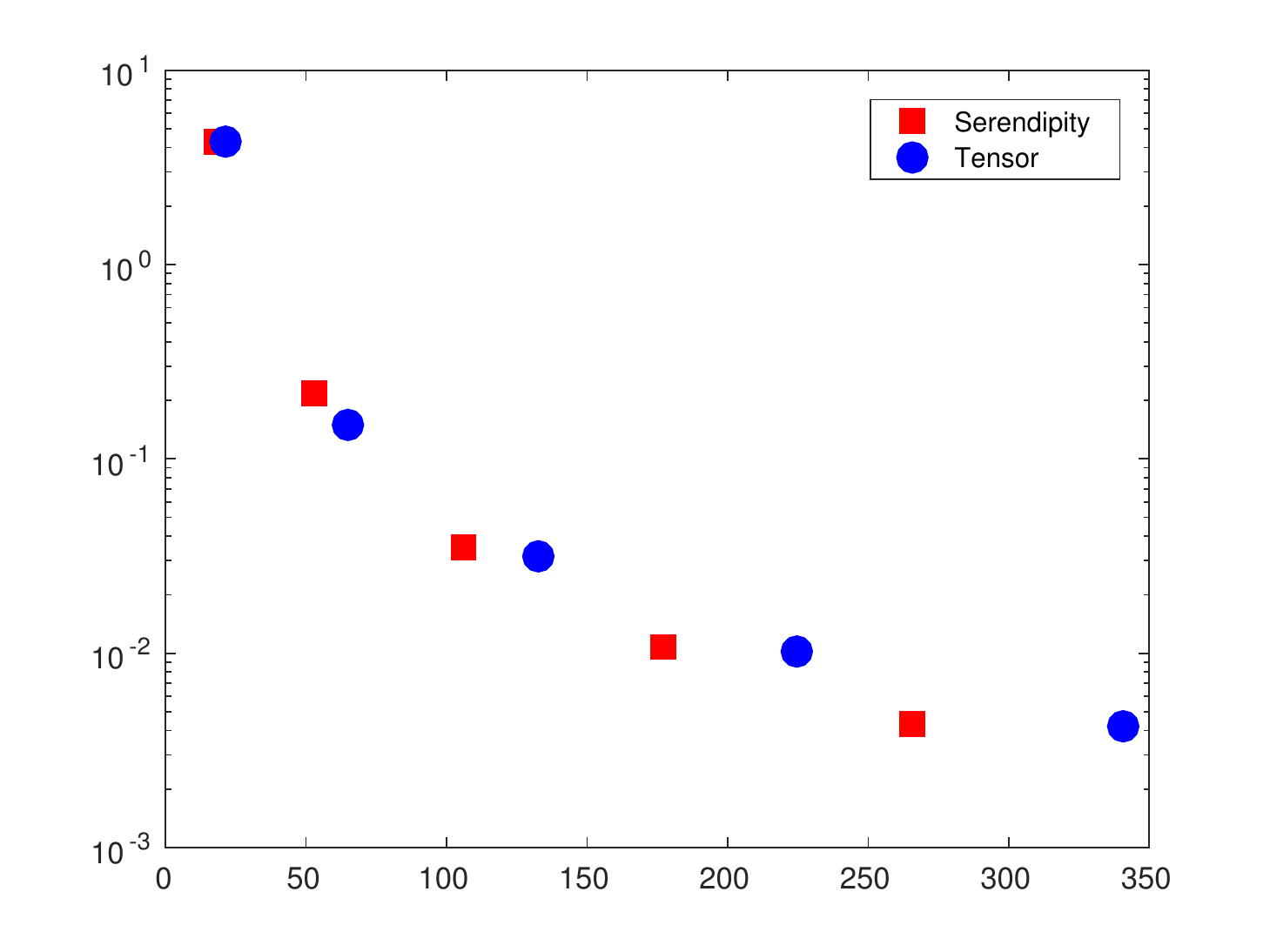} &
\includegraphics[width=.31\textwidth]{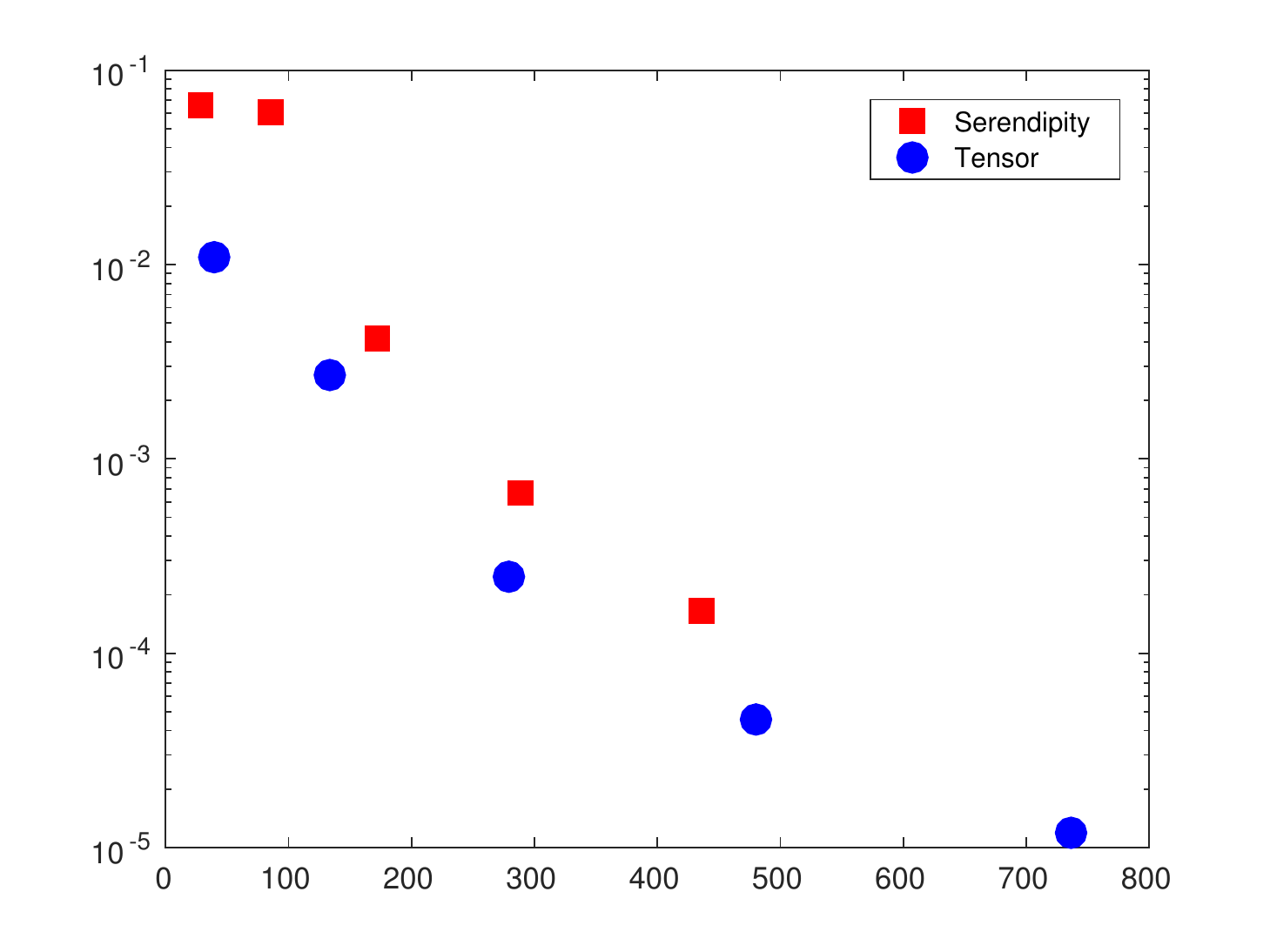} \\[-4mm]
p=1 & 
p=2 &
p=3 \\[2mm]
\includegraphics[width=.31\textwidth]{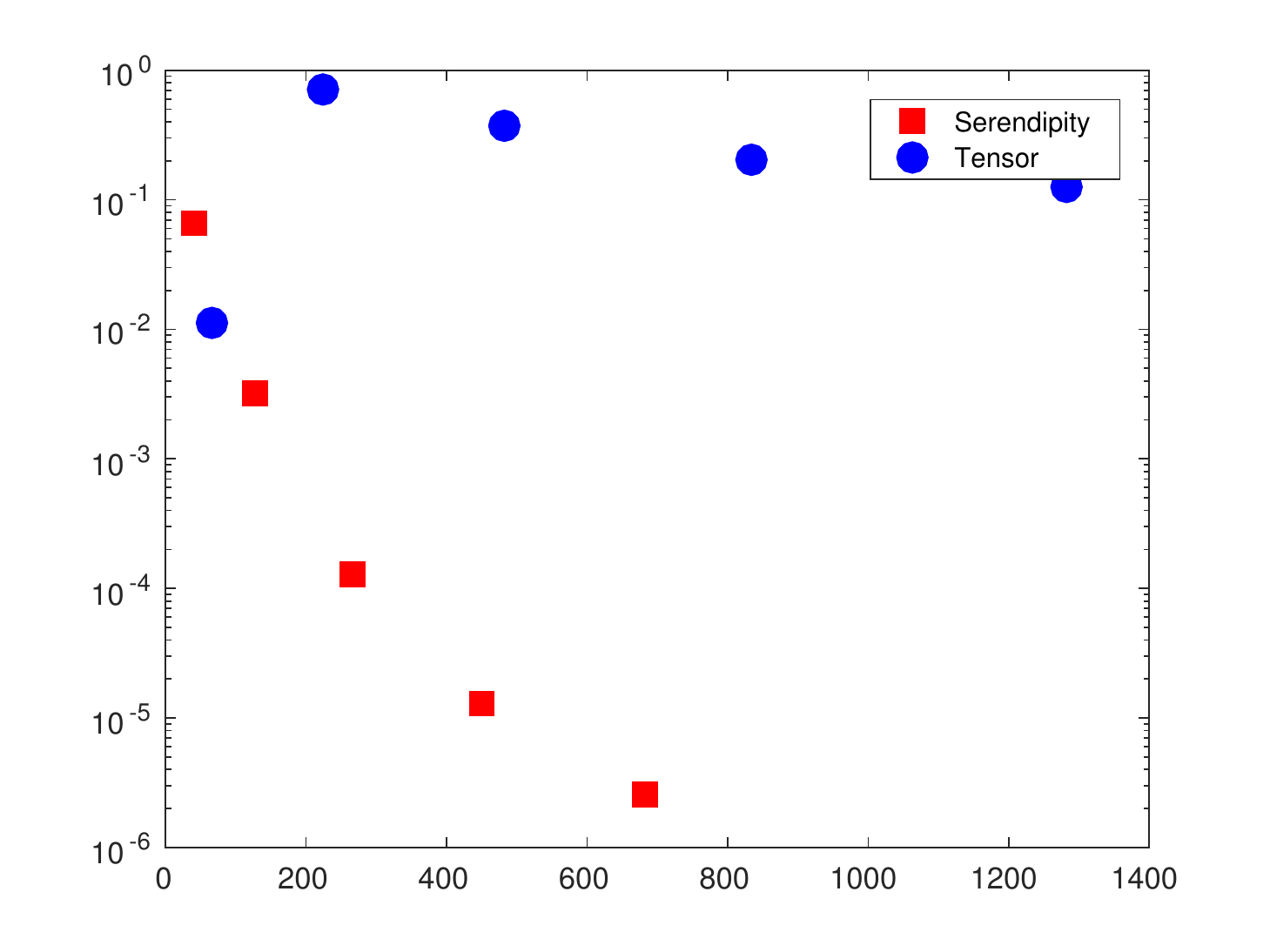} &
\includegraphics[width=.31\textwidth]{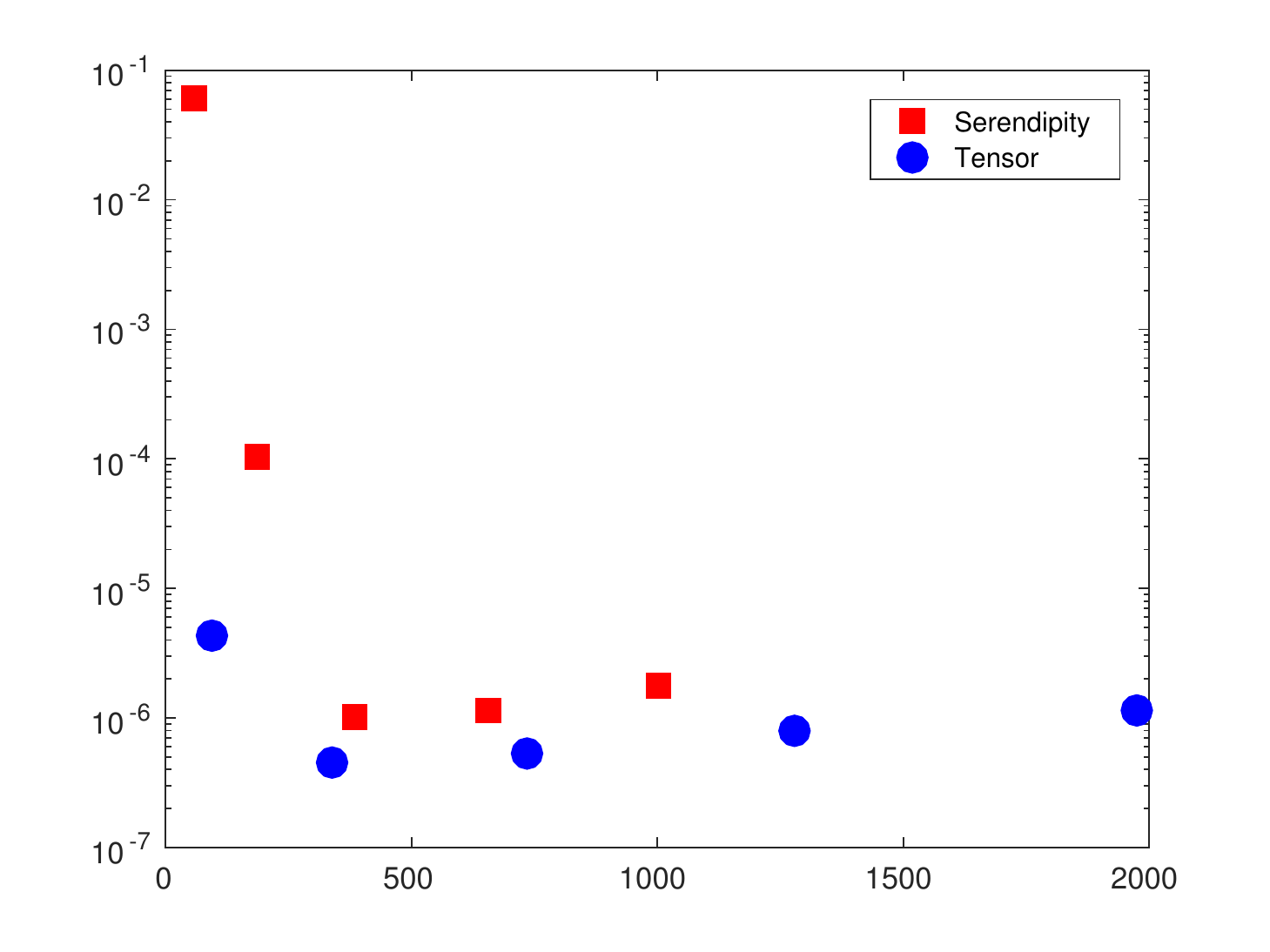} &
\includegraphics[width=.31\textwidth]{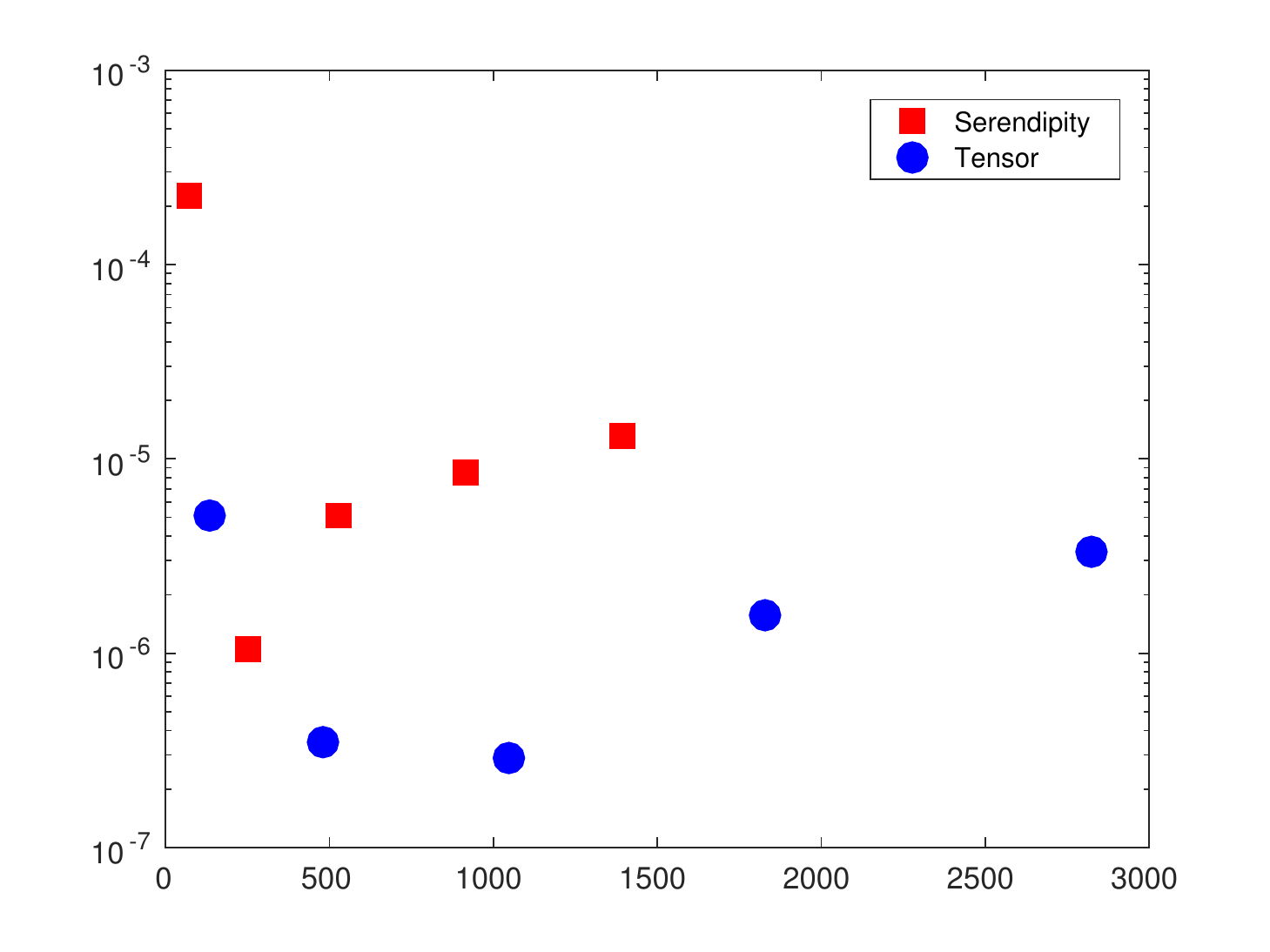} \\[-4mm]
p=4 &
p=5 &
p=6
\end{array}$
    \caption{
L-shaped domain, Neumann conditions, $\lambda=2\pi^2$, $h$-refinement experiments.
}
    \label{fig:LNeuTwoPiP}
\end{figure}

\begin{figure}[ht]
$\begin{array}{ccc}
\includegraphics[width=.31\textwidth]{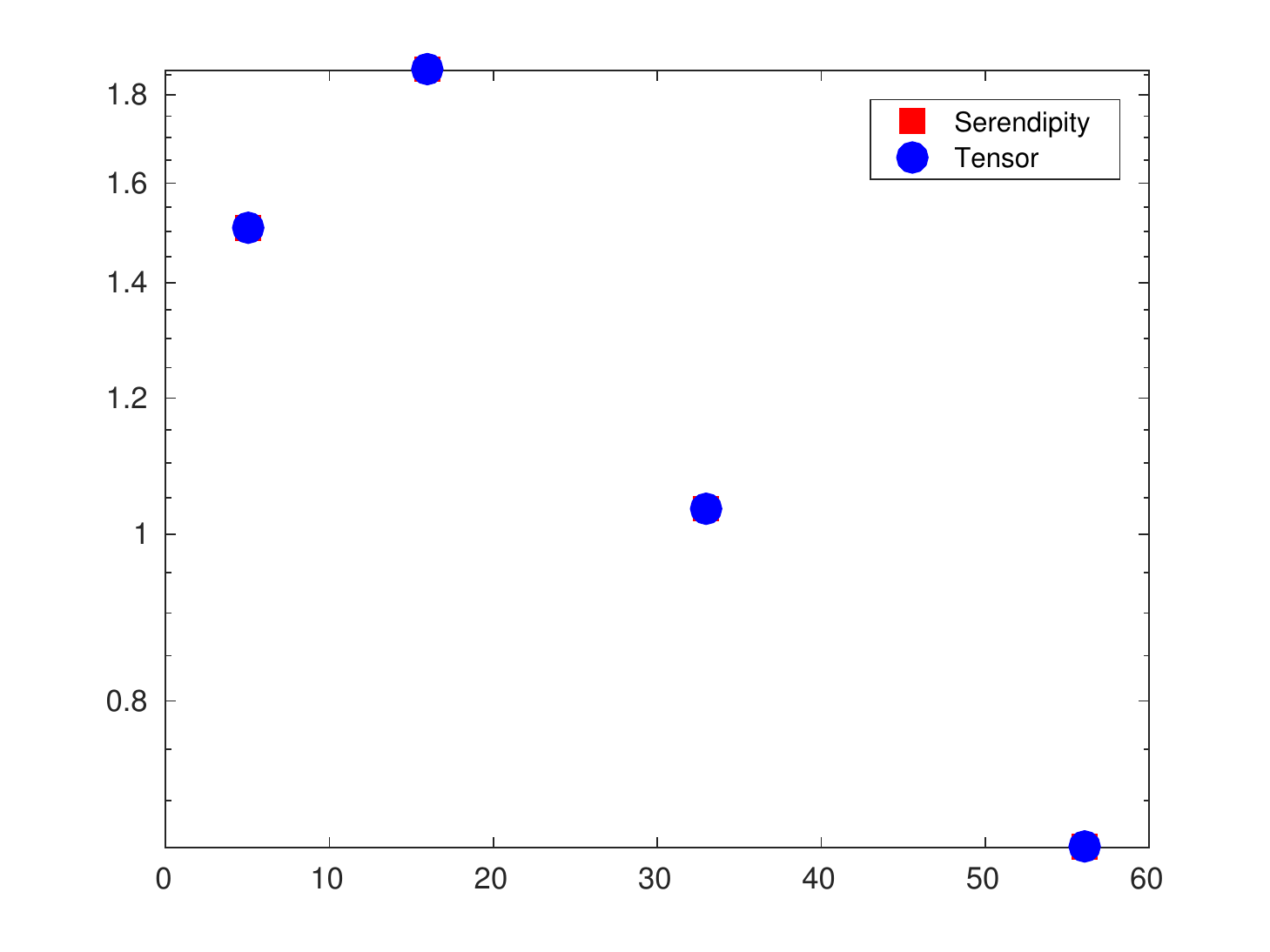} &
\includegraphics[width=.31\textwidth]{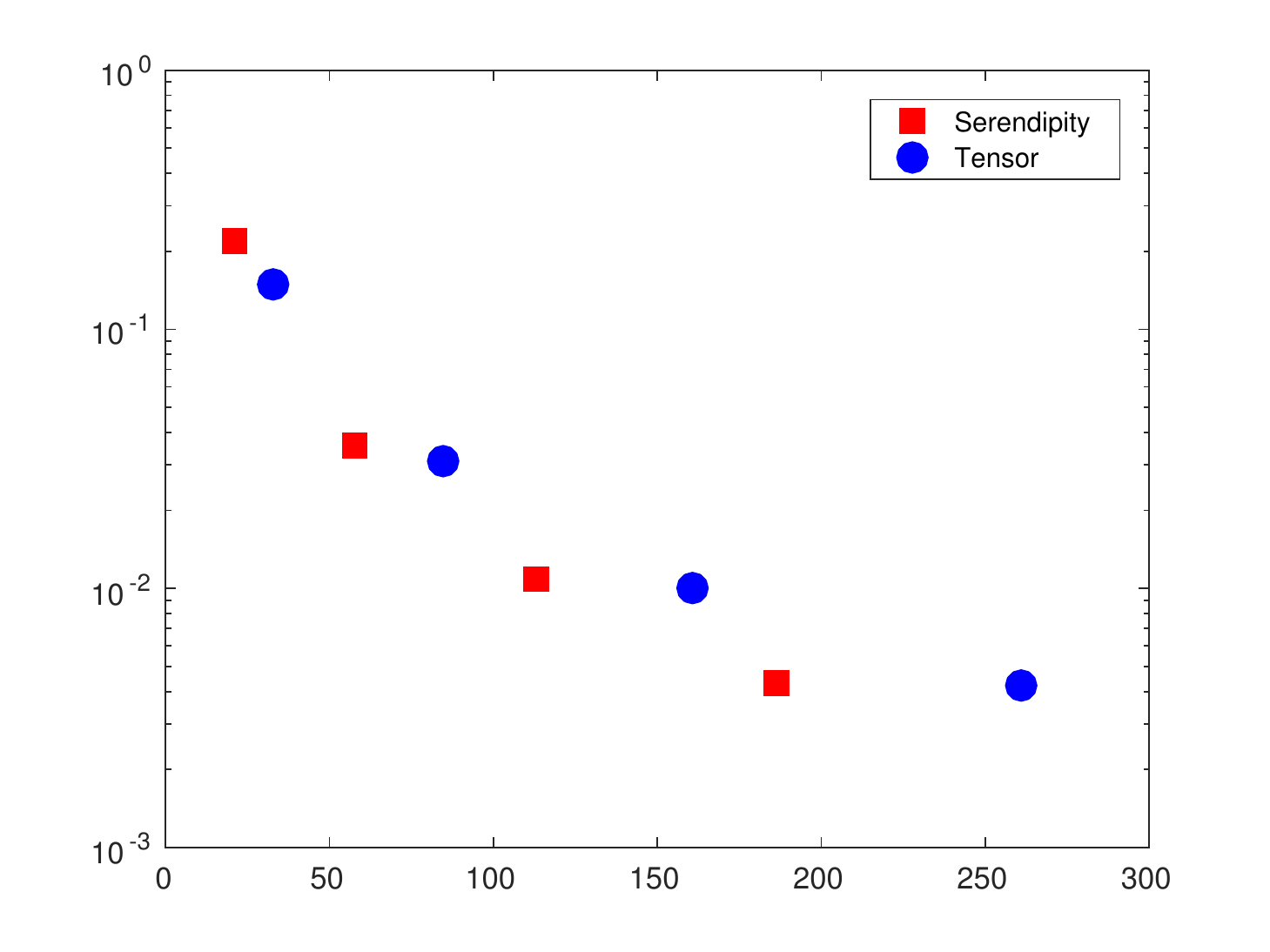} &
\includegraphics[width=.31\textwidth]{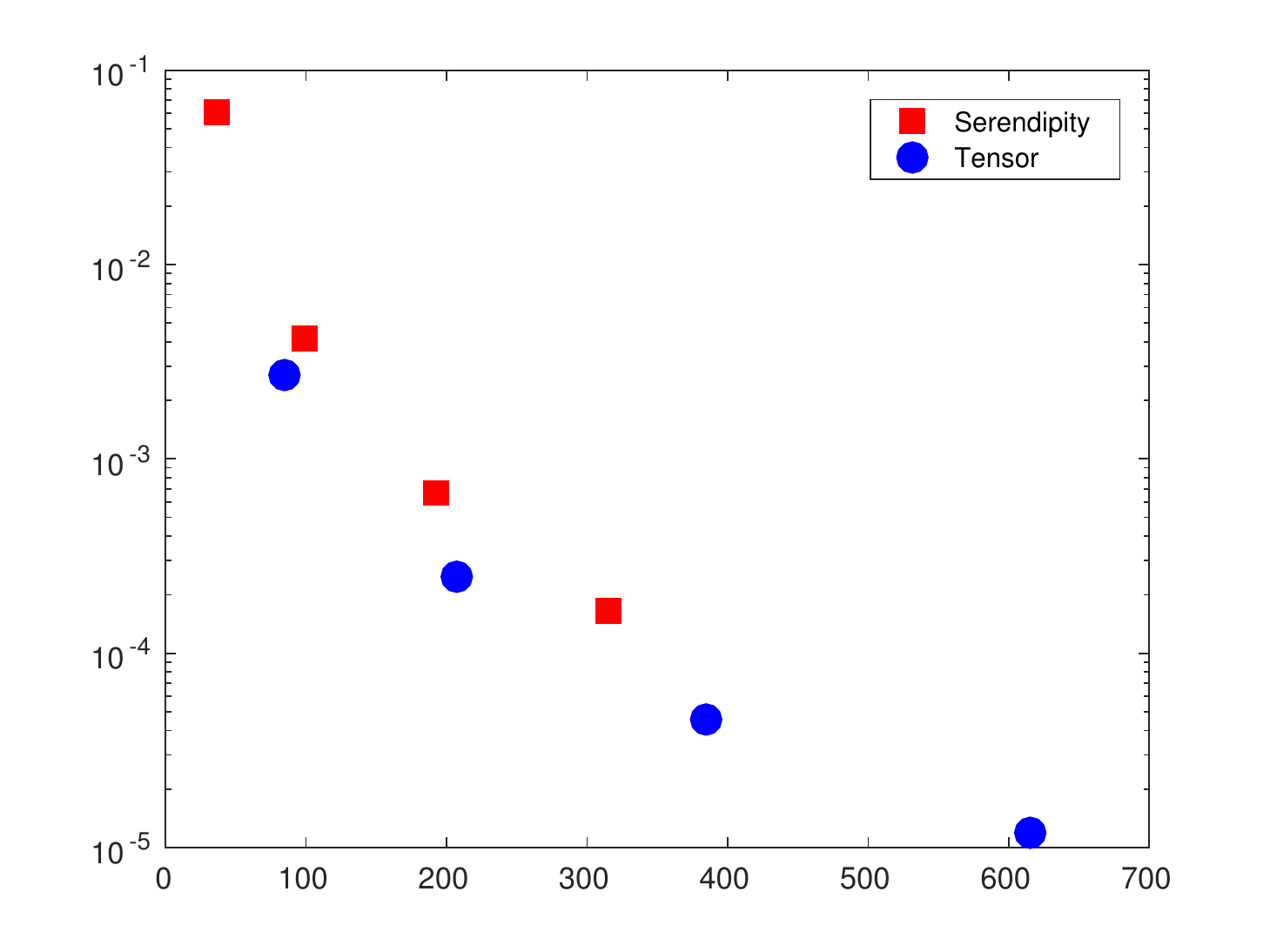} \\[-4mm]
p=1 & 
p=2 &
p=3 \\[2mm]
\includegraphics[width=.31\textwidth]{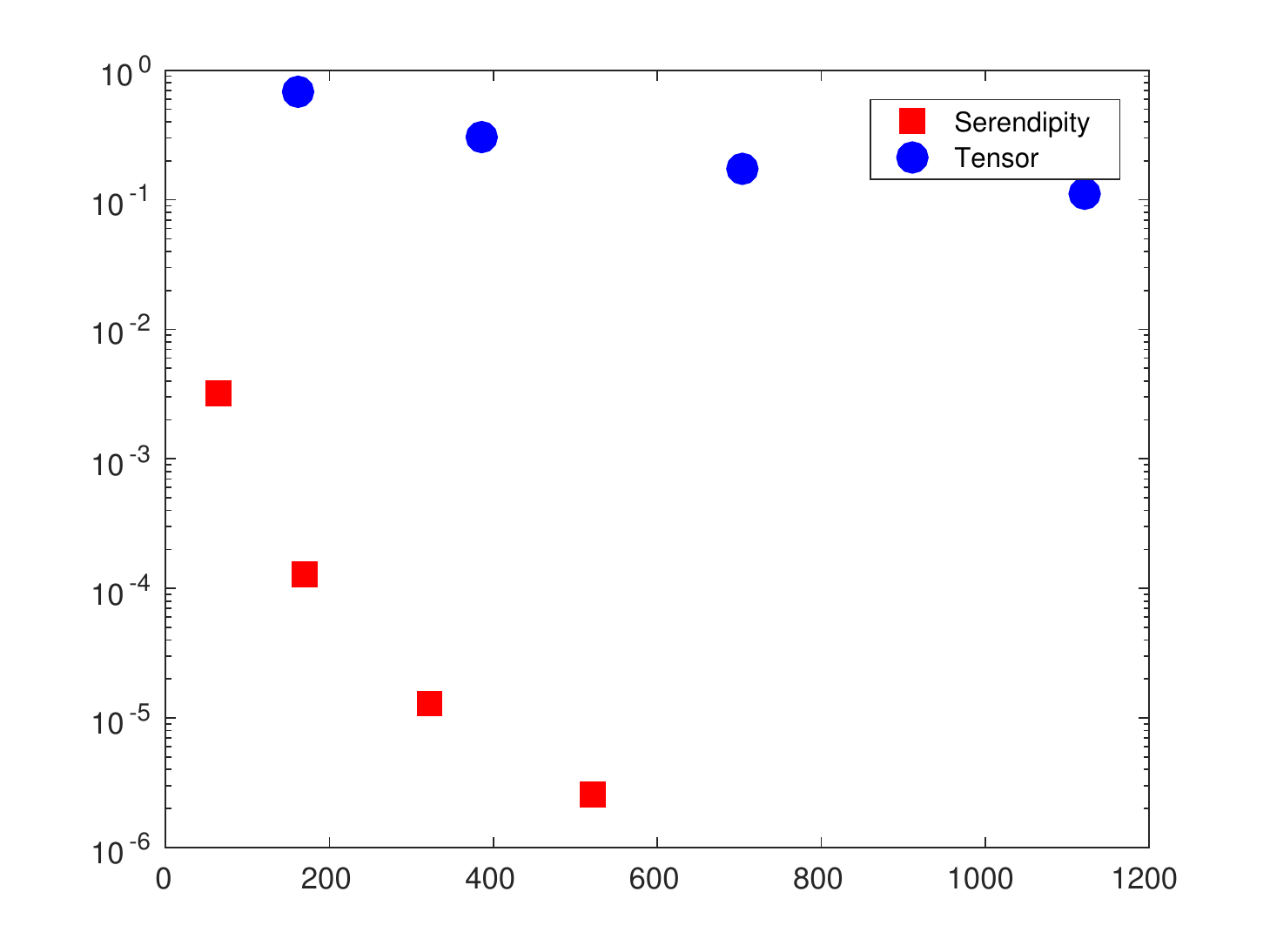} &
\includegraphics[width=.31\textwidth]{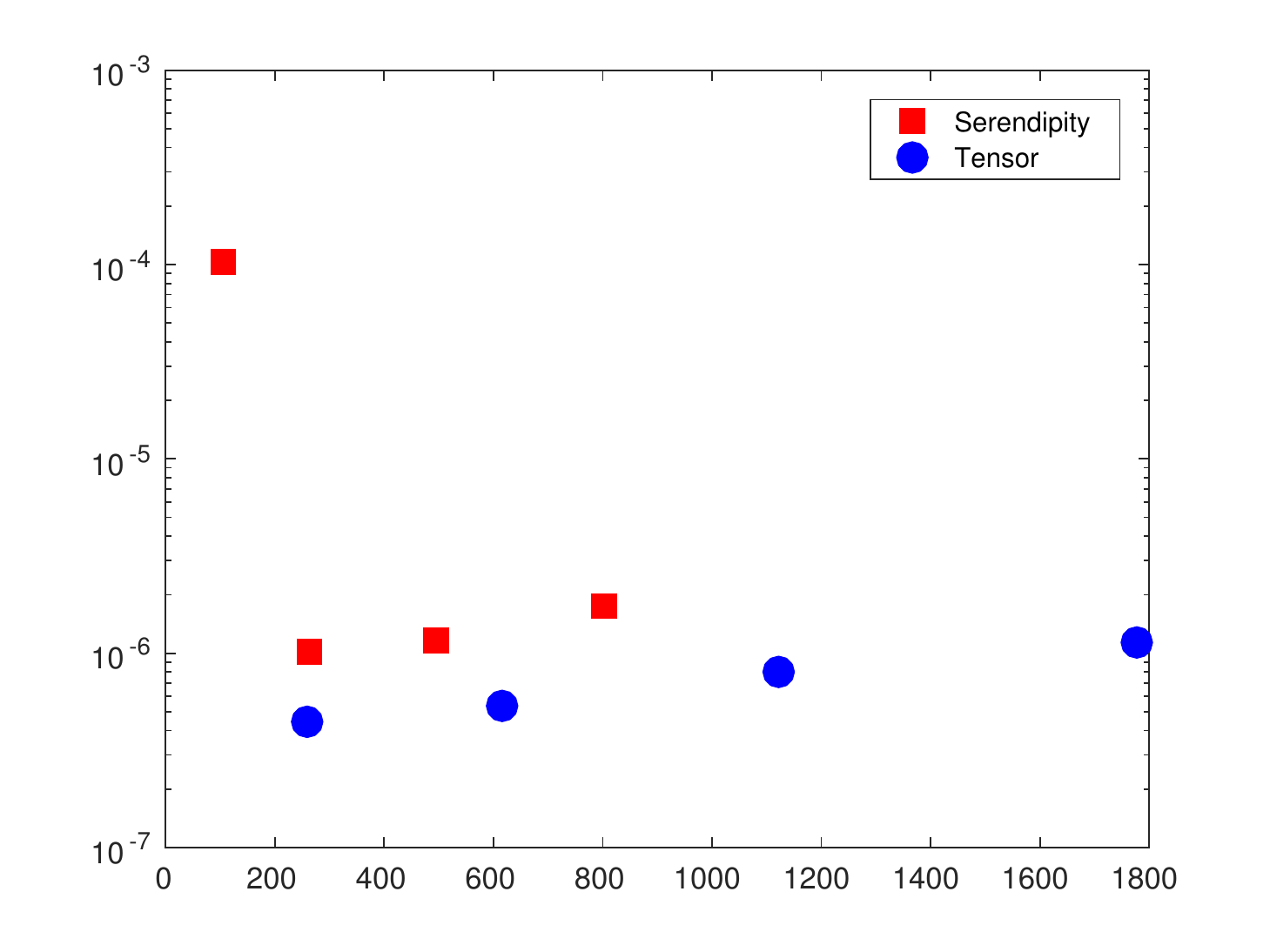} &
\includegraphics[width=.31\textwidth]{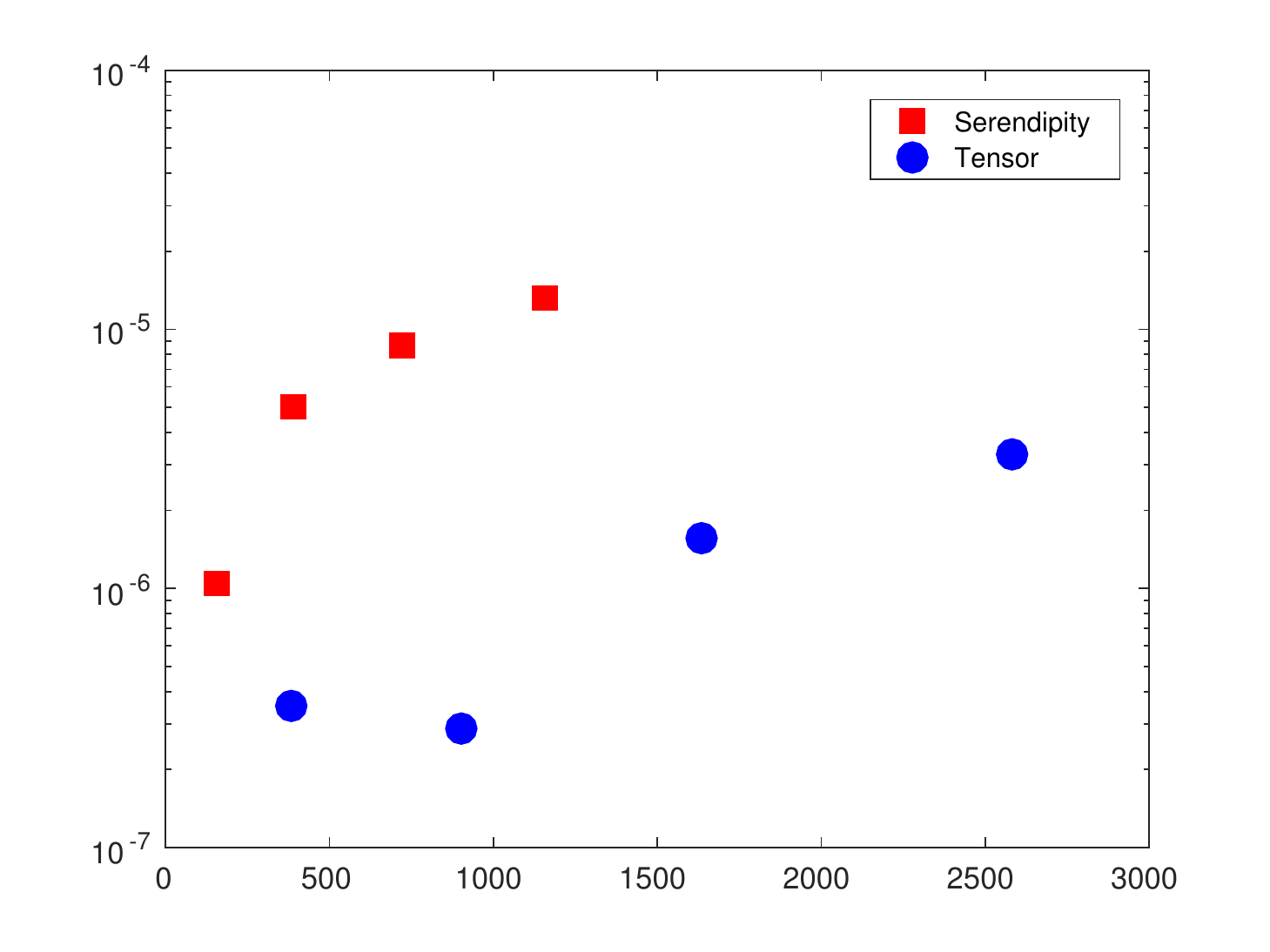} \\[-4mm]
p=4 &
p=5 &
p=6
\end{array}$
    \caption{
L-shaped domain, Dirichlet conditions, $\lambda=2\pi^2$, $h$-refinement experiments.
}
    \label{fig:LDirTwoPiP}
\end{figure}

\begin{figure}[ht]
$\begin{array}{ccc}
\includegraphics[width=.31\textwidth]{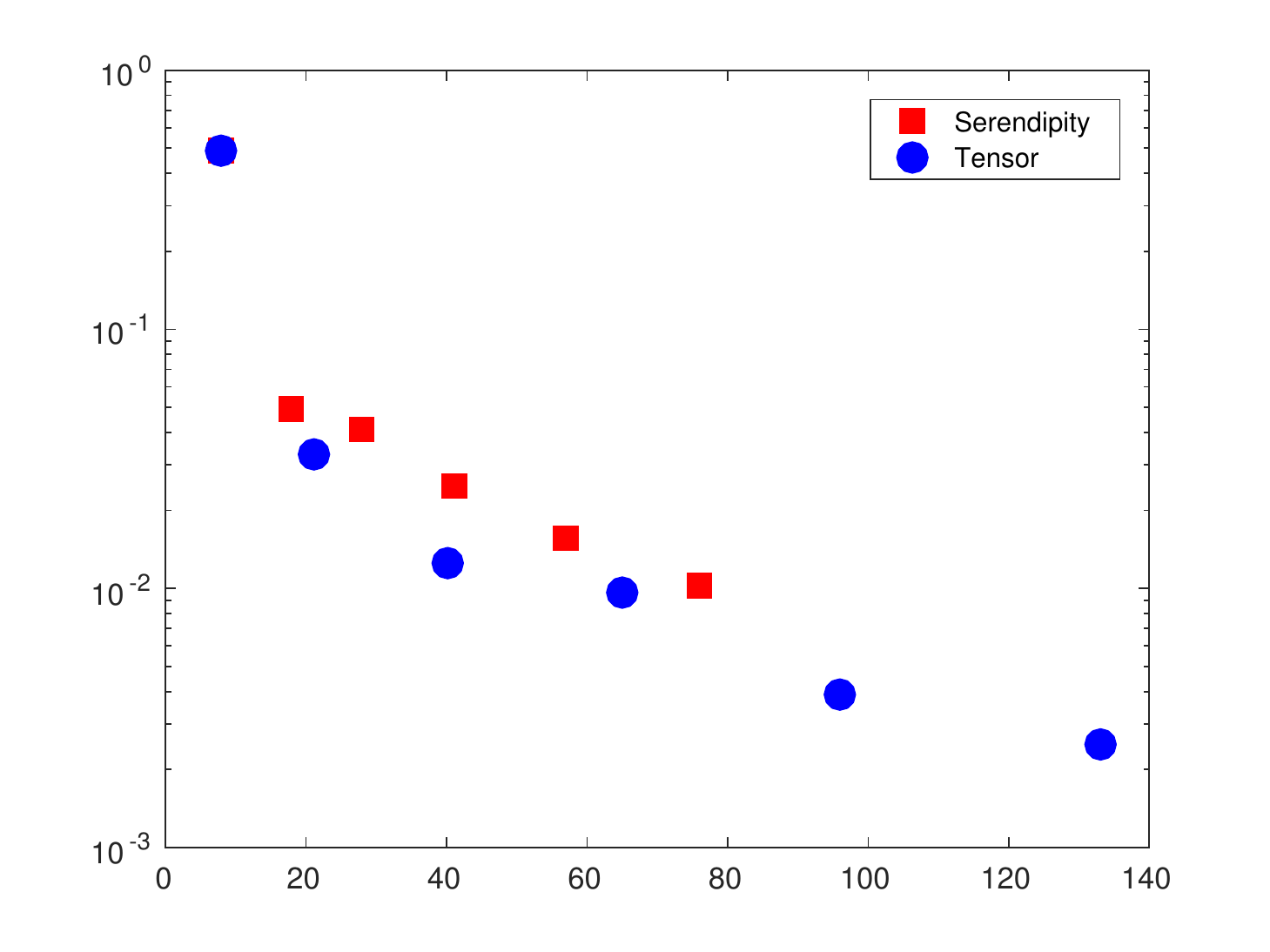} &
\includegraphics[width=.31\textwidth]{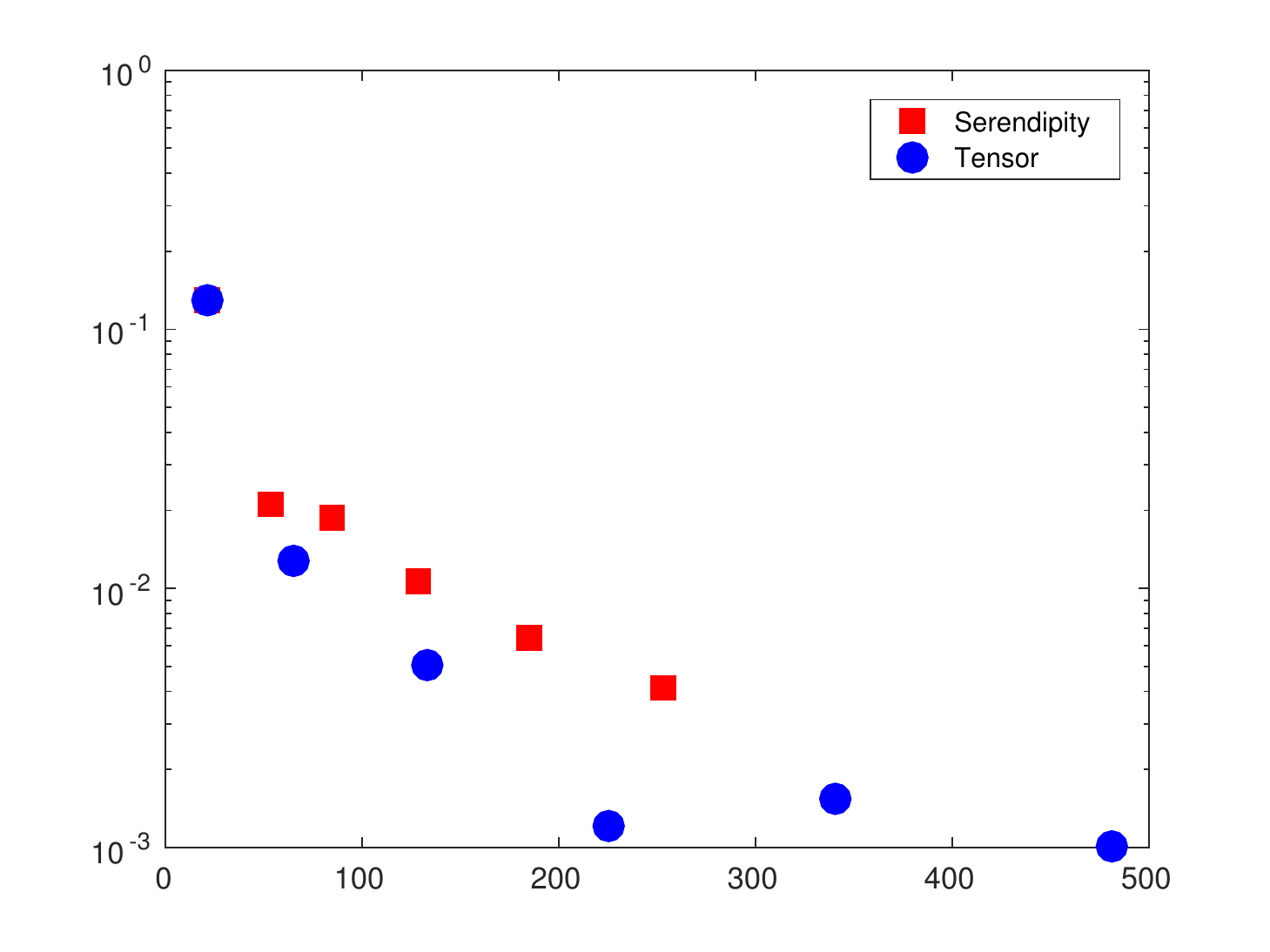} &
\includegraphics[width=.31\textwidth]{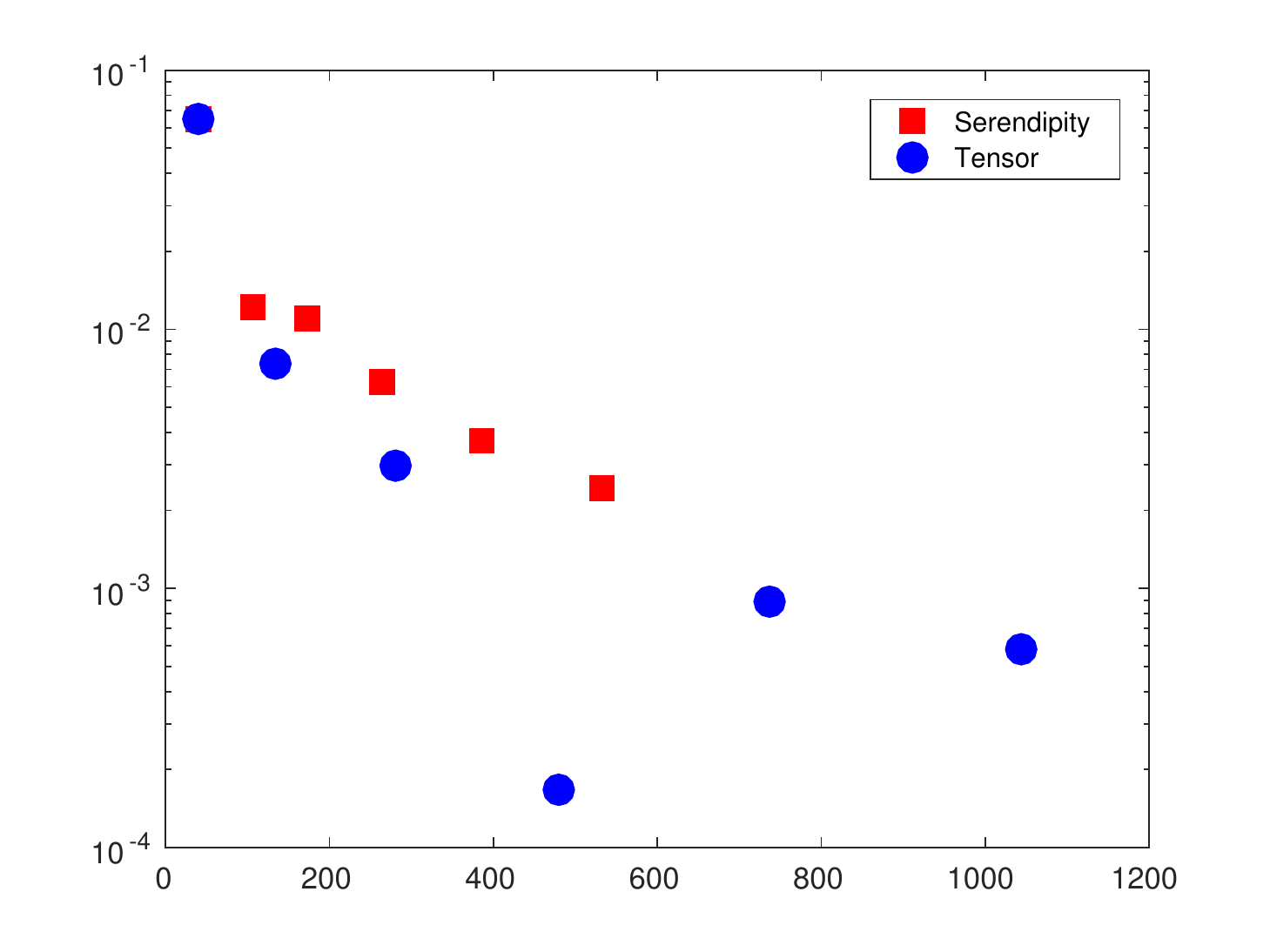} \\[-4mm]
h=1 & 
h=1/2 &
h=1/3 \\[2mm]
\includegraphics[width=.31\textwidth]{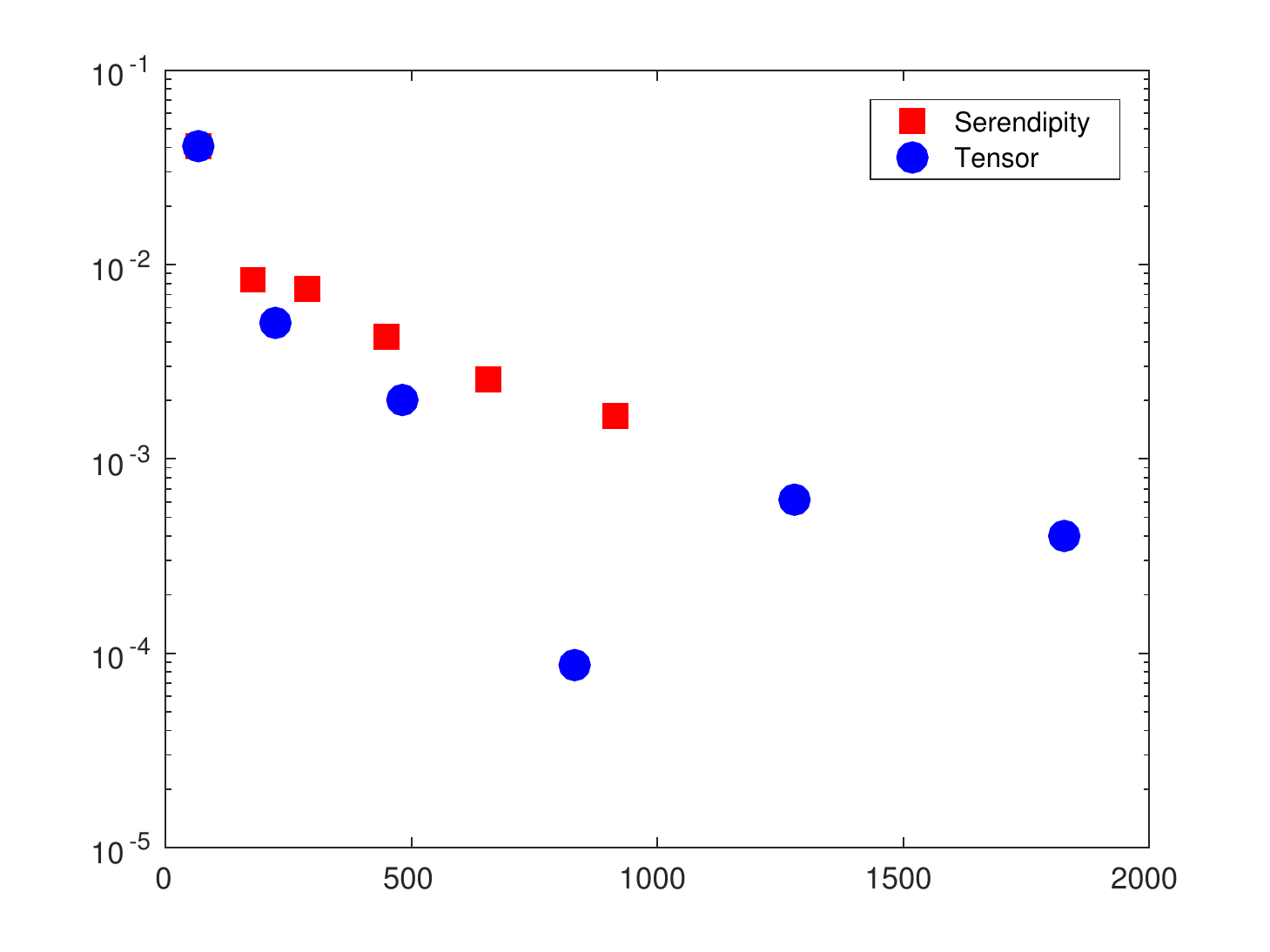} &
\includegraphics[width=.31\textwidth]{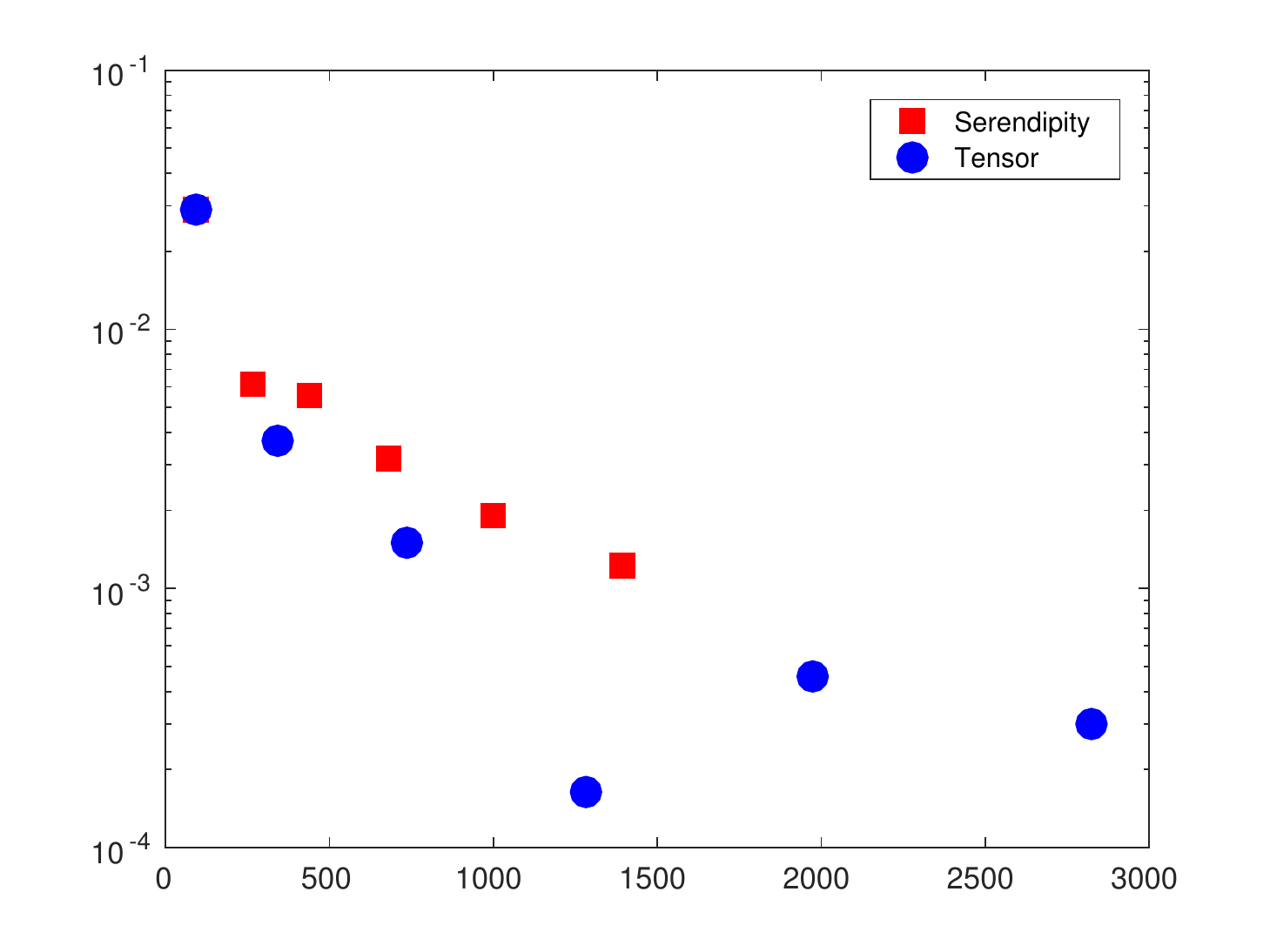} \\[-4mm]
h=1/4 &
h=1/5 
\end{array}$
\caption{
L-shaped domain, Neumann conditions, $\lambda=1.4756218450$, $p$-refinement experiments.
}
    \label{fig:LNeuOneK}
\end{figure}

\begin{figure}[ht]
$\begin{array}{ccc}
\includegraphics[width=.31\textwidth]{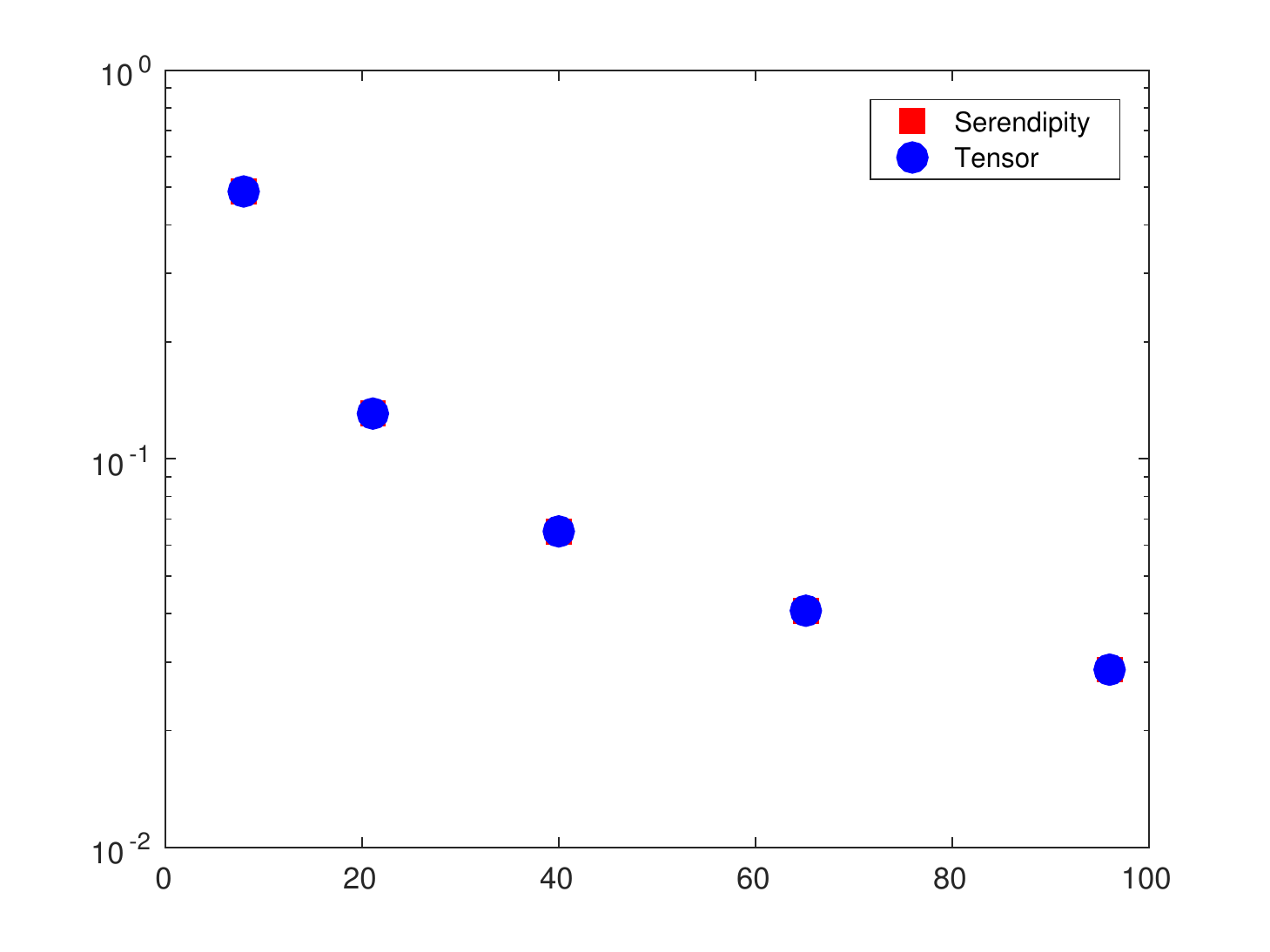} &
\includegraphics[width=.31\textwidth]{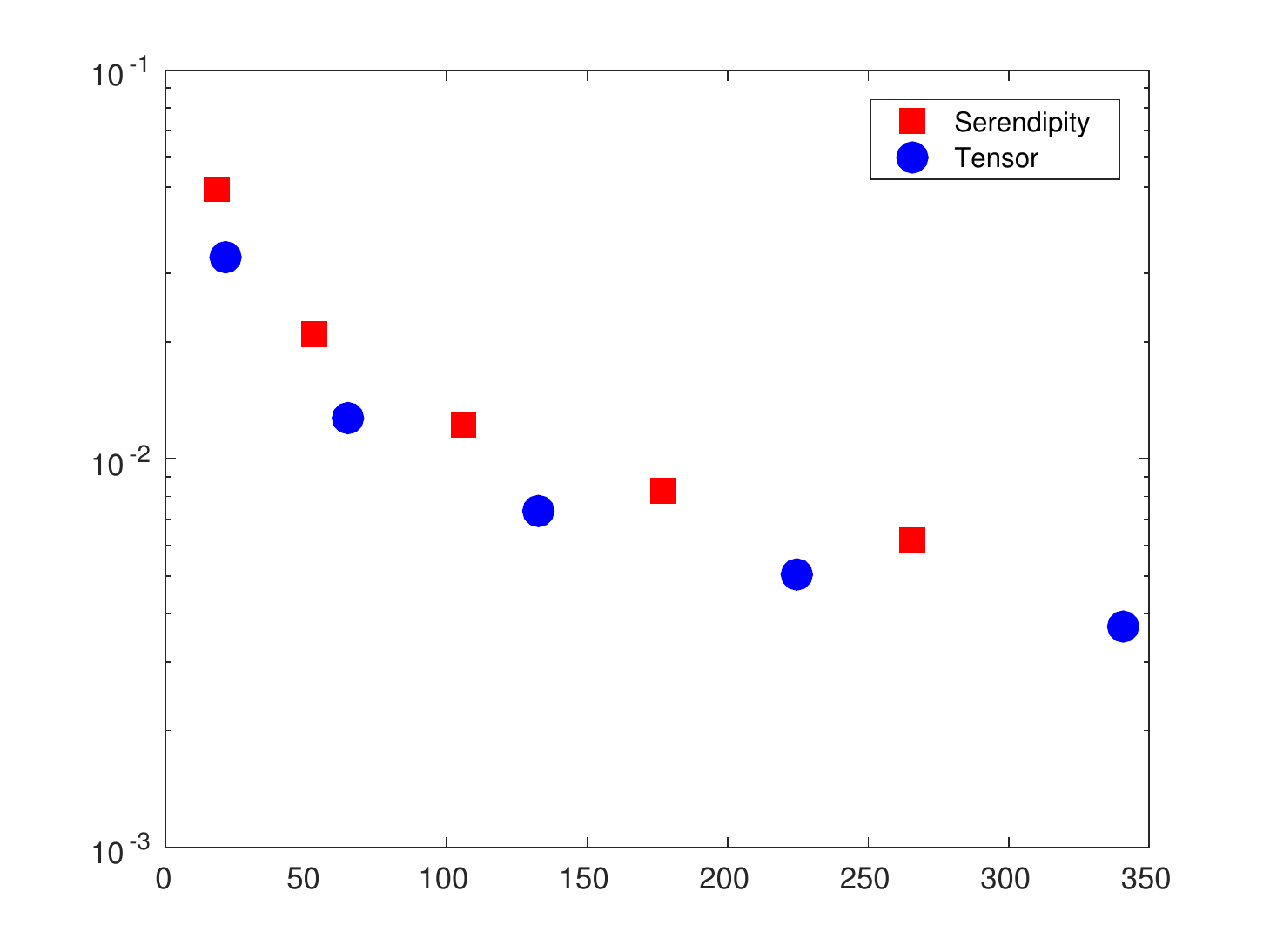} &
\includegraphics[width=.31\textwidth]{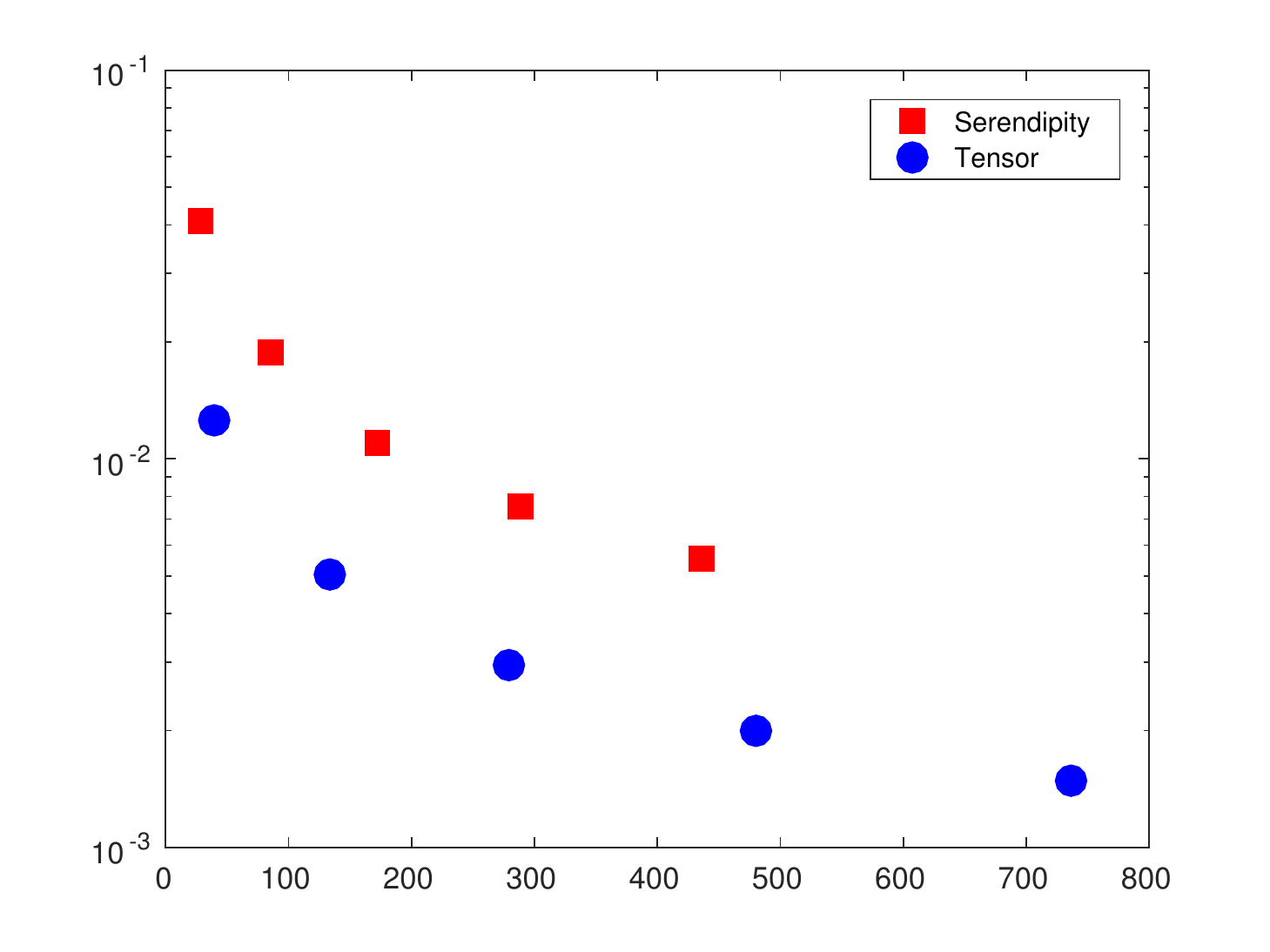} \\[-4mm]
p=1 & 
p=2 &
p=3 \\[2mm]
\includegraphics[width=.31\textwidth]{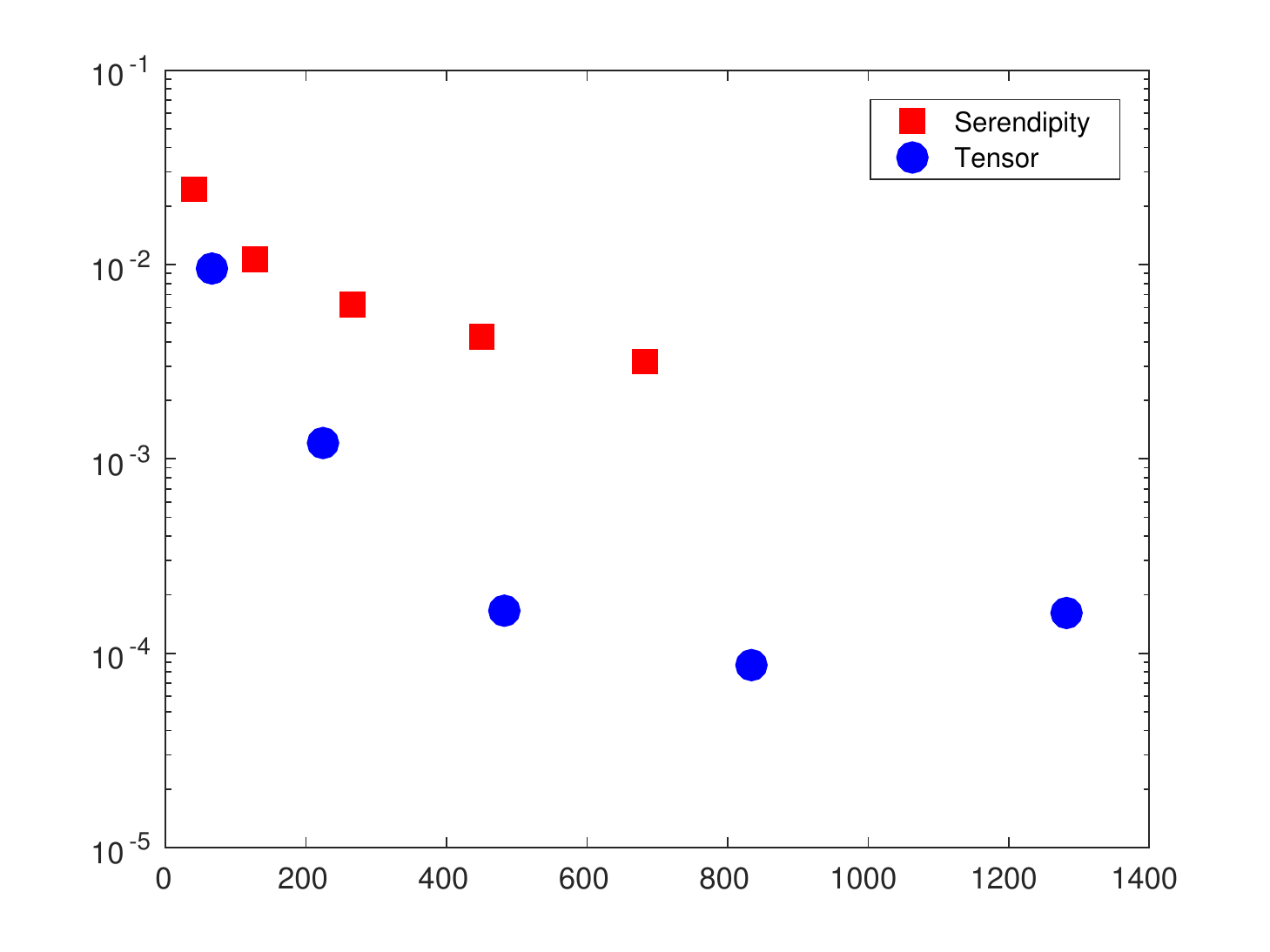} &
\includegraphics[width=.31\textwidth]{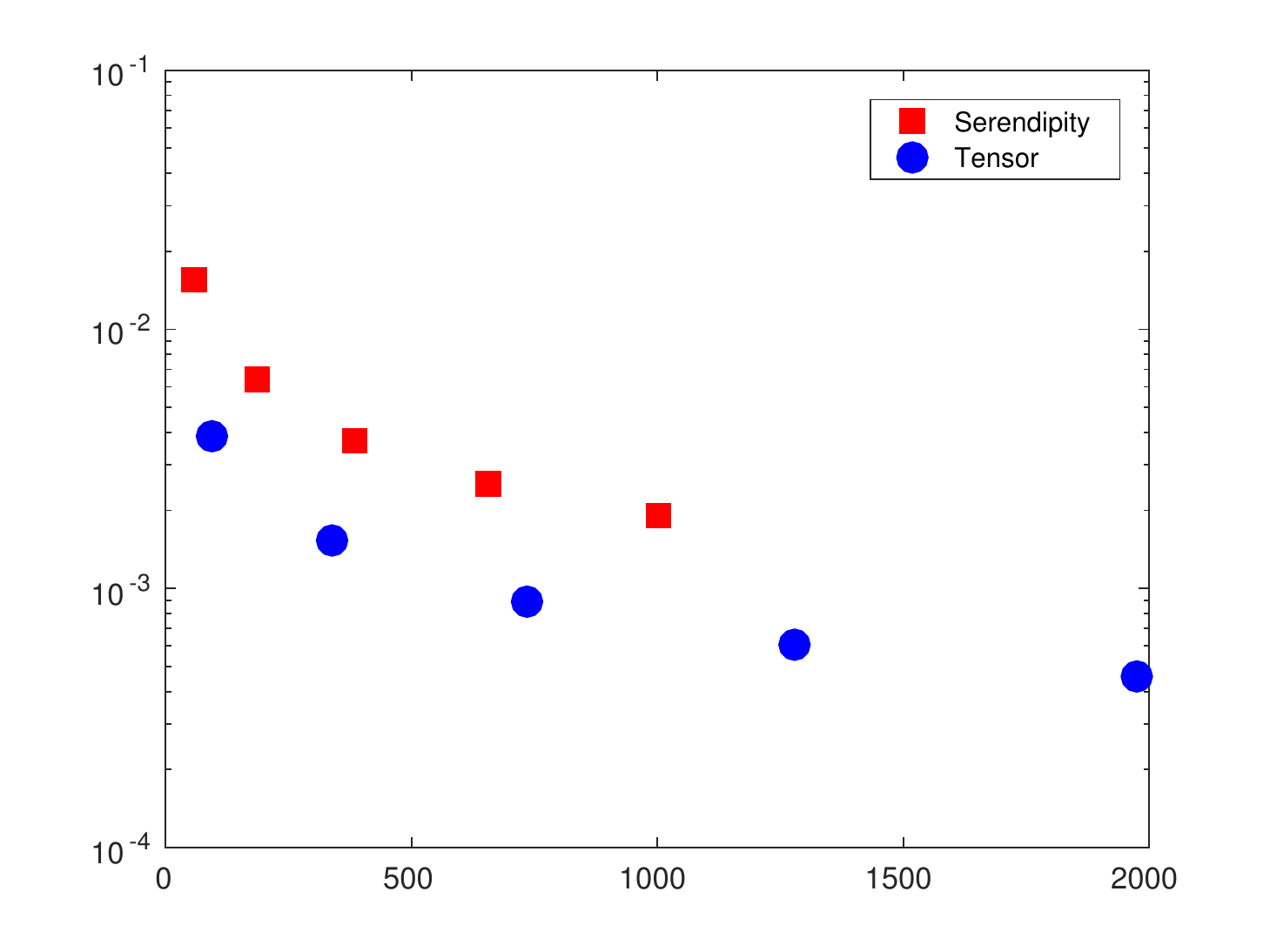} &
\includegraphics[width=.31\textwidth]{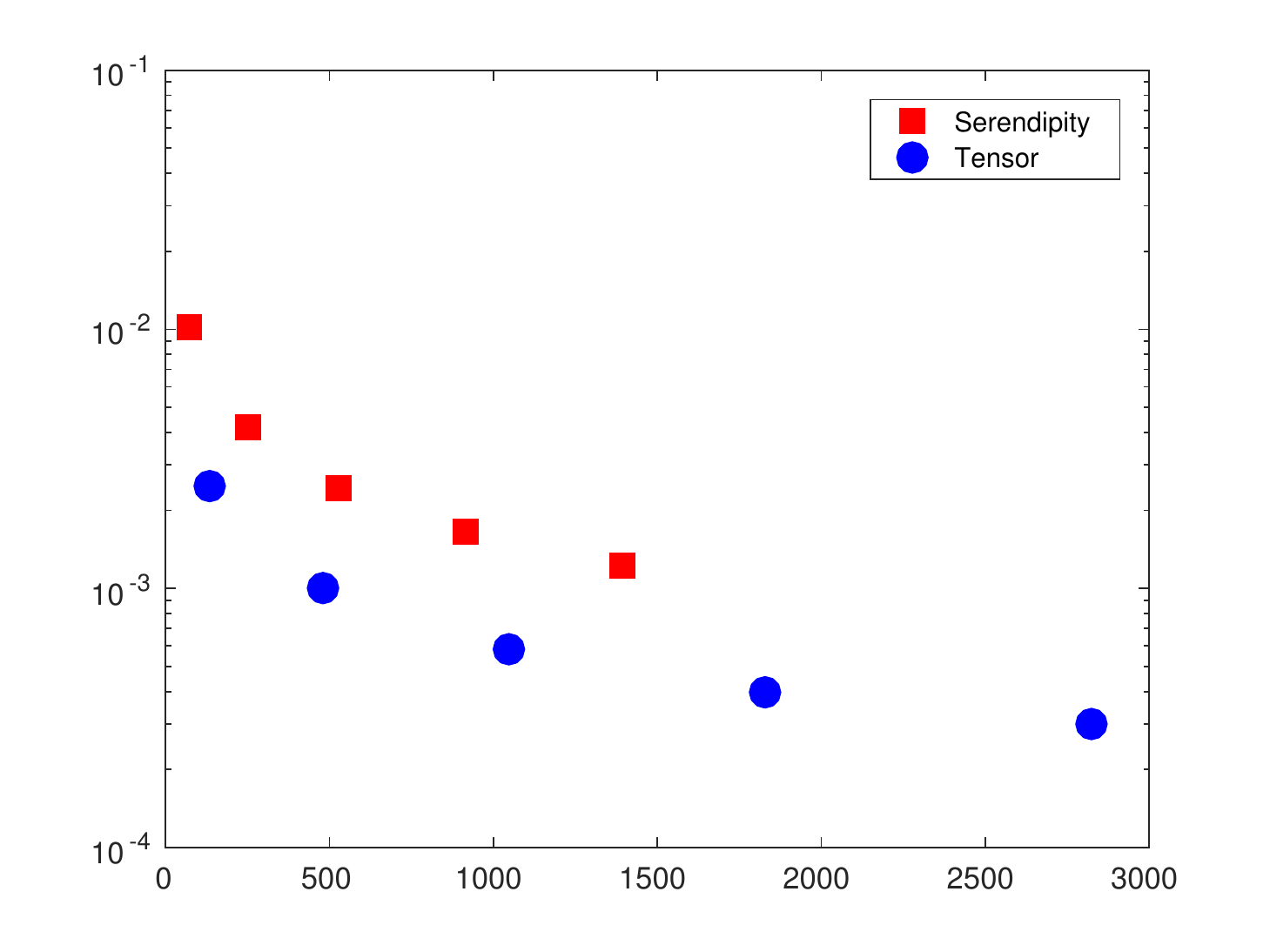} \\[-4mm]
p=4 &
p=5 &
p=6
\end{array}$
    \caption{
L-shaped domain, Neumann conditions, $\lambda=1.4756218450$, $h$-refinement experiments.
}
    \label{fig:LNeuOneP}
\end{figure}

Our numerical experiments are characterized by four choices: domain (square or L-shaped), boundary conditions (Dirichlet or Neumann), eigenvalue $\lambda$ being approximated, and refinement strategy ($p$-refine with $h$ fixed or $h$-refine with $p$ fixed).
For each choice, we report the error in the numerical approximation of $\lambda$ as a function of the number of global degrees of freedom, i.e.\ the dimension of $V_{h,p}$.
Two data series are generated in this fashion: one for tensor product elements and one for serendipity elements, using $p=1$ through 6 for a fixed $h$ value ($p$-refinement), or for $h=1$, $1/2$, $\ldots$, $1/5$ for a fixed $p$ value ($h$-refinement).
The results of these experiments are shown in Figures~\ref{fig:NeuTwoPiK}-\ref{fig:LNeuOneP}, collected at the end of this section, and are discussed below.

\paragraph*{Square Domain}
Our first comparison of the tensor product and serendipity elements is on a square domain with Neumann boundary conditions.
Figure~\ref{fig:NeuTwoPiK} shows the error in approximating the eigenvalue $2\pi^2$ when we fix $h$ and allow $p$ to vary.
Ignoring for now the outlier corresponding to one of the tensor product solutions, we see that in nearly every case, using serendipity elements can match the accuracy of the eigenvalue obtained by tensor product elements with much fewer degrees of freedom. 
For example, in Figure~\ref{fig:NeuTwoPiK}: $h=1/4$, we see that we can obtain an approximate eigenvalue which differs by about $10^{-6}$ from the theoretical using both serendipity and tensor elements. 
However, when using the serendipity elements, we see a reduction in number of degrees of freedom by approximately half compared to the tensor product element.

We see similar behavior in the Dirichlet problem, depicted in Figure~\ref{fig:DirTwoPiK}, with the obvious difference of an overall reduction in the number of degrees of freedom. 
Note that since we remove the degrees of freedom corresponding to the boundary, not discretizing the mesh at all (i.e. when $h = 1$) results in having too few equations to properly solve for a nonzero eigenvalue for small $p$.

When we consider the Neumann problem with $p$ fixed and $h$ varied, we see the results depicted in Figure~\ref{fig:NeuTwoPiP}. 
In contrast to the previously discussed results, we see that, in nearly every case, the tensor product elements achieve better accuracy than serendipity while using fewer degrees of freedom. 
The only exception is when $p = 4$, also depicted in Figure~\ref{fig:NeuTwoPiP}. 
Here, we note a large increase in error when using tensor product elements. 
This effect can be seen in nearly every plot for $h$-refinements and accounts for the large jumps in the tensor product results where $h$ is fixed. 
The reason for this error was undetermined in our experiments, but will be revisited when the L-shaped results are discussed. 
We see the same behavior for the Dirichlet problem in Figure~\ref{fig:DirTwoPiP}.
In results not displayed here, we analyze $p$- and $h$- refinements in approximating the non-simple Neumann eigenvalue $5\pi^2$.
The results are qualitatively similar to the previously discussed results.

We also note strange behavior when using elements of order $5$ and $6$. Exhibited in the Neumann case on the square in Figure \ref{fig:NeuTwoPiP}: $p=5$ and  $p=6$, we see that as we refine our mesh further, the error increases. 
The error sometimes increases higher than lower order elements solving the same problem, as seen in many of the plots when $h$ is fixed; the trend in error seems to ``flair up'' towards the end. 
The reason for this behavior is likely due to numerical roundoff errors.

\paragraph*{L-Shaped Domain}
On the L-shaped domain, we see in Figures \ref{fig:LNeuTwoPiK}-\ref{fig:LDirTwoPiP} nearly the same patterns described above when approximating the eigenvalue $2\pi^2$. We note that when $h$ is fixed, the savings achieved by serendipity elements is increased even further. For example, with Neumann boundary conditions and $h = 1/4$ (Figure~\ref{fig:LNeuTwoPiK}), for the p=5 case, both the serendipity and tensor product elements exhibit an error of about $10^{-6}$. The number of degrees of freedom used in the serendipity case however is less than half of that of the tensor case. With Dirichlet boundary conditions as seen in Figure~\ref{fig:LDirTwoPiK}: $h=1/4$, this savings is further increased, with serendipity elements using nearly a third of the degrees of freedom used by tensor product elements.

In addition to the plots described above, we have also added plots depicting the results of approximating the Neumann eigenvalue numerically approximated as $1.4756218450$. Figure~\ref{fig:LNeuOneK} and Figure~\ref{fig:LNeuOneP} show that these results mostly correspond to the previously exhibited behavior with the exception that in Figure \ref{fig:LNeuOneK}, the tensor product elements also achieve better approximations when refining $p$. We also note that in Figure~\ref{fig:LNeuOneP}: $p=4$, the order $4$ tensor product elements have a large decrease in error. This behavior contrasts the increase we saw when approximating $2\pi^2$ with order $4$ tensor product elements over the square and is, again, unexplained.

\begin{figure}[ht]
	\[\includegraphics[width=.7\textwidth]{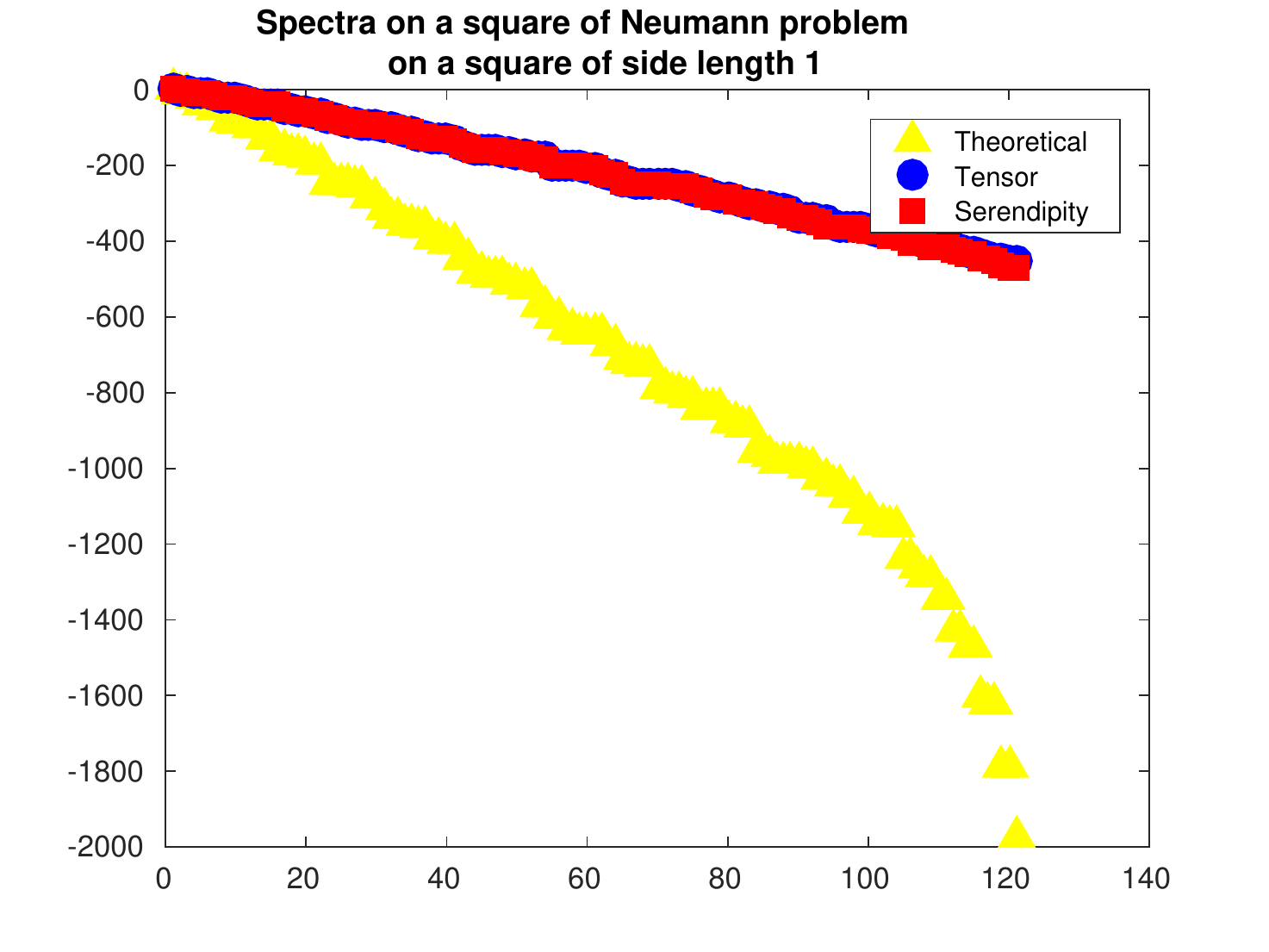}\]
	\caption{Spectra for p = 3, h = 1 / 5, tensor and serendipity bases along with theoretical spectrum over square domain.}
	\label{fig:squareSpectrum}
\end{figure}

\paragraph*{Spectrum comparison}
We also compare the spectrum of eigenvalues that are computed by the tensor product and serendipity elements on the square versus the theoretical spectrum. The results are shown in Figure \ref{fig:squareSpectrum}. We see that the eigenvalues calculated by the tensor product and serendipity elements are nearly the same, and, as expected, as we attempt to approximate larger eigenvalues, the results become less accurate.

%

\section*{Plots from Numerical Experiments}

Plots begin on the next page to aid in comparison between similar cases.

\newpage

\clearpage
\section{Conclusion and Future Directions}
\label{sec:conc}

A key takeaway message from our numerical experiments is that when seeking eigenvalue estimates on a fixed mesh of squares, serendipity elements do appear to fulfill their promise of producing as accurate a result as tensor product elements, despite having roughly 50\% the number of degrees of freedom.
Since many application contexts require a fixed domain mesh, it would be advantageous computationally to use serendipity elements in such circumstances.

Various additional experiments are planned.
First, there are questions in regards to differing behavior on the square versus the L-shaped domain, and the Neumann versus Dirichlet boundary conditions.
A study of serendipity elements for the Poisson equation (i.e.\ with non-zero boundary conditions) or for more general eigenvalue problems might help explain our results.
A second issue is to resolve the dramatic aberrations in the results for the case of tensor product basis functions for the case $p=4$. 
Further investigation into the pattern observed in the convergence behavior depending on mesh discretization is in progress.

We also plan to investigate the observation that mesh discretization for high degree polynomial basis functions sometimes results in less accurate approximations. 
We suspect that this arises from numerical roundoff errors, as the results became worse only after reaching a threshold on the order of $10^{-8}$.

As discussed in the Mathematica code accompanying this paper, similar constructions for serendipity basis functions in three dimensions were also determined.
In future work, we plan to extend the implementation of our finite element solver to allow for three-dimensional domains, and implement these three-dimensional serendipity basis functions in order to produce similar analysis and comparisons as those that we have found for two dimensions.

\begin{appendix}
\section{Serendipity Basis Functions}
\label{app:mathematica}

The following are the serendipity element basis functions in two-dimensions from order $1$ to $4$.
The basis functions are organized as they are calculated in Mathematica, i.e.\ as the sum of re-indexed arrays of basis functions as discussed Section~\ref{sec:implement}.
The Mathematica code that was used to generate these functions is available at
\url{http://math.arizona.edu/~agillette/research/srdpBasisFns/}

\clearpage
$\calS_1\Lambda^0(\square_2)~\text{basis}~ = $\\

\makebox[\textwidth][c]{\scalebox{1}{
$ \left(
\begin{array}{cc}
 \frac{1}{4} (1-x) (1-y) & \frac{1}{4} (1-x) (y+1) \\
 \frac{1}{4} (x+1) (1-y) & \frac{1}{4} (x+1) (y+1) \\
\end{array}
\right) $}}\\[3ex]

$\calS_2\Lambda^0(\square_2)~\text{basis}~ = $\\

\makebox[\textwidth][c]{\scalebox{.9}{
$ \left(
\begin{array}{ccc}
 -\frac{1}{4} (x-1) (y-1) (x+y+1) & \frac{1}{2} (x-1) \left(y^2-1\right) & \frac{1}{4} (x-1) (x-y+1) (y+1) \\
 \frac{1}{2} \left(x^2-1\right) (y-1) & 0 & -\frac{1}{2} \left(x^2-1\right) (y+1) \\
 \frac{1}{4} (y-1) \left(-x^2+y x+y+1\right) & -\frac{1}{2} (x+1) \left(y^2-1\right) & \frac{1}{4} (x+1) (y+1) (x+y-1) \\
\end{array}
\right) $}}\\[3ex]

$\calS_3\Lambda^0(\square_2)~\text{basis}~ = $\\

\makebox[\textwidth][c]{\scalebox{.7}{
$ \left(
\begin{array}{cccc}
 \frac{1}{4} (x-1) (y-1) \left(x^2+y^2-1\right) & \frac{1}{2} (x-1) \left(y^2-1\right) & \frac{1}{2} (x-1) y \left(y^2-1\right) & -\frac{1}{4} (x-1) (y+1) \left(x^2+y^2-1\right) \\
 \frac{1}{2} \left(x^2-1\right) (y-1) & 0 & 0 & -\frac{1}{2} \left(x^2-1\right) (y+1) \\
 \frac{1}{2} x \left(x^2-1\right) (y-1) & 0 & 0 & \frac{1}{2} \left(x-x^3\right) (y+1) \\
 -\frac{1}{4} (x+1) (y-1) \left(x^2+y^2-1\right) & -\frac{1}{2} (x+1) \left(y^2-1\right) & \frac{1}{2} (x+1) \left(y-y^3\right) & \frac{1}{4} (x+1) (y+1) \left(x^2+y^2-1\right) \\
\end{array}
\right)$}}\\[3ex]

$\calS_4\Lambda^0(\square_2)~\text{basis}~ = $\\

\makebox[\textwidth][c]{\scalebox{.45}{
$ \left( \begin{array}{ccccc}
 -\frac{1}{4} (x-1) (y-1) \left(x^3-(y+1) x+y \left(y^2-1\right)\right) & \frac{1}{2} \left(y^2-1\right) \left(-x^2+y^2 x+x-y^2\right) & \frac{1}{2} (x-1) y \left(y^2-1\right) & \frac{1}{4} (x-1) (y-1) y^2 (y+1) & \frac{1}{4} (x-1) (y+1) \left(x^3+(y-1) x-y^3+y\right) \\
 \frac{1}{2} \left(x^2-1\right) \left(x^2-y\right) (y-1) & \left(x^2-1\right) \left(y^2-1\right) & 0 & 0 & -\frac{1}{2} \left(x^2-1\right) (y+1) \left(x^2+y\right) \\
 \frac{1}{2} x \left(x^2-1\right) (y-1) & 0 & 0 & 0 & \frac{1}{2} \left(x-x^3\right) (y+1) \\
 \frac{1}{4} (x-1) x^2 (x+1) (y-1) & 0 & 0 & 0 & -\frac{1}{4} (x-1) x^2 (x+1) (y+1) \\
 \frac{1}{4} (x+1) (y-1) \left(-x^3+y x+x+y^3-y\right) & -\frac{1}{2} \left(y^2-1\right) \left(x^2+y^2 x+x+y^2\right) & \frac{1}{2} (x+1) \left(y-y^3\right) & -\frac{1}{4} (x+1) (y-1) y^2 (y+1) & \frac{1}{4} (x+1) (y+1) \left(x^3+(y-1) x+y \left(y^2-1\right)\right) \\
\end{array} \right) $ } }
\end{appendix}

\text{}\\
\textbf{Acknowledgements.}
All authors were supported in part by NSF Award 1522289.
\vspace{-3mm}

\bibliography{../bib-srdp-eig}{}
\bibliographystyle{abbrv}

\end{document}